\documentclass[10pt, a4paper]{article}

\usepackage{amssymb}
\usepackage{amsthm}
\usepackage{amsmath}
 
\allowdisplaybreaks
    
\newcommand{\R}{\mathbf{R}} 
\newcommand{\C}{\mathbf{C}}

\newcommand{\eps}{\varepsilon} 
\newcommand{\Rtn}{\mathbf{R}^{2n}}
\newcommand{\me}{\mathrm{e}} 
\newcommand{\mi}{\mathrm{i}}
\newcommand{\dif}{\mathrm{d}} 
\newcommand{\Dif}{\mathrm{D}}

\newcommand{\imdz}{\mathbf{Im} \, \dif z }
\newcommand{\redz}{\mathbf{Re} \, \dif z}
\newcommand{\elbow}{\, \rfloor \, }
\newcommand{\Mapprox}{\bar{M}_{\alpha}}
\newcommand{\Mapproxj}{\bar{M}_{\alpha_{j}}}
\newcommand{\Sapprox}{\bar{S}_{\alpha}}
\newcommand{\Sapproxj}{\bar{S}_{\alpha_{j}}}
\newcommand{\ctbr}{C^{2,\beta}_{\rho}(\Mapprox)}
\newcommand{\cobr}{C^{0,\beta}_{\rho}(\Mapprox)}
\newcommand{\psizero}{\psi_{0,\alpha}}
\newcommand{\psione}{\psi_{1,\alpha}}
\newcommand{\Bonealpha}{\mathcal{B}_{1,\alpha}}
\renewcommand{\Re}{\mathbf{Re}} 
\renewcommand{\Im}{\mathbf{Im}}

\newcommand{\happrox}{\bar{h}_{\alpha}}
\renewcommand{\det}{\mathrm{det}}
\newcommand{\thapprox}{\theta_{\Mapprox}}
\newcommand{\meancurv}{\vec{H}_{\Mapprox}}
\newcommand{\Bonealphaj}{\mathcal{B}_{1, \alpha_{j}}}    
\newcommand{\ctbrj}{C^{2,\beta}_{\rho_{j}}(\Mapproxj)}
\newcommand{\cobrj}{C^{0,\beta}_{\rho_{j}}(\Mapproxj)}
\newcommand{\psizeroj}{\psi_{0,\alpha_{j}}}
\newcommand{\psionej}{\psi_{1,\alpha_{j}}}
\newcommand{\meancurvj}{\vec{H}_{\Mapproxj}}

\newcommand{\Cnuone}{C_{1}} % Bound on the first eigenfunction
\newcommand{\Casym}{C_{0}}     % Asymptotics of Lawlor neck
\newcommand{\Rasym}{R_{0}}     % Asymptotics of Lawlor neck
\newcommand{\Ceps}{C}      % Relation b/w \eps and \alpha
\newcommand{\Rinit}{A_{1}}     % This initial size of \alpha
\newcommand{\Cimdz}{C} % Size of Im dz
\newcommand{\Rimdz}{A_{2}}     % Size of required \alpha
\newcommand{\Cvoltr}{C}    % Size of vol of neck + trans region
\newcommand{\Rvol}{A_{2}}      % Size of \alpha to get vol bound
\newcommand{\Rweight}{R} % Radius in which metric is uniform
\newcommand{\Cell}{C_{Ell}}     % Elliptic estimate constant
\newcommand{\Rms}{A_{3}}       % Size of \alpha for MS inequality
\newcommand{\Cms}{E_{MS}}          % M-S constant
\newcommand{\Cnutwo}{C_{2}}   % 2nd eigenvalue bound
\newcommand{\Rnew}{A_{4}}     % A number less than them all so far.
\newcommand{\Cczero}{C} % C0 injectivity bound
\newcommand{\Crho}{K_{\eps}} % Bound on the derivative of the weight
\newcommand{\Cinj}{C}     % C^{k,\beta}_{\rho} injectivity bound
\newcommand{\Crest}{C}   % Constant involving lower order terms
\newcommand{\Rrest}{A_{5}}   % Size of \alpha for the lower order terms
\newcommand{\Capprox}{C}  % Approx of first eigenfunction
\newcommand{\Rinj}{A_{6}}      % Size of \alpha for injectivity
\newcommand{\Rfull}{A_{7}}     % Size of \alpha for full injectivity
\newcommand{\Ctrnsvsl}{C} % Lower bound for <\psione, S_{\alpha}>
\newcommand{\Rsurj}{A_{7}} % Size of surjectivity
\newcommand{\CsizeE}{C} % Size of C^{0,\beta}_{\rho} norm of E

\newtheorem*{mainthm}{Main Theorem}
\newtheorem{thm}{Theorem}
\newtheorem{lemma}[thm]{Lemma}
\newtheorem{cor}[thm]{Corollary}
\newtheorem{prop}[thm]{Proposition}

\theoremstyle{definition}
\newtheorem{defn}[thm]{Definition}

\begin{document}
    
\title{Regularizing a Singular Special Lagrangian Variety}

\author{Adrian Butscher \\ Max Planck Institute for Gravitational
Physics \\ email: \ttfamily butscher@aei-potsdam.mpg.de}

\maketitle

\begin{abstract}
    Suppose $M_{1}$ and $M_{2}$ are two special Lagrangian
    submanifolds with boundary of $\Rtn$, $n \geq 3$, that intersect
    transversally at one point $p$.  The set $M_{1} \cup M_{2}$ is a
    singular special Lagrangian variety with an isolated singularity
    at the point of intersection.  Suppose further that the tangent
    planes at the intersection satisfy an angle criterion (which
    always holds in dimension $n=3$).  Then, $M_{1} \cup M_{2}$ is
    regularizable; in other words, there exists a family of smooth,
    minimal Lagrangian submanifolds $M_{\alpha}$ with boundary that
    converges to $M_{1} \cup M_{2}$ in a suitable topology.  This
    result is obtained by first gluing a smooth neck into a
    neighbourhood of $M_{1} \cap M_{2}$ and then by perturbing this
    approximate solution until it becomes minimal and Lagrangian.
\end{abstract}

\tableofcontents

% \renewcommand{\baselinestretch}{1.25}
% \normalsize

\section{Introduction}

A minimal Lagrangian submanifold of a symplectic manifold $\Sigma$ is
at once minimal with respect to the metric of $\Sigma$ and Lagrangian
with respect to the symplectic structure of $\Sigma$.  Furthermore,
when $\Sigma$ is a Calabi-Yau manifold, Harvey and Lawson showed in
their seminal paper \cite{hl1} that minimal Lagrangian submanifolds
are also \emph{calibrated} by the real part, up to phase, of the
canonical, non-vanishing holomorphic $(n,0)$-form $\zeta$ of
$\Sigma$.  (I.~e.  $\zeta$ is bounded above by the $n$-volume form
and equality holds, up to a phase angle, only on the tangent spaces of
minimal Lagrangian submanifolds.)  A consequence of this property is
that minimal Lagrangian submanifolds satisfy a relatively simple
geometric PDE (simple relative to the equations of vanishing mean
curvature, which they would satisfy by virtue of minimality alone). 
Namely, $M \subset \Sigma$ is minimal Lagrangian if and only if
\begin{equation}
    \begin{gathered}
	\left.\Im \left( \me^{\mi \theta} \zeta \right) \right|_{M} =
	0 \\
	\left. \omega \right|_{M} = 0 \, ,
    \end{gathered}
    \label{eqn:minlageqn}
\end{equation}
for some real number $\theta$.  Here, $\omega$ is the symplectic form
of $\Sigma$.  The calibration form of $M$ is in this case $\Re \left(
\me^{\mi \theta} \zeta \right)$.  The term \emph{special} Lagrangian
refers to those $M$ whose \emph{calibrating angle} $\theta$ vanishes.

Many researchers have exploited the geometric structure implicit in
the calibration condition in order to tackle questions related to the
existence and properties of minimal Lagrangian submanifolds.  Harvey
and Lawson themselves produced several examples of minimal Lagrangian
submanifolds and gave certain general constructions of such objects. 
More recently, Schoen and Wolfson \cite{sandw} have been working on
questions of existence of special Lagrangian submanifolds using
variational techniques.  Current developments in the mathematical
foundations of string theory, in the form of mirror symmetry and the
Strominger-Yau-Zaslow conjecture \cite{morrison}, have greatly
stimulated further investigation into special Lagrangian submanifolds
in the hopes of understanding the moduli spaces of these objects. 
McLean \cite{mclean} began the study of deformations of smooth special
Lagrangian submanifolds and Hitchin \cite{hitchin} deduced certain
properties of the moduli space of such deformations.  Compactification
of this moduli space (an important ingredient in mirror symmetry),
however, requires understanding the singularities that can arise
amongst special Lagrangian submanifolds.  Harvey and Lawson
contributed a number of examples of singular special Lagrangian
varieties, while Haskins \cite{haskins} extended these results by
constructing new examples of special Lagrangian cones and their
desingularizations.  Many further advances on the topic of singular
special Lagrangian submanifolds are presented in papers by Joyce (for
a summary, see \cite{joyce1}).

Despite these advances, still relatively little is known about
singularities of special Lagrangian submanifolds and research in this
domain continues.  This paper investigates singular Lagrangian
geometry in a setting \emph{related} to the one in which mirror
symmetry takes place.  That is, this paper studies the nature of the
singular special Lagrangian submanifold \emph{with boundary} of
$\Rtn$, $n \geq 3$, formed by the union of two otherwise smooth
special Lagrangian submanifolds $M_{1}$ and $M_{2}$ that intersect
transversally in one point $p$, and whose tangent cone at $p$
satisfies an \emph{angle criterion} to be explained shortly.  The
union $M$ of the two smooth submanifolds is singular at the
intersection point (which is isolated for reasons of transversality)
and the tangent cone there is just the union of the two tangent planes
of the constituent smooth submanifolds.  The question this paper
partially answers in the affirmative in the Main Theorem, to be stated
precisely below, is whether $M$ can be realized as the limit of a
sequence of smooth minimal Lagrangian submanifolds boundary, thus
placing $M$ in the compactification of the moduli space of minimal
Lagrangian submanifolds with boundary of $\Rtn$.  It is hoped that
this result as well as the methods used in its proof will shed some
light on the analogous question of intersections of special Lagrangian
submanifolds in compact Calabi-Yau manifolds.

Before stating precisely the result that will be proved in
this paper, it is necessary to mention that the boundary behaviour of
the regularizing family will play a critical role in the forthcoming
analysis.  The way in which the boundary of $M$ will be controlled is
to use a submanifold of $\Rtn$ that will be called a \emph{scaffold}.

\begin{defn} \label{defn:scaffold}
    Let $M$ be a submanifold of $\Rtn$ with boundary $\partial M$ and
    inward unit normal vector field $Z \in \Gamma \big( T_{\partial
    M}M \big)$.  A \emph{scaffold} for $M$ is a smooth submanifold $W$
    of $\Rtn$ with the following properties:
    \begin{enumerate}
	\item $\partial M \subset W$; 
	
	\item $Z \in \Gamma \big( T_{\partial M} W \big)^{\omega}$
	($V^{\omega}$ is the symplectic orthogonal complement of $V$);
	
	\item The normal bundle of $W$ is trivial.
    \end{enumerate}
\end{defn}

\noindent The idea of the scaffold is first introduced in \cite{me2}
and \cite{me} in the context of a different problem concerning
deformations of smooth minimal Lagrangian submanifolds with boundary
and is discussed in detail there.  Condition (2) might seem slightly
unnatural, but it considerably simplifies the proof of the Main
Theorem.  Nevertheless, one expects the Main Theorem to hold under a
more general transversality condition.

The precise statement of the theorem to be proved in this paper that
incorporates the boundary behaviour of the regularization is the
following.

\begin{mainthm}
    Suppose $M_{1}$ and $M_{2}$ are two special Lagrangian
    submanifolds with boundary of $\Rtn$, $n \geq 3$, that intersect
    transversally at one point $p$.  Furthermore, suppose that the
    tangent planes of $M_{1}$ and $M_{2}$ satisfy the \emph{angle
    criterion} (described below) at the point $p$ and let $W$ be a
    scaffold for $M$ that is also a codimension 2, symplectic
    submanifold of $\Rtn$.  Define $M = M_{1} \cup M_{2}$.  Then there
    exists a family $M_{\alpha}$ of smooth, minimal Lagrangian
    submanifolds with boundary and a family of symplectic,
    codimension-two submanifolds $W_{\alpha}$ such that the following
    results hold.
    \begin{enumerate}
	\item $M_{\alpha} \rightarrow M$ in some suitable topology;
	\item $W_{\alpha} \rightarrow W$ in the same topology;
	\item $\partial M_{\alpha} \subseteq W_{\alpha}$ for every
	$\alpha$.
   \end{enumerate}    
\end{mainthm}

\medskip \noindent \scshape Remarks: \upshape (1) If $M$ is any
Lagrangian submanifold of $\Rtn$ with boundary then there
\emph{always} exists a symplectic, codimension-two scaffold for $M$:
for instance, it can be constructed by exponentiating a sufficiently
small neighbourhood $\mathcal{U}$ of the zero section of the bundle
$T(\partial M) \oplus J T (\partial M)$.  (2) Due to technical reasons
that make their appearance only in Section 6, the argument used to
prove the Main Theorem fails in dimension $n=2$.  It is unfortunately
not yet clear to the Author how to proceed in the $n=2$ case using the
techniques developed herein; however, the $n=2$ case is true in
certain circumstances for other reasons.  See the end of this section
for an explanation.\medskip

The angle criterion that $M_{1}$ and $M_{2}$ must satisfy can be
explained as follows.  It can be shown that if $P_{1}$ and $P_{2}$ are
any two transversally intersecting Lagrangian planes in $\Rtn$, then,
up to a global unitary transformation of $\Rtn$, there exists an
orthonormal Darboux basis $E_{1}, \ldots, E_{n}$ and $JE_{1}, \ldots,
JE_{n}$ of $\Rtn$ and a unique set of angles $\theta_{i}$, where
\begin{equation}
    \theta_{1} \in [\frac{\pi}{2}, \pi) \qquad \mbox{and} \qquad
    \theta_{i} \in [0, \frac{\pi}{2} ) \quad \mbox{for $i = 2, \ldots, 
    n$} \, , 
    \label{eqn:range}
\end{equation}
such that
\begin{gather*}
    P_{1} = \mathrm{span} \big\{ E_{1}, \ldots, E_{n} \big\} \\
    \intertext{and}
    P_{2} = \mathrm{span} \big\{ \cos \theta_{1} E_{1} + \sin 
    \theta_{1} JE_{1}, \ldots, \cos \theta_{n} E_{n} + \sin \theta_{n} 
    JE_{n} \big\} \, .
\end{gather*}
The existence of this basis follows from standard symplectic linear
algebra and can be found in \cite{me,ynging} as well as \cite{harvey},
for example.  Since an orthonormal Darboux basis is by definition
holomorphic, the form $\dif z$ becomes $(E_{1}^{\ast} + \mi
JE_{1}^{\ast} ) \wedge \cdots \wedge (E_{n}^{\ast} + \mi JE_{n}^{\ast}
)$ in these coordinates, where $E_{i}^{\ast}$ denotes the dual 1-form
corresponding to the vector $E_{i}$.  Thus by applying the calibrating
form $\imdz$ to the pair of planes $P_{1}$ and $P_{2}$, one sees that
they are special Lagrangian if and only if $\theta_{1} + \cdots +
\theta_{n} = m \pi$ for some positive integer $m$.  The \emph{angle
criterion} for $P_{1}$ and $P_{2}$ is that the characteristic angles
satisfy $\theta_{1} + \cdots + \theta_{n} = \pi$.  Note that the angle
criterion need not hold for a general pair of intersecting special
Lagrangian planes, though it is \emph{always} satisfied in dimensions
$n=2$ and $3$ for numerical reasons (i.~e.\ as a result of
\eqref{eqn:range}).  The angle criterion that the pair of submanifolds
$M_{1} \cup M_{2}$ must satisfy in order for the Main Theorem to hold
is that the characteristic angles of the pair of special Lagrangian
planes that form the tangent cone $T_{p}M_{1} \cup T_{p} M_{2}$ at the
singularity $p$ satisfy the angle criterion explained above.

The Author would like to thank Yng-Ing Lee for pointing out the
necessity and importance of the angle criterion in the proof of the
Main Theorem.  The angle criterion is a necessary requirement on
$M_{1}$ and $M_{2}$ because only when it holds, is it known that there
exists a \emph{local regularization} for the tangent cone $T_{p}M_{1}
\cup T_{p} M_{2}$ --- that is, a smooth submanifold $N$ with two ends,
asymptotic to $T_{p} M_{1}$ and $T_{p} M_{2}$ respectively.  This
object will be called a \emph{Lawlor neck}, and was discovered by Lawlor
in \cite{lawlor2}, but see also \cite{harvey} for a more thorough
treatment.  The first step of the proof of the Main Theorem is to glue
a rescaled piece of this Lawlor neck into the neighbourhood of $p$ in
order to obtain a smooth, approximately minimal Lagrangian
submanifold.  Thus without the angle criterion, this first step can
not be undertaken and the proof of the Main Theorem does not get off
the ground.  When the angle criterion fails for $M_{1}$ and $M_{2}$,
one will first have to find a new local regularization of the tangent
cone at the singularity (one which is necessarily not of Lawlor's
type, as indicated by Lee) before being able to apply the techniques
of this paper to obtain a version of the Main Theorem valid in this
case.  The reader should consult Lee's paper \cite{ynging} for a more
complete discussion of this issue, and how it relates to her work on
immersed, self-intersecting special Lagrangian submanifolds without
boundary in general Calabi-Yau manifolds.

\medskip \noindent \scshape Remark on the $n=2$ case: \upshape Suppose
$M_{1}$ and $M_{2}$ are two special Lagrangian submanifolds with
boundary in $\R^{4}$ that intersect transversally at one point.  As
will become clear in Section 6.3, the techniques used to prove the
Main Theorem fail in this situation.  However, in certain cases,
another method provides a solution.  Let $x^{1}, x^{2}, y^{1},
y^{2}$ denote the standard coordinates of $\R^{4}$.  It is easy to
verify that the coordinate transformation
\begin{equation*}
    (x^{1}, x^{2}, y^{1}, y^{2}) \mapsto (x^{1}, y^{1}, -x^{2}, y^{2})
\end{equation*}
applied to $M_{1}$ and $M_{2}$ produces submanifolds $M_{1}'$ and
$M_{2}'$ intersecting transversally at one point whose tangent spaces
are invariant under the complex structure of $\R^{4}$.  In other
words, $M_{1}'$ and $M_{2}'$ are intersecting Riemann surfaces.  Thus
the question of regularizing $M_{1} \cup M_{2}$ becomes a problem in
two complex variables.  If $M_{1}'$ and $M_{2}'$ can be represented as
the zero sets of holomorphic functions $h_{1}$ and $h_{2}$, then the
variety $M' = M_{1}' \cup M_{2}'$ can be regularized using complex
analytic techniques.  Indeed, $M'$ can be represented as the zero set
of the holomorphic function $h = h_{1} \cdot h_{2}$, and due to the
nature of the singularity (without loss of generality, occuring at the
origin) of $M'$, there are coordinates $(z,w)$ for a neighbourhood of
the origin in which $h(z,w) = z \cdot w \cdot \tilde{h}$, where
$\tilde{h}(0,0) \neq 0$.  Consequently, the function $h_{\eps} = h +
\eps$ is holomorphic and has smooth zero section, provided $\eps$ is
sufficiently small.  The Riemann surfaces $M'_{\eps} =
h_{\eps}^{-1}(0)$ thus regularize $M'$.  Changing coordinates back to
the original ones then provides the desired regularization of $M$. 
This method is not always available since $M_{1}'$ and $M_{2}'$ may
not be representable \emph{globally} as the zero sets of holomorphic
functions.  In this case, new methods will have to be developed to
solve the regularization problem, and these are, at the moment, beyond
the scope of this paper.

\section{Preliminaries}

\subsection{Outline of the Proof}

The regularization $M_{\alpha}$ of the singular variety $M = M_{1}
\cup M_{2}$ will be constructed by applying gluing techniques
\cite{joyce, nikos2, nikos1, rafe} in combination with the Inverse
Function Theorem to a system of differential equations associated to
the geometric equations \eqref{eqn:minlageqn} whose solutions
correspond to minimal Lagrangian submanifolds near $M$.  The purpose
of this section is to describe this procedure in general terms.

The preliminary step in the proof of the Main Theorem is to replace
$M$ by a family of smooth, embedded, Lagrangian submanifolds
$\Mapprox$ \emph{approximating} $M$, i.~e.\ by submanifolds $\Mapprox$
that are \emph{almost} special Lagrangian in a precise sense and that
converge to $M$ in a suitable topology.  These approximating
submanifolds will be constructed in detail in Section 3.  The reason
for replacing $M$ by $\Mapprox$ is that the point of non-smoothness of
$M$ makes it difficult to apply PDE techniques directly to $M$.  The
question of finding a minimal Lagrangian submanifold near $M$ can be
converted to solving a PDE on $\Mapprox$ as follows.

Denote by $\happrox : \Mapprox \rightarrow \Rtn$ the canonical
embedding of $\Mapprox$ and let $\mathcal{B}$ be a Banach space that
parametrizes a subset of $C^{1,\beta}$ Lagrangian embeddings of
$\Mapprox$ into $\Rtn$ that are near $\happrox$ in some suitable
topology.  In other words, suppose that there is a continuous map
$\Phi : \mathcal{B} \rightarrow C^{1,\beta} \big( \mathit{Emb}
(\Mapprox, \Rtn) \big)$ into the immersions of $\Mapprox$ into $\Rtn$
with $\Phi(0) = \happrox$ and such that $\Phi(x) (\Mapprox)$ is
Lagrangian for every $x \in \mathcal{B}$.  The parametrization and the
Banach space will be specified in Section 4 and will be such that if
$\Vert x \Vert_{\mathcal{B}}$ is sufficiently small, then $\Phi(x)
(\Mapprox)$ is smooth as well as embedded.  The first of the two
minimal Lagrangian equations \eqref{eqn:minlageqn}, namely that
$\Phi(x) ^{\ast} \omega = 0$, is automatically satisfied because
each $\Phi(x)$ is a Lagrangian embedding, by definition.  Hence it is
sufficient to analyze only the second of the two equations.  Recall
that the canonical, non-vanishing holomorphic $(n,0)$-form of $\Rtn$
is $\zeta = \dif z^{1} \wedge \cdots \wedge \dif z^{n}$ in standard
holomorphic coordinates.  This form will henceforth be abbreviated by
$\dif z$.  Consider now the differential operator between Banach
spaces defined by the second equation of \eqref{eqn:minlageqn}:
namely, the map
$$F_{\alpha} : \mathcal{B} \times \R \longrightarrow
C^{0,\beta}(\Mapprox)$$
given by
\begin{equation}
    F_{\alpha}(x, \theta) = \left\langle \Phi(x)^{\ast} \big( \Im
    (\me^{\mi \theta} \dif z) \big), \mathrm{Vol}_{\Mapprox}
    \right\rangle_{\Mapprox} 
    \label{eqn:minlagmap} 
\end{equation}
for any $x \in \mathcal{B}$.  Here, $\mathrm{Vol}_{\Mapprox}$ is the
induced volume form of $\Mapprox$ and $\langle \cdot, \cdot
\rangle_{\Mapprox}$ is the induced inner product on the $n$-forms of
$\Mapprox$.

A solution of the equation $F_{\alpha} (x, \theta) = 0$ for some $(x,
\theta)$ near the origin in $\mathcal{B} \times \R$ corresponds to an
embedded minimal Lagrangian submanifold $M_{\alpha} = \Phi(x)
(\Mapprox) \subset \Rtn$ with calibration angle $\theta$ that is near
$\Mapprox$, and thus near $M$ as well.  Solutions are found by
invoking the following version of the Inverse Function Theorem (see
\cite{amr} for the proof).

\begin{thm}[Inverse Function Theorem] \label{thm:ift}
    Let $F : \mathcal{B} \longrightarrow \mathcal{B}'$ be a $C^{1}$
    map between Banach spaces and suppose that the linearization $\Dif
    F(0)$ of $F$ at $0$ is an isomorphism.  Moreover, suppose $F$
    satisfies the estimates:
    \begin{enumerate}
	\item $\Vert \Dif F(0) x \Vert_{\mathcal{B}'} \geq C_{L} \Vert
	x \Vert_{\mathcal{B}}$ for any $x \in \mathcal{B}$,
	\item $\Vert \Dif F(0) y - \Dif F(x) y \Vert_{\mathcal{B}'}
	\leq C_{N} \Vert x \Vert_{\mathcal{B}} \, \Vert y
	\Vert_{\mathcal{B}}$ for all $x$ sufficiently near $0$ and
	for any $y \in \mathcal{B}$,
    \end{enumerate}
    where $C_{L}$ and $C_{N}$ are constants independent of $x$ and
    $y$.  Let $r \leq C_{L} / 2 C_{N}$.  Then there exist neighbourhoods
    $\mathcal{U}$ of $0$ and $\mathcal{V}$ of $F(0)$ so that
    $F:\mathcal{U} \longrightarrow \mathcal{V}$ is a $C^{1}$
    diffeomorphism and $\mathcal{V}$ contains the ball $B_{r C_{L}
    /2}(F(0))$.  Furthermore, $B_{r C_{L} /2}(F(0))$ is in the image
    of the ball $B_{r}(0)$ under $F$.  
\end{thm}

\noindent The submanifold $\Mapprox$ only approximates $M$ and is not
necessarily minimal Lagrangian; thus $F_{\alpha} (0, 0)$ is \emph{not}
equal to $0$.  But according to the theorem above, the equation
$F_{\alpha} (x, \theta) = 0$ can be solved for $(x, \theta)$ near
$(0,0)$ if $0$ belongs to the ball $B_{r C_{L}/2}(F(0))$; that is, if
\begin{equation}
\Vert F_{\alpha} (0, 0) \Vert_{C^{0,\beta}(\Mapprox)} \leq r C_{L} /2 
\, ,
\label{eqn:fundest}
\end{equation}
and the theorem asserts that the solution satisfies $\Vert (x,\theta)
\Vert_{\mathcal{B} \times \R} \leq r$.  

The Main Theorem will be proved by verifying that the inequality
\eqref{eqn:fundest} holds for the choices of approximating submanifold
$\Mapprox$ and parametrization of nearby Lagrangian embeddings made in
Sections 2 and 3.  In Section 4, the constant $C_{L}(\alpha)$
(depending of course on $\alpha$) will be estimated by analyzing the
linearization at $(0,0)$ of the operator $F_{\alpha}$ and in Section
5, the proof will be completed by deriving the constant $C_{N}
(\alpha)$ (also depending on $\alpha$), invoking the Inverse Function
Theorem, and using the estimate $\Vert (x,\theta) \Vert_{ \mathcal{B}
\times \R} \leq r$ to be used to guarantee the smoothness of the
solution.

\subsection{The Linearized Operator}

The linearized operator $\Dif F_{\alpha}(0,0)$ will clearly play a
central role in the forthcoming analysis.  It will be helpful, at this
point, to derive a general expression for this operator in terms of
the parametrization $\Phi$, since its specific form will guide the
course of the proof.  The following calculations are based on McLean's
work on perturbations of smooth, special Lagrangian submanifolds
\cite{mclean}.

Begin with two preliminary observations.  First, recall that
$\Phi(tx)$ is a family of Lagrangian embeddings.  Hence,
$\Phi(tx)^{\ast} \omega = 0$ for every $t$.  Differentiating this
expression in $t$ and using the Lie derivative formula
$\mathcal{L}_{V} \eta = \dif( V \elbow \eta) + V \elbow \dif \eta$ as
in McLean's paper produces the expression
\begin{equation*}
    \dif \, \happrox^{\, \ast} (V_{x} \elbow \omega) = 0 \, ,
\end{equation*}
where $V$ is the
\emph{deformation vector field} associated to $\Phi(tx)$ given by
\begin{equation}
    V_{x} = \left.  \frac{\dif}{\dif t} \Phi(tx) \right|_{t=0} \, .
    \label{eqn:deffield}
\end{equation}
Second, if the parametrization $\Phi$ is chosen such that the
closedness of the form $\happrox^{\, \ast} (V_{x} \elbow \omega)$
implies exactness, then there is a function $H_{x}:\Mapprox
\rightarrow \R$ so that $\happrox^{\, \ast} (V_{x} \elbow \omega) =
\dif H_{x}$.  Refer to $H_{x}$ as the \emph{Hamiltonian function
associated to $x$}.

Next, it is necessary to remind the reader of the \emph{Lagrangian
angle function} that is defined on \emph{any} Lagrangian submanifold
of a C-Y manifold $\Sigma$ with canonical, non-vanishing holomorphic
$(n,0)$-form $\zeta$.  First, let $e_{1}, \ldots, e_{n}$ be an
orthonormal basis for any Lagrangian subspace of $T_{x}\Sigma$.  It
can be shown that the complex number $\zeta(e_{1}, \ldots, e_{n})$ is
independent of the particular basis and always has modulus equal to
one.  See \cite[page 89]{hl1} for details.  Consequently, if $M$ is a
Lagrangian submanifold of $\Sigma$ and $e_{1}, \ldots, e_{n}$ is an
orthonormal basis for $T_{x}M$, then there is an angle $\theta(x)$
that satisfies
\begin{equation}
    \me^{- \mi \theta(x)} \equiv \zeta (e_{1}, \ldots, e_{n}) \,
    .
    \label{eqn:lagangle2}
\end{equation}
If $M$ is minimal Lagrangian and has calibration angle $\theta$ then
$\theta(x) = \theta$ for all $x \in M$.  The significance of the
possibly multi-valued function $\theta$ lies in the fact that its
differential $\dif \theta$ is well-defined and satisfies the relation
\begin{equation}
    \dif \theta = - \vec{H}_{M} \elbow \omega \, ,
    \label{eqn:angleHrel}
\end{equation}
where $\vec{H}_{M}$ is the mean curvature vector of $M$.

The linearization of the operator $F_{\alpha}$ can now be phrased
using this terminology.

\begin{prop} \label{prop:lin}
    Let $H_{x}$ and $V_{x}$ be the Hamiltonian function and the
    deformation vector field associated to $x$ via the parametrization
    $\Phi$ as discussed above.  The linearization of the operator
    $F_{\alpha}: \mathcal{B} \times \R \rightarrow C^{0,\beta}
    (\Mapprox)$ at $(0,0)$ in the direction $(x,a) \in \mathcal{B}
    \times \R$ is given by
    \begin{equation}
	\Dif F_{\alpha}(0,0) (x, a) = - \cos ( \theta_{\Mapprox} ) \,
	\Delta_{\Mapprox} H_{x} - \sin(\theta_{\Mapprox} ) \,
	\big\langle \vec{H}_{\Mapprox} , V_{x} \big\rangle_{\Mapprox}
	+ a \cos(\thapprox)
	\label{eqn:lin}
    \end{equation}
    where $\vec{H}_{\Mapprox}$ and $\theta_{\Mapprox}$ are the mean
    curvature and the Lagrangian angle of $\Mapprox$.
\end{prop}

\begin{proof}
    
The linearization of $F_{\alpha}$ at $(0,0)$ in the direction of $x
\in \mathcal{B}$ is defined as the quantity
\begin{equation*}
    \Dif F_{\alpha}(0,0) x \equiv \left.  \frac{\dif}{\dif t}
    F_{\alpha}(tx,0) \right|_{t=0} \, .
\end{equation*}
Thus the calculations performed in Mclean's paper \cite{mclean},
modified for the case when the Lagrangian angle function is not
identically zero, imply that
\begin{align*}
    \Dif F_{\alpha}(0,0) (x,0) &= \left\langle \left.  \frac{\dif}
    {\dif t} \Phi(tx)^{\ast} \imdz \right|_{t=0} ,
    \mathrm{Vol}_{\Mapprox} \right\rangle_{\Mapprox} \\
    &= \left\langle \dif \star \Big( \cos( \theta_{\Mapprox}) \,
    \happrox^{\, \ast} (V_{x} \elbow \omega ) \, \Big),
    \mathrm{Vol}_{\Mapprox} \right\rangle_{\Mapprox} \\
    &= \left\langle \dif \Big( \cos( \theta_{\Mapprox}) \star \dif
    H_{x} \, \Big), \mathrm{Vol}_{\Mapprox} \right\rangle_{\Mapprox}
    \\
    &= - \cos( \theta_{\Mapprox}) \, \Delta_{\Mapprox} H_{x} + \sin (
    \theta_{\Mapprox}) \, \left\langle \dif \theta_{\Mapprox} \wedge
    \star \dif H_{x} , \mathrm{Vol}_{\Mapprox}
    \right\rangle_{\Mapprox} \\
    &= - \cos( \theta_{\Mapprox}) \, \Delta_{\Mapprox} H_{x} + \sin (
    \theta_{\Mapprox}) \, \left\langle \vec{H}_{\Mapprox} , V_{x}
    \right\rangle_{\Mapprox}
\end{align*}
using the relationship \eqref{eqn:angleHrel}.  Here, $\star$ and
$\Delta_{\Mapprox}$ are the Hodge star operator and the Laplacian of
$\Mapprox$ in the induced metric.  The linearization of $F_{\alpha}$
in the $\theta$-direction can be calculated similarly.  In fact, if $a
\in \R$, then
\begin{align*}
   \Dif F_{\alpha}(0,0) (0,a) &= \left\langle \left. 
   \frac{\dif}{\dif t} \Phi(0)^{\ast} \Im (\me^{\mi t a} \dif z)
   \right|_{t=0} , \mathrm{Vol}_{\Mapprox} \right\rangle_{\Mapprox} \\
   &= a \big\langle \happrox^{\, \ast} \Im (\mi \, \dif z) ,
   \mathrm{Vol}_{\Mapprox} \big\rangle_{\Mapprox} \\
   &= a \big\langle \happrox^{\, \ast} \Re (\dif z) ,
   \mathrm{Vol}_{\Mapprox} \big\rangle_{\Mapprox} \\
   &= a \cos(\thapprox) 
\end{align*}
by definition of the angle function.  Combining the two results above
yields the desired expression for the linearization.
\end{proof}

\subsection{Boundary Conditions}

Since the linearization of $F_{\alpha}$ at $(0,0)$ is equal to the
Laplacian plus a lower order term, it is elliptic only when the
functions in $\mathcal{B}$ satisfy appropriate boundary conditions. 
The purpose of this section is to show how the scaffold introduced in
the statement of the Main Theorem brings this about.

Suppose that the boundary of $M$ lies on a scaffold $W$.  Since the
submanifold $\Mapprox$ will differ from $M$ only in a small
neighbourhood disjoint from the boundary, then $\partial \Mapprox =
\partial M$ and thus $W$ is a scaffold for $\Mapprox$ as well.  Let
$\Phi : \mathcal{B} \rightarrow \C^{1,\beta} \big( \mathit{Emb}
(\Mapprox, \Rtn) \big)$ parametrize a subset of Lagrangian embeddings
near $\happrox$ as above.  This time, however, suppose that each
embedding confines the boundary of $\Mapprox$ to $W$; in other words,
suppose that $\Phi(x)(\partial \Mapprox) \subset W$ for every $x \in
\mathcal{B}$.

Consider now a one-parameter family $\Phi(tx)$ of such embeddings. 
Since $\Phi(tx) (\partial \Mapprox) \subset W$ for every $t$, the
deformation vector field $V_{x}$ from \eqref{eqn:deffield} must be
tangent to $W$, which leads to the following result.

\begin{prop} \label{prop:symplscaffold}
    Let $M$ be a Lagrangian submanifold of $\Rtn$ and let $W$ be a
    scaffold for $M$.  If $\phi^{t} : M \rightarrow \Rtn$ is any
    Hamiltonian deformation of $M$ with Hamiltonian $H$ that satisfies
    $\phi_{H}^{t} (\partial M) \subset W$ for all $t$, then $H$
    satisfies $Z(H) \big|_{\partial M} = 0$ where $Z$ is the inward
    unit normal vector field of $\partial M$.
\end{prop}

\begin{proof}

The deformation vector field $V$ of $\phi$ must be parallel to $W$
along $\partial M$.  But according to the definition of a scaffold, $Z
\in (T_{x}W)^{\omega}$ for every $x \in \partial L$.  Therefore
$\omega(Z, V) \big|_{\partial M} = 0$.  Since $V$ is a Hamiltonian
vector field, this equality is equivalent to $Z(H) \big|_{\partial M}
= 0$.
\end{proof}

\noindent Thus if functions in $\mathcal{B}$ confine $\partial
\Mapprox$ to the scaffold $W$, then they satisfy Neumann boundary
conditions.

\section{The Approximate Solution}

\subsection{The Local Regularization}

The proof of the Main Theorem begins with the explicit construction of
the approximate submanifolds $\Mapprox$.  These will be constructed by
gluing an appropriate interpolation between $M_{1}$ and $M_{2}$ in a
neighbourhood of small radius about the singular point of $M_{1} \cup
M_{2}$.  As mentioned in the Introduction, the interpolating
submanifold that will be used is the special Lagrangian submanifold
asymptotic to the tangent cone of $M$ at the singularity which is
known as the Lawlor neck and exists whenever the planes comprising the
tangent cone satisfy the angle criterion.  The purpose of the present
section is to describe the Lawlor submanifold.  The actual gluing will
be carried out in Section 3.2 and the relevant properties of the
resulting submanifold $\Mapprox$ will be derived in Section 3.3.
 
Without loss of generality, the singularity of $M$ is located at the
origin in $\Rtn$.  Let $P_{1}$ and $P_{2}$ denote the tangent planes
$T_{0}M_{1}$ and $T_{0}M_{2}$ respectively, and denote by $P$ the cone
$P_{1} \cup P_{2}$ --- this is the tangent cone of $M$ at $0$.  The
Lawlor neck is an embedded cylinder of the form $N = \Psi (\R \times
\textbf{S} ^{n-1})$, where $\Psi$ is a special Lagrangian embedding
into $\Rtn$, and the two ends of this embedding, namely $E_{1} \equiv
\Psi \left( (\lambda , \infty) \times \textbf{S}^{n-1} \right)$ and
$E_{2} \equiv \Psi \left( (- \infty, - \lambda) \times
\textbf{S}^{n-1} \right)$ tend towards $P_{1}$ and $P_{2}$,
respectively, in a pointwise sense as $\lambda \rightarrow \infty$
(i.~e.\ $E_{1}$ can be written as a graph over $P_{1}$ outside a large
enough ball, and the graphing function tends to zero as the radius of
the ball increases).  Because each rescaled submanifold $\eps N$ is
still special Lagrangian and asymptotic to $P$, the family of
homotheties $\eps N$ is a special Lagrangian regularization of the
singular variety $P$.  The idea behind the gluing construction of
$\Mapprox$ is to use a sufficiently small rescaling of $N$ as the
interpolating submanifold.

The precise definition of the Lawlor neck proceeds as follows.  First,
let $a_{1}, \ldots, a_{n}$ be positive real numbers and let $P :
\R^{n} \times \R \longrightarrow \R$ be the function given by
\begin{equation}\label{eqn:P}
    P(a, \lambda) \equiv \frac{\left( 1+a_{1} \lambda^{2} \right)
    \cdots \left( 1+a_{n} \lambda^{2} \right) - 1}{\lambda^{2}} \, .
\end{equation}
Next, set
\begin{equation}\label{eqn:angle}
    \theta_{k}(a, \lambda) \equiv \int_{0}^{\lambda} \frac{- \dif s}
    {(\frac{1}{a_{k}} + s^{2}) \sqrt{P(a, s)}}
\end{equation}
where $a = (a_{1}, \ldots, a_{n})$.  It is easy to see that the
integrals \eqref{eqn:angle} converge as $\lambda \rightarrow \pm
\infty$.  Let $\theta_{k}(a)$ denote the asymptotic values
$\lim_{\lambda \rightarrow \infty} \theta_{k}( a, \lambda)$; then
$\lim_{\lambda \rightarrow - \infty} \theta_{k}( a, \lambda) =
-\theta_{k}(a)$.  This terminology sets the stage for the definition
of the Lawlor embeddings.

% Definition of the Lawlor embeddings
\begin{subequations}

\begin{defn}\label{defn:lawlor}
    For every $a \in \R^{n}$ with $a_{k} > 0$ for all $k$, the map
    $\Psi_{a} : \R \times \mathbf{S}^{n-1} \longrightarrow \R^{2n}$
    defines a Lawlor neck $N_{a} = \Psi_{a}(\R \times \mathbf{S}
    ^{n-1})$ according to the following prescription.  Let
    \begin{equation}
	\Psi_{a} (\lambda, \mu^{1}, \ldots, \mu^{n}) \equiv \left(
	x^{1}(\lambda, \mu), \ldots , x^{n}(\lambda, \mu); \;
	y^{1}(\lambda, \mu), \ldots , y^{n}(\lambda, \mu) \right)
    \end{equation}
    where 
    \begin{equation}\label{eqn:lawlor}
	\begin{split}
	    &x^{k}(\lambda, \mu) = \mu^{k} \, \sqrt{ \tfrac{1}{a_{k}}
	    + \lambda^{2} } \, \cos \!  \left( \tfrac{\pi}{2}
	    \delta_{1k} + \theta_{k} (a, \lambda) \right) \\
	    &y^{k}(\lambda, \mu) = \mu^{k} \, \sqrt{ \tfrac{1}{a_{k}}
	    + \lambda^{2} } \, \sin \!  \left( \tfrac{\pi}{2}
	    \delta_{1k} + \theta_{k}(a, \lambda) \right) \, .
	\end{split}
    \end{equation}
    and $\mu = (\mu^{1}, \ldots, \mu^{n}) \in \R^{n}$ satisfies $\sum
    (\mu^{k})^{2} = 1$ and thus $\mu$ represents a point in
    $\mathbf{S} ^{n-1}$.  Here, $\delta_{1k}$ is the Kronecker symbol,
    defined to equal zero unless $k=1$, in which case it equals 1.
\end{defn}	
\end{subequations}

It will be necessary to have precise numerical estimates of the degree
of closeness between $N_{a}$ and $P_{1} \cup P_{2}$, but suitable
coordinates must be found to perform the calculations.  It is true
that sufficiently far from the origin, the nearest point projection of
$N_{a}$ onto one or the other of its asymptotic planes will be a
diffeomorphism.  Consequently, coordinates can be chosen so that the
ends of $N_{a}$ are graphs over the corresponding asymptotic planes. 
The desired estimates on the asymptotics of the Lawlor neck $N_{a}$
will be phrased in terms of these graphing functions.

% Theorem: asymptotics of the Lawlor necks

\begin{thm} \label{thm:lawlorasym}
    Suppose $n \geq 3$.  There exists a positive, real number $\Rasym$
    so that $N_{a} \cap \big(B_{\Rasym}(0) \big)^{c}$ consists of two
    ends $E_{1}$ and $E_{2}$ that are graphs over $P_{1} \cap
    \big(B_{\Rasym}(0) \big)^{c}$ and $P_{2} \cap \big(B_{\Rasym}(0)
    \big)^{c}$, respectively, of the gradient of a single function $g
    : P_{i} \cap \big(B_{\Rasym}(0) \big)^{c} \longrightarrow
    P_{i}^{\perp}$.  Furthermore, the function $g$ has the property
    that there exists some constant $\Casym$ depending only on $a_{1},
    \ldots, a_{n}$ and $n$ so that
    \begin{equation*}
	\big\Vert \nabla g (p) \big\Vert + \Vert p \Vert \, 
	\big\Vert \mathrm{Hess} \, g(p) \big\Vert + \Vert p \Vert^{2}
	\big\Vert \nabla^{3} g (p) \big\Vert \leq
	\frac{\Casym}{\hspace{1ex} \left\Vert p \right\Vert^{n-1}} 	
    \end{equation*}
    for any $p \in P_{i}$ with $\Vert p \Vert \geq \Rasym$; and the
    $\beta$-H\"older coefficient of $\nabla^{3}g$ satisfies
    \begin{equation*}
	\big[ \nabla^{3} g \big]_{\beta, (B_{R}(0))^{c}} \leq \frac{C_{0}} 
	{R^{n+1 + \beta}}
    \end{equation*}
    for any radius $R \geq \Rasym$; and finally, $g$ can be chosen so
    that
    \begin{equation*}
	\vert g (p) \vert \leq \frac{\Casym}{\hspace{1ex} \left\Vert p
	\right\Vert^{n-2}}
    \end{equation*}    
    for any $p \in P_{i}$ with $\Vert p \Vert \geq \Rasym$.  The norms
    and derivatives used here are those associated with the standard
    Euclidean metric on the planes $P_{i}$.
\end{thm}

\begin{proof}
    
Only the end $E_{1}$, asymptotic to the plane $P_{1}$ and
corresponding to large positive $\lambda$, needs to be developed in
detail since the calculations for $E_{2}$ are identical to those of
$E_{1}$.  Begin with a series of preliminary estimates.

Let $A = \min \{a_{1}, \ldots, a_{n} \}$ and put $\Rasym = \sqrt{2/A}$. 
First, a calculation reveals that
\begin{equation}
    \left\vert P(a, \lambda) \right\vert \geq \min \{ A^{n}
    \lambda^{2n-2}, nA \} \, . \label{eqn:Pbelow}
\end{equation}
Thus one can estimate
\begin{equation}
    \left\vert \theta_{k}(a, \lambda) - \theta_{k} (a) \right\vert =
    \int_{\lambda}^{\infty} \frac{\dif s} {(\frac{1}{a_{k}} + s^{2} )
    \sqrt{P(a, s)}} \leq \frac {1}{n (\sqrt{A})^{n}} \cdot \frac{1}
    {\lambda^{n}} \, .  \label{eqn:thetaabove}
\end{equation}
Finally, a relation between $\lambda$ and $\Vert p \Vert$ for $p \in
E_{1}$ can be found from the equation $\Vert p \Vert = \Vert
\Psi_{a}(\lambda, \mu) \Vert = \sum \big( (x^{k}(\lambda, \mu))^{2} +
(y^{k}(\lambda, \mu))^{2} \big)$ and the fact that $\sum (\mu^{k})^{2}
= 1$.  A simple calculation shows that when $\Vert p \Vert \geq
\Rasym$, then
\begin{equation}
    \frac{\lambda}{\sqrt{2}} \leq \Vert p \Vert \leq \sqrt{2} \lambda
    \, .  \label{eqn:rels}
\end{equation}

In order to study the asymptotics of $E_{1}$ to $P_{1}$, it is most
convenient to choose new coordinates in which $P_{1}$ is transformed
into $R^{n} \times \{ 0 \}$.  For simplicity, let $\theta_{k}$ denote
the angles $\frac{\pi}{2} \delta_{1k} + \theta_{k} (a)$
and choose new coordinates $(s,t)$ for $\Rtn$ according to the
formulae:
\begin{align*}
    s^{k} &= x^{k} \cos(\theta_{k}) + y^{k} \sin(\theta_{k}) \\
    t^{k} &= -x^{k} \sin(\theta_{k}) + y^{k} \cos(\theta_{k}) \, .
\end{align*}
The coordinate functions \eqref{eqn:lawlor} of the embedding of the
Lawlor neck are given in the new coordinates by the equations:
\begin{subequations} \label{eqn:newcoord} 
    \begin{align}
	s^{k} &= \mu^{k} \sqrt{\tfrac{1}{a_{k}} + \lambda^{2}} 
	\cos \left( \theta_{k}(a, \lambda) - 
	\theta_{k} (a) \right) \\
	t^{k} &= \mu^{k} \sqrt{\tfrac{1}{a_{k}} + \lambda^{2}} 
	\sin \left( \theta_{k}(a, \lambda) - 
	\theta_{k} (a) \right) 
    \end{align}
\end{subequations}
where $(\lambda, \mu) \in \R \times \mathbf{S}^{n-1}$ and $\sum
(\mu^{k})^{2} = 1$.

The Lawlor embedding given in \eqref{eqn:newcoord} converges to the
standard polar coordinate embedding of the plane $P_{1}$, that is to
the embedding given by
\begin{align*}
    s^{k} &= \mu^{k} \lambda  \\
    t^{k} &= 0 \, ,
\end{align*}
where, again $\sum \mu_{k}^{2} = 1$.  Thus it can be shown that
$E_{1}$ is a graph over $P_{1}$ outside a sufficiently large ball ---
and in fact that a ball of radius $\Rasym$ suffices.  The coordinates
$t^{k}$ restricted to the Lawlor neck can thus be written as functions
of $s^{k}$ in the region $\Vert s \Vert \geq \Rasym$.  Since the
Lawlor neck is Lagrangian, it is the graph of the gradient of a
function $g : \R^{n} \cap \big( B_{\Rasym}(0) \big)^{c}\longrightarrow
\R$ over the asymptotic plane in this region.  Thus $t^{k} =
\frac{\partial g}{\partial s^{k}}$ are the partial derivatives of the
function $g$.

The norm of the gradient of $g$ can now be estimated in the asymptotic
region.  Divide the first of equations \eqref{eqn:newcoord} by the
second.  This gives the relation
\begin{equation} \label{eqn:sktk}
    \frac{\partial g}{\partial s^{k}} \equiv t^{k} = s^{k} \tan
    \left(\theta_{k}(a,\lambda) -
    \theta_{k}(a) \right) \, .
\end{equation}
Suppose $\Vert s \Vert \geq \Rasym$.  Now use the preliminary
estimates \eqref{eqn:Pbelow} to \eqref{eqn:rels} to estimate:
\begin{equation*}
    \left\vert \frac{\partial g}{\partial s^{k}} \right\vert \leq
    \frac{2}{n} \left( \frac{2}{A} \right)^{ n/2} \frac{1}
    {\hspace{1ex} \left\Vert s \right\Vert^{n-1}} \, .
\end{equation*}
Consequently,
\begin{equation*}
    \Vert \nabla g \Vert = \Bigg( \sum_{k} \left| \frac{\partial
    g}{\partial s^{k}} \right|^{2} \Bigg)^{\scriptscriptstyle 1/2}
    \leq \frac{\Casym}{\hspace{1ex} \left\Vert s
    \right\Vert^{n-1}}
\end{equation*}
where $\Casym$ is some constant depending only on $n$ and $A$.

The bounds on the higher derivatives of $g$ come from differentiating
equation \eqref{eqn:sktk} and estimating all quantities that appear
using the preliminary estimates once again.  Furthermore, it is
possible to derive the estimate $\Vert \nabla ^{4} g \Vert \leq C_{0}
\Vert p \Vert^{-n-2}$ using similar calculations and the estimate on
the H\"older coefficient of the third derivatives of $g$ follows from
this in the standard way.  Finally, the function $g$ itself is defined
only up to an arbitrary constant; thus it is possible to choose
\begin{equation} \label{eqn:g}
    g(s) = - \lim_{r \rightarrow \infty}
    \int_{\gamma(s, rs)} \dif g \, ,
\end{equation}
where $\gamma(s, rs)$ is the line segment between $s$ and $rs$.  This
integral is well defined because the derivatives of $g$ decay
sufficiently rapidly near infinity.  The bound on the size of $g$ now
follows from the bound on the gradient of $g$ when $\Vert s \Vert$ is
sufficiently large.
\end{proof}

\subsection{Construction of $\Mapprox$}

% The idea of the approximate submanifold -- the right annulus

The approximate submanifold $\Mapprox$ will be constructed by removing
a neighbourhood of the singularity of $M$ and smoothly reconnecting
the pieces by a rescaling $\eps N_{a}$ of the appropriate Lawlor neck
$N_{a}$.  The graphing functions of $N_{a}$ and $M$ over the tangent
cone $P_{1} \cup P_{2}$ of $M$ at the singularity will be used to
formulate a numerically precise version of this construction.  Let
$\pi_{i}$ be the orthogonal projection onto the plane $P_{i}$ and
denote by $\mathit{Ann}_{\delta}$ the annulus $\big( B_{\delta/2}(0)
\big)^{c} \cap B_{\delta}(0)$.  The graphical property of $M$ over its
tangent cone can now be phrased in the following way.

% M as a graph near the singularity

\begin{prop} \label{prop:Mgraph}
    There is a constant $K$ and a number $\delta_{0} > 0$ which depend
    only on the geometry of $M$ such that the following is true. 
    There is a function $f_{i} : P_{i} \cap Ann_{\delta}
    \longrightarrow \R$ such that if $\delta \leq \delta_{0}$, then
    $M_{i} \cap \pi_{i}^{-1} (P_{i} \cap Ann_{\delta})$ is the graph
    of $\nabla f_{i}$ over $P_{i} \cap Ann_{\delta}$.  In addition,
    the function $f_{i}$ satisfies
    \begin{equation} \label{eqn:fequalszero}
	\frac{\partial^{2} f_{i}}{\partial x^{k} \partial x^{l}}
	(0) = \frac{\partial f_{i}}{\partial x^{k}}
	(0) = f_{i}(0) = 0
    \end{equation}
    for all $k$ and $l$, along with the bounds: 
    \begin{equation} \label{eqn:fbounds}
	\big\vert f_{i} \big\vert_{0, B_{ \delta}(0)} + \delta
	\big\Vert \nabla f_{i} \big\Vert_{0, B_{\delta}(0)} +
	\delta^{2} \big\Vert \nabla^{2} \, f_{i} \big\Vert_{0,
	B_{\delta}(0)} + \delta^{3} \big\Vert \nabla^{3} f_{i}
	\big\Vert_{0, B_{\delta}(0)} + \delta^{3+\beta} \big[
	\nabla^{3} f_{i} \big]_{\beta, B_{\delta}(0)} \leq K
	\delta^{3}
    \end{equation}
    The norms and derivatives are those associated to the standard
    Euclidean metric on the planes $P_{i}$.
\end{prop}

\begin{proof}
There is some neighbourhood of the origin in which each $M_{i}$ is
graphical over its tangent plane; and in this neighbourhood, the
extrinsic curvature of the submanifolds $M_{1}$ and $M_{2}$ is bounded
in the $C^{0,\beta}$ norm by \emph{some} number because of
compactness.  For a gradient graph, the extrinsic curvature is
expressed in terms of the derivatives up to order three of the
graphing function.  The curvature condition thus translates into the
bound 
$$\Vert \nabla^{3} f_{i} \Vert_{0, B_{\delta}(0)} + \delta^{\beta} [
\nabla^{3} f_{i} ]_{\beta, B_{\delta}(0)} \leq K$$
on the graphing functions $f_{i}$.  The estimates on the second and
lower derivatives follow by integration, and use
\eqref{eqn:fequalszero}.
\end{proof}

% \eps N_{a} as a graph near the singularity

Analogous estimates can be found for the graphing function of $\eps
N_{a}$ over the components of $P$.  Recall that $N_{a}$ is the graph
of the gradient of some function $g$ outside the ball of radius
$\Rasym$.  

\begin{prop} \label{prop:epslawlorasym}
    Each asymptotic end $\eps N_{a} \cap \pi_{i}^{-1} \bigl( P_{i}
    \cap (B_{\eps \Rasym}(0)) \bigr)^{c}$ of the rescaled Lawlor neck
    is the graph of the gradient of a function $g_{\eps} : P_{i} \cap
    \bigl( B_{\eps \Rasym}(0) \bigr)^{c} \longrightarrow \R$ over the
    appropriate asymptotic plane $P_{i}$.  Moreover,
    $$g_{\eps}(x) = \eps^{2}g \left( \frac{x}{\eps} \right) \, .$$
    Thus the function $g_{\eps}$ satisfies the asymptotic inequalities
    \begin{subequations} \label{eqn:epslawlorasym}
    \begin{equation}
	\bigr\vert \, g_{\eps}(x) \bigr\vert + \Vert x \Vert
	\bigl\Vert \nabla g_{\eps}(x) \bigr\Vert + \Vert x \Vert^{2}
	\bigl\Vert \nabla^{2} \, g_{\eps}(x) \bigr\Vert + \Vert x
	\Vert^{3} \bigl\Vert \nabla^{3} g_{\eps}(x) \bigr\Vert \leq
	\frac{\Casym \eps^{n}} {\left\Vert x \right\Vert^{n-2}}
  \end{equation}
    for $x \in P_{i}$ with $\Vert x \Vert \geq \eps \Rasym$, as well 
    as the estimate
    \begin{equation}
	\big[ \nabla^{3} g_{\eps} \big]_{\beta, (B_{R}(0))^{c}} \leq 
	\frac{\Casym \eps^{n}}{R^{n+1 + \beta}}
    \end{equation}
    \end{subequations}
    for any $R \geq \eps \Rasym$.  Here, both $\Rasym$ and $\Casym$
    are as in Theorem \ref{thm:lawlorasym} and the norms and
    derivatives are those associated to the standard Euclidean metric
    on the planes $P_{i}$.
\end{prop}

\begin{proof}
The identity $\nabla g_{\eps}(x) = \eps \nabla g \left( \frac{x}{\eps}
\right)$ for points on $P$ outside the ball of radius $\eps \Rasym$
follows from scaling arguments.  Consequently, $g_{\eps} (x) =
\eps^{2} g \left( \frac{x}{\eps} \right)$ (once the constant of
integration is set to zero), and
\begin{equation}
    \bigl\vert g_{\eps}(x) \bigr\vert = \eps^{2} \left\vert g \left(
    \frac{x}{\eps} \right) \right\vert \leq \frac{\Casym \eps^{n}}{\,
    \left\Vert x \right\Vert^{n-2}} \, .
\end{equation}
The remaining estimates follow from those of Theorem
\ref{thm:lawlorasym} in a similar way.
\end{proof}

% The choice of \eps and \delta in terms of \alpha

In order to ensure that both $M$ and the correct rescaling $\eps
N_{a}$ of the Lawlor neck are close to the tangent cone $P$ in the
annulus $Ann_{\delta}$, the numbers $\eps$ and $\delta$ will be chosen
to produce
$$ \bigl\Vert \nabla^{2} \, f_{i}(x) \bigr\Vert \leq \alpha \qquad
\mbox{and} \qquad \bigl\Vert \nabla^{2} \, g_{\eps}(x) \bigr\Vert
\leq \alpha $$
for any $x \in P$ with $\frac{\delta}{2} \leq \Vert x \Vert \leq
\delta$, where $\alpha$ is any sufficiently small, positive number.

\begin{prop} \label{prop:epschoice}
    There is a number $\Rinit > 0$ and a constant $\Ceps$ depending
    only on the geometry of $M$ and $N_{a}$ such that if $0 < \alpha
    \leq \Rinit$ and the values
    \begin{equation*}
	\delta = \frac{\alpha}{K} \qquad \mbox{and} \qquad \eps \leq
	\Ceps \, \alpha^{\scriptscriptstyle 1 + 1/n}
    \end{equation*}
    are chosen, then $M_{i} \cap \pi_{i}^{-1} \left( P_{i} \cap
    Ann_{\delta} \right)$ for $i = 1$ and $2$, and $\eps N_{a} \cap
    \pi_{i}^{-1} \left( P_{i} \cap Ann_{\delta} \right)$ are graphs of
    $\nabla f_{i}$ and $\nabla g_{\eps}$, respectively, over the
    annulus $P_{i} \cap Ann_{\delta}$.  In addition, the norms of the
    Hessians of these functions are bounded above by $\alpha$; that is
    \begin{equation}
	\bigl\Vert \nabla^{2} \, f_{i}(x) \bigr\Vert \leq \alpha
	\qquad \mbox{and} \qquad \bigl\Vert \nabla^{2} \, g_{\eps}(x)
	\bigr\Vert \leq \alpha
        \label{eqn:hessbound}
    \end{equation}
    for any $x \in P_{i}$ satisfying $\frac{\delta}{2} \leq \Vert x
    \Vert \leq \delta$.
\end{prop}

\begin{proof}

By Proposition \ref{prop:Mgraph}, it is sufficient to choose $\delta =
\frac{\alpha}{K}$ in order to control the Hessian of each $f_{i}$.  To
achieve the second inequality, choose $\eps$ small enough to bring the
asymptotic region of $\eps N_{a}$ into the annulus $Ann_{\delta}$. 
Thus
$$\frac{\alpha}{2 K} \geq \eps \Rasym$$ 
is needed.  Now the bounds of Proposition \ref{prop:epslawlorasym} are
valid and thus for $x \in P$ with $\Vert x \Vert \geq
\frac{\alpha}{2K}$,
$$\bigl\Vert \nabla^{2} \, g_{\eps}(x) \bigr\Vert \leq \frac{\Casym
\eps^{n}}{\, \left\Vert x \right\Vert^{n}} \leq 2^{n} K^{n} \Casym
\frac{\eps^{n}}{\alpha^{n}} \, .$$
If $\eps \leq \frac{\alpha^{\scriptscriptstyle 1 + 1/n}}{2 K
\sqrt[n]{\Casym}}$, then $\bigl\Vert \nabla^2 \, g_{\eps}(x)
\bigr\Vert \leq \alpha$, as desired.  If $\alpha$ is sufficiently
small to begin with, then both choices of $\eps$ can be made
simultaneously.
\end{proof}

The estimates of Proposition \ref{prop:epslawlorasym} as well as the
bounds on the functions $f_{i}$ from equation \eqref{eqn:fbounds} can
now be reformulated in terms of the parameter $\alpha$.

\begin{cor} \label{cor:epschoice}
    The functions $f_{i}$ and the function $g_{\eps}$ from 
    Proposition \ref{prop:epschoice} also satisfy the bounds:
    \begin{subequations} \label{eqn:morebounds}
    \begin{equation}
	\big[ \nabla^{3} f_{i} \big]_{\beta, P_{i} \cap Ann_{\delta}}
	\leq K^{1+\beta} \alpha^{-\beta} \qquad \mbox{and} \qquad
	\big[ \nabla^{3} g_{\eps} \big]_{\beta, P_{i} \cap
	Ann_{\delta}} \leq (2K)^{1+\beta} \alpha^{-\beta}
    \end{equation}
    along with    
    \begin{equation}
        \begin{gathered}
	    \bigl\Vert \nabla^{3} f_{i} (x) \bigr\Vert \leq K \\
            \bigl\Vert \nabla f_{i} (x) \bigr\Vert \leq
            \frac{\alpha^{2}}{K} \\
            \bigl\vert f_{i}(x) \bigr\vert \leq \frac{\alpha^{3}}{K^{2}}
        \end{gathered}
        \qquad \textrm{and} \qquad
        \begin{gathered}
	    \bigl\Vert \nabla^{3} g_{\eps} (x) \bigr\Vert \leq 2K\\
            \bigl\Vert \nabla g_{\eps} (x) \bigr\Vert \leq
            \frac{\alpha^{2}}{2K} \\
            \bigl\vert g_{\eps}(x) \bigr\vert \leq
            \frac{\alpha^{3}}{4K^{2}}
        \end{gathered}
    \end{equation}
    \end{subequations}
    for any $x \in P_{i}$ satisfying $\frac{\delta}{2} \leq \Vert x 
    \Vert \leq \delta$.
\end{cor}

\begin{proof}
    These estimates can be verified by substituting for $\delta$ and
    $\eps$ in the appropriate equations \eqref{eqn:fbounds} or
    \eqref{eqn:epslawlorasym} and by using $\frac{\delta}{2} \leq \Vert x
    \Vert \leq \delta$ .
\end{proof}

% Construction of \Mapprox

Suppose now that $\alpha \leq \Rinit$ and the quantities $\delta$ and
$\eps$ have been chosen according to Proposition \ref{prop:epschoice}. 
It remains to be seen how to glue $M$ and $\eps N_{a}$ together in the
annulus $Ann_{\delta}$ in order to build the smooth submanifold
$\Mapprox$.  The `stickiness' is provided by a cut-off function: let
$\eta : \R^{2n} \rightarrow \R$ denote a positive, $C^{\infty}$
function which is equal to zero outside the ball of radius $\delta$,
one inside the ball of radius $\frac{\delta}{2}$ and interpolates
smoothly in between such that the supremum norm bounds (for the
Euclidean metric on $\R^{2n}$)
\begin{equation}
    \vert \eta \vert_{0, \R^{2n}} + \delta \, \bigl\Vert \nabla \eta
    \bigr\Vert_{0, \R^{2n}} + \delta^{2} \, \bigl\Vert \nabla^2 \,
    \eta \bigr\Vert_{0, \R^{2n}} + \delta^{3} \, \bigl\Vert \nabla^{3}
    \eta \bigr\Vert_{0, \R^{2n}} \leq C
    \label{eqn:etabounds}
\end{equation}
hold in the annulus $Ann_{\delta}$, where $C$ is some geometric
constant depending only on $n$.  The approximate submanifold
$\Mapprox$ is the union of five pieces that overlap smoothly thanks to
the cut-off function $\eta$.

\begin{defn} \label{defn:Mapprox} 
    Suppose $\alpha \leq \Rinit$ and $\eps$ and $\delta$ have been chosen
    as in Proposition \ref{prop:epschoice}.  Define the following
    submanifolds.
    \begin{align*}
	1) &\qquad M_{i}' = M_{i} \setminus \pi_{i}^{-1} ( P_{i} \cap 
	B_{\delta}(0) ) \qquad \mbox{for $i=1,2$} \\
	2) &\qquad T_{i} = \left\{ \Bigl( x, \nabla \bigl( (1-\eta)
	f_{i} + \eta g_{\eps} \bigr) (x) \Bigr) \in P_{i} \times
	P_{i}^{\perp} : \, \frac{\delta}{2} \leq \Vert x \Vert \leq
	\delta \right\} \qquad \mbox{for $i=1,2$} \\
	3) &\qquad N' = \eps N_{a} \cap \Bigl(\pi_{2}^{-1} \bigl(P_{2}
	\cap B_{\frac{\delta}{2}}(0) \bigr) \cup \pi_{2}^{-1}
	\bigl(P_{2} \cap B_{\frac{\delta}{2}}(0) \bigr) \Bigr)
    \end{align*}
    The \emph{approximate solution} to the deformation problem is the 
    submanifold
    \begin{equation*}
	\Mapprox = M_{1}' \cup T_{1} \cup N' \cup T_{2} \cup M_{2}' \, .
    \end{equation*}
    The submanifold $M_{1}' \cup M_{2}'$ is called the \emph{exterior
    region} of $\Mapprox$, the submanifold $T_{1} \cup T_{2}$ is
    called the \emph{transition region} of $\Mapprox$ and the
    submanifold $T_{1} \cup N' \cup T_{2}$ is called the \emph{neck
    region} of $\Mapprox$.
\end{defn}

The following theorem shows that the $\Mapprox$ are indeed smooth
Lagrangian approximations of $M$.

\begin{thm} \label{thm:Mapprox}
    The submanifolds $\Mapprox$, with $\alpha \leq \Rinit$, that have
    been constructed above are smooth, Lagrangian submanifolds of
    $\Rtn$ which converge to the singular submanifold $M$ in a
    pointwise sense as $\alpha \rightarrow 0$.
\end{thm}
   
\begin{proof}
The submanifold $\Mapprox$ is smooth because each $T_{i}$ overlaps
smoothly with its neighbours as a result of the way in which the
graphing functions $f_{i}$ and $g_{\eps}$ were combined. 

Furthermore, both $M$ and $\eps N_{\alpha}$ are Lagrangian and thus
$\Mapprox \setminus \bigl(T_{1} \cup T_{2} \bigr)$ is Lagrangian.  But
the transition regions $T_{i}$ are gradient graphs over Lagrangian
planes in the symplectic coordinates given by the splitting $P_{i}
\times P_{i}^{\perp}$ and thus are Lagrangian as well.

Convergence is easily verified.  The distance between $M$ and
$\Mapprox$ is actually zero outside a neighbourhood of radius on the
order of $\alpha$.  Since this neighbourhood is shrinking and 
$M_{\alpha}$ is smooth, $\Mapprox \rightarrow M$ as $\alpha
\rightarrow 0$.
\end{proof}

\subsection{Properties of $\Mapprox$}

The remainder of this section is devoted to deriving the properties of
$\Mapprox$ that will be used in the sequel.  The first such property
concerns how well $\Mapprox$ approximates a special Lagrangian
submanifold near $M$.  The upcoming estimates will be explicitly of a
global nature, unlike the previous estimates, so they must be phrased
using the induced metric of $\Mapprox$, and will thus make use of the
following bounds on the metric components and their derivatives.

\begin{lemma} \label{lemma:transmetric}
    Let $g_{kl}^{i}$ denote the coefficients of the induced metric in 
    the graphical coordinates for the transition region $T_{i}$.  
    There is a number $A_{g} > 0$ so that if $\alpha \leq A_{g}$ then 
    the functions $g_{kl}^{i}$ satisfy
    $$\left\Vert \big( g_{kl}^{i}(x) \big) - I \right\Vert \leq 1
    \quad \mbox{and} \quad \frac{1}{2} \leq \left\vert \det \big(
    g_{kl}^{i}(x) \big) \right\vert \leq 2$$
    along with the derivative bounds
    $$\left\vert \frac{\partial g_{kl}^{i}}{\partial x^{m}}
    \right\vert_{0, Ann_{\delta}} + \delta^{\beta} \left[
    \frac{\partial g_{kl}^{i}} {\partial x^{m}} \right]_{\beta,
    Ann_{\delta}} \leq 1$$
    in the annulus $P_{i} \cap Ann_{\delta}$.
\end{lemma}

\begin{proof}
    
This is a straightforward calculation using the expression for
the metric of the transition region in the graphical coordinates.
\end{proof}

It is now possible to measure the extent to which $\Mapprox$ deviates
from being special Lagrangian in a manner independent of coordinates. 
Of course, by definition, $\Mapprox$ is exactly special Lagrangian
outside the transition region.

% Im dz and Re dz bounds

\begin{prop} \label{prop:transition}
    Let $\Mapprox$ be the submanifold constructed in Theorem
    \ref{thm:Mapprox} and suppose that $\bar{h}_{\alpha}$ embeds this
    submanifold into $\Rtn$.  There exist a constant $\Cimdz$ and a
    number $\Rimdz$ with $0 < \Rimdz \leq \min \{ \Rinit, A_{g} \}$,
    both independent of $\alpha$, such that whenever $\alpha \leq
    \Rimdz$, the pull back of the form $\dif z$ satisfies the
    estimates
    \begin{subequations} \label{eqn:imdzbound}
	\begin{equation}
	   \begin{gathered}
	       \bigl\Vert \bar{h}_{\alpha}^{\, \ast} \left(
	       \mathbf{Im} \, \dif z \right) \bigr\Vert_{0, T_{1} \cup
	       T_{2}} \leq \Cimdz \alpha \\
	       1 - \Cimdz^{2} \, \alpha^{2} \leq \bigl\Vert
	       \bar{h}_{\alpha}^{\, \ast} \left( \mathbf{Re} \, \dif z
	       \right) \bigr\Vert_{0, T_{1} \cup T_{2}}^{2} \leq 1
	   \end{gathered}
       \end{equation}
       and
       \begin{equation}
	   \begin{gathered}
	       \bigl\Vert \nabla \, \bar{h}_{\alpha}^{\, \ast} \left(
	       \mathbf{Im} \, \dif z \right) \bigr\Vert_{0, T_{1} \cup
	       T_{2}} + \alpha^{ \beta} \bigl[ \nabla \,
	       \bar{h}_{\alpha}^{\, \ast} \left( \mathbf{Im} \, \dif z
	       \right) \bigr]_{\beta, T_{1} \cup T_{2}} \leq \Cimdz \\
	       \bigl\Vert \nabla \, \bar{h}_{\alpha}^{\, \ast} \left(
	       \mathbf{Re} \, \dif z \right) \bigr\Vert_{0, T_{1} \cup
	       T_{2}} + \alpha^{ \beta} \bigl[ \nabla \,
	       \bar{h}_{\alpha}^{\, \ast} \left( \mathbf{Re} \, \dif z
	       \right) \bigr]_{\beta, T_{1} \cup T_{2}} \leq \Cimdz
	   \end{gathered}
       \end{equation}
   \end{subequations}
   within the transition region.
\end{prop}

\begin{proof}
Only the estimates for $T_{1}$ need to be performed since the
computations for $T_{2}$ are identical.  Choose coordinates so that
$T_{1}$ is the graph of the Euclidean gradient of $(1-\eta) f_{1} +
\eta g_{\eps}$ over the annular region $Ann_{\delta}$ in $\R^{n}
\times \{ 0 \}$.  The operators $\nabla^2$ and $\nabla$ in the
following calculations refer to the derivatives associated to the
Euclidean metric in this coordinate system.

Begin with the bound on the size of $\bar{h}_{\alpha}^{\, \ast} \left(
\mathbf{Im} \, \dif z \right)$ in the transition region.  In their
paper, Harvey and Lawson compute the pull back of the form $\Im \,
\dif z$ to a graphical Lagrangian submanifold of $\Rtn$ \cite{hl1}. 
What they obtain is the expression
\begin{equation} \label{eqn:dethessmapprox}
    \bar{h}_{\alpha}^{\, \ast} \bigl( \mathbf{Im} \, \dif z \bigr) =
    \mathbf{Im} \left[ \mathrm{det}_{\C} \Bigl( I + \mi \, \nabla^2
    \bigl( (1-\eta) f_{1} + \eta g_{\eps} \bigr) \Bigr) \right] \dif x^{1}
    \wedge \cdots \wedge \dif x^{n} \, .
\end{equation}
By the Taylor expansion of the function $t \mapsto \det_{\C}(I+\mi t
A)$, there exists a number $r$ and of a constant $C$ such that if
$\Vert A \Vert \leq r$, then $\bigl\vert \mathbf{Im} \,
\mathrm{det}_{\C} \bigl(I + \mi \, A \bigr) \bigr\vert \leq C
\big\Vert A \big\Vert$.  Given the bounds on the functions $\eta$,
$g_{\eps}$ and $f_{1}$ and their derivatives, it is clear that there
is a number $\Rimdz$ (which should be chosen smaller than $\Rinit$ and
$A_{g}$) so that if $\alpha \leq \Rimdz$, then
\begin{align}
    \left\vert \mathbf{Im} \, \mathrm{det}_{\C} \bigl( I + \mi \,
    \nabla^2 \bigl( (1-\eta) f_{1} + \eta g_{\eps} \bigr) \bigr)
    \right\vert_{0, \R^{n}} &\leq C \left\vert \nabla^2 \bigl(
    (1-\eta) f_{1} + \eta g_{\eps} \bigr) \right\vert_{0, \R^{n}}
    \notag \\
    &\hspace{-15ex} \leq \vert 1 - \eta \vert \Vert
    \nabla^2 \, f_{i} \Vert + \vert \eta \vert
    \Vert \nabla^2 \, g_{\eps} \Vert + \Vert
    \nabla^2 \, \eta \Vert \vert f_{i} + g_{\eps}
    \vert \notag \\
    &\hspace{-15ex} \qquad + 2 \Vert \nabla \eta \Vert
    \Bigl( \Vert \nabla f_{i} \Vert + \Vert \nabla
    g_{\eps} \Vert \Bigr) \notag \\
    &\hspace{-15ex} \leq 1 \cdot \alpha + 1 \cdot \alpha + \frac{C
    n}{\delta^{2}} \cdot \frac{2 \alpha^{3}}{K^{2}} + 2 \cdot
    \frac{C}{\delta} \cdot \frac{2 \alpha^{2}}{K} \notag \\
    &\hspace{-15ex} = C \, \alpha
    \label{eqn:dethessmapprox3}
\end{align}
where the constant $C$ depends only on $\Rimdz$ and $n$.  The fact
that $\delta = \frac{\alpha}{K}$ has been used here.

Equation \eqref{eqn:dethessmapprox3} is an estimate for the desired
quantity in the Euclidean norm.  This must now be converted to a
global estimate.  Let $g_{\Mapprox}$ denote the induced metric on
$\bar{M}_ {\alpha}$.  Then,
\begin{align*} \label{eqn:here}
    \bigl\Vert \bar{h}_{\alpha}^{\, \ast} \bigl( \mathbf{Im} \, \dif z
    \bigr) \bigr\Vert_{0, T_{1} \cup T_{2}} &= \bigl\vert \mathbf{Im}
    \, \mathrm{det}_{\C} \bigl( I + \mi \, \nabla^2 \bigl(
    (1-\eta) f_{1} + \eta g_{\eps} \bigr) \bigr) \bigr\vert_{0,
    \R^{n}} \bigl\Vert \dif x^{1} \wedge \cdots \wedge \dif x^{n}
    \bigr\Vert_{0, T_{1} \cup T_{2}} \\
    &\leq \frac{C \alpha}{\left( \det \left( g_{\Mapprox} \right)
    \right)^{\scriptscriptstyle 1 / 2}} \, .
\end{align*}
For the desired estimate, invoke Lemma \ref{lemma:transmetric} to
bound the metric term in this expression from below.

The covariant derivative of $\imdz$ also only needs to be analyzed in
$T_{1}$ and once again, the local, graphical coordinates can be used
for this purpose.  The calculations are similar to the ones above,
though far more tedious because they involve the third derivatives of
$\eta$ and of the graphing functions along with their H\"older
coefficients.  Finally, the bounds on $\bar{h}_{\alpha}^{\, \ast}
\left( \mathbf{Re} \, \dif z \right)$ follow in the same way from the
identity $\Vert \dif z \Vert = 1$ proved in Harvey and Lawson's paper
\cite{hl1}.
\end{proof}

\noindent \scshape Remark: \upshape The quantities $\happrox^{\,
\ast}(\imdz)$ and $\happrox^{\, \ast}(\redz)$ estimated in the
previous proposition can be related to the Lagrangian angle function
of $\Mapprox$.  According to the defining equation
\eqref{eqn:lagangle2} of this function,
\begin{equation}
    \begin{aligned}
	\sin(\thapprox) &= \big\langle \happrox^{\, \ast} \big( \imdz
	\big) , \mathrm{Vol}_{\Mapprox} \big\rangle \\
	\cos(\thapprox) &= \big\langle \happrox^{\, \ast} \big( \redz
	\big), \mathrm{Vol}_{\Mapprox} \big\rangle \, .
    \end{aligned}
    \label{eqn:trigrel}
\end{equation}
Consequently, the trigonometric functions of $\thapprox$ satisfy the 
following estimates:
\begin{subequations} \label{eqn:trigbounds}
\begin{gather}
    \vert \sin(\thapprox) \vert_{0,T} \leq \Cimdz \alpha \qquad
    \mbox{and} \qquad \sqrt{1 - \Cimdz^{2} \alpha^{2}} \leq \vert
    \cos(\thapprox) \vert_{0,T} \leq 1 \, \\
    \intertext{as well as}
    \begin{gathered}
	\vert \nabla \sin (\thapprox) \vert_{0,T} + \alpha^{\beta}
	[\nabla \sin (\thapprox) ]_{\beta, T} \leq \Cimdz \\
	\vert \nabla \cos (\thapprox) \vert_{0,T} + \alpha^{\beta}
	[\nabla \sin (\thapprox) ]_{\beta, T} \leq \Cimdz  \, .
    \end{gathered}\\
    \intertext{These last two equations imply}
    [ \sin(\thapprox) ]_{\beta,T} \leq \Cimdz \alpha^{1-\beta} \qquad
    \mbox{and} \qquad [ \cos(\thapprox) ]_{\beta,T} \leq \Cimdz
    \alpha^{1-\beta} \, .
\end{gather}
\end{subequations}
These estimates will be used repeatedly throughout the remainder of
the paper.

\medskip The submanifold $\Mapprox$ is \emph{almost} special
Lagrangian since the quantities $\imdz$ and $\redz $ are close to
their special Lagrangian values as $\alpha \rightarrow 0$.  Since
special Lagrangian submanifolds are minimal, the mean curvature vector
of $\Mapprox$ should thus also be controllable.

\begin{prop} \label{prop:meancurv}
    Whenever $\alpha \leq \Rimdz$, the mean curvature vector of
    $\Mapprox$ satisfies the estimate
    $$\big\Vert \vec{H}_{\Mapprox} \big\Vert_{0, T_{1} \cup T_{2}} +
    \alpha^{\beta} \big[ \vec{H}_{\Mapprox} \big]_{\beta, T_{1} \cup
    T_{2}} \leq
    \Cimdz$$
    within the transition region, where $\Cimdz$ is independent of
    $\alpha$.
\end{prop}

\begin{proof}
    
Recall that the mean curvature vector $\vec{H}_{\Mapprox}$ of
$\Mapprox$ is related to the Lagrangian angle function by $
\happrox^{\, \ast} \big( \vec{H}_{\Mapprox} \elbow \omega \big) = \dif
\, \theta_{\Mapprox}$.  In the transition region, the angle is
non-constant and $\dif \, \theta_{\Mapprox} = \frac{\dif
\sin(\thapprox)} {\cos(\thapprox)}$.  The estimates of
\eqref{eqn:trigbounds} and calculations similar to those of
Proposition \ref{prop:transition} can be used to obtain
$$\Vert \dif \, \theta_{\Mapprox} \Vert_{0, T_{1} \cup T_{2}} +
\alpha^{\beta} [ \dif \, \theta_{\Mapprox} ]_{\beta, T_{1} \cup T_{2}}
\leq C \, ,$$
where $C$ is independent of $\alpha$.  The correspondence $\vec{H}
\mapsto \vec{H} \elbow \omega$ is an isometry, hence $\Vert
\vec{H}_{\Mapprox} \Vert_{0} + \alpha^{\beta} [ \vec{H}_{\Mapprox}
]_{\beta}$ in the transition region is bounded by the same geometric
constant given above.
\end{proof}

% Volume bounds

In addition to these bounds on the form $\dif z$ and the mean
curvature $\vec{H}_{\Mapprox}$, estimates on the volume of the neck
region of $\Mapprox$ will also be required in the proof of the Main
Theorem.  The volume will be calculated in two stages: first the
volume of the transition region will be estimated and then the volume
of the rescaled Lawlor neck itself will be estimated.  Begin in the
transition region $T_{1} \cup T_{2}$ by using the results of Lemma
\ref{lemma:transmetric}.

\begin{prop} \label{prop:transvol}
    If $\alpha \leq \Rimdz$, then the volume of the transition region
    satisfies $\mathit{Vol}(T_{1} \cup T_{2}) \leq C \alpha^{n}$,
    where $C$ is a constant independent of $\alpha$.
\end{prop}

\begin{proof}
  
The computation of the volume is best carried out in local coordinates
in each component $T_{i}$.  Let $g_{ij}$ denote the components of the
induced metric $g_{\Mapprox}$ in the local, graphical coordinates for
$T_{1}$ (as usual, the results for the transition region $T_{2}$ are 
identical).  Then,
\begin{align*}
    \mathit{Vol}(T_{1}) &= \int_{P_{1} \cap Ann_{\delta}}
    \sqrt{\det(g_{ij})} \, \dif x^{1} \wedge \cdots \wedge \dif x^{n}
    \\
    &\leq \sqrt{2} \int_{P_{1} \cap Ann_{\delta}} \dif x^{1} \wedge
    \cdots \wedge \dif x^{n}
\end{align*}
according to Lemma \ref{lemma:transmetric}.  Since $\delta =
\frac{\alpha}{2K}$, the last integral above is bounded by a number
proportional to $\alpha^{n}$.  Combining this fact with the inequality
above (and adding the contribution to the volume from $T_{2}$) proves
the proposition.
\end{proof}

The upper bound on the volume of $N'$, the remaining portion of the
neck region of $\Mapprox$, is slightly more subtle.  The
submanifold $N'$ is a truncation of a scaled Lawlor neck; thus its
volume exhibits the following scaling property:
\begin{align} \label{eqn:volscale}
    \mathrm{Vol} \bigl( N' \bigr) &= \mathrm{Vol} \left( \eps N_{a}
    \cap \Bigl( \pi_{1}^{-1} \bigl( P_{1} \cap B_{\frac{\delta}{2}}(0)
    \bigr) \cup \pi_{2}^{-1} \bigl( P_{2} \cap B_{\frac{\delta}{2}}(0)
    \bigr) \Bigr) \right) \notag \\
    &= \eps^{n} \mathrm{Vol} \left( N_{a} \cap \Bigl( \pi_{1}^{-1}
    \bigl( P_{1} \cap B_{\frac{\delta}{2 \eps}}(0) \bigr) \cup
    \pi_{2}^{-1} \bigl( P_{2} \cap B_{\frac{\delta}{2 \eps}}(0) \bigr)
    \Bigr) \right) \, .
\end{align}
It is now fairly easy to use the definition of the Lawlor embeddings
to conclude that
$$N_{a} \cap \Bigl( \pi_{1}^{-1} \bigl( B_{\frac{ \delta}{2 \eps}} (0)
\bigr) \cup \pi_{2}^{-1} \bigl( B_{\frac{\delta}{2 \eps}} (0) \bigr)
\Bigr) \subseteq \Psi_{a} \Bigl( \left[ \scriptstyle{-}
\tfrac{\delta}{\eps}, \tfrac{\delta}{\eps} \right] \times
\mathbf{S}^{n-1} \Bigr) \, .$$
Recall that $\frac{\delta}{\eps} = C \alpha^{-1/n}$ where $C$ is
independent of $\alpha$.  The volume of $N'$ can now be estimated in
terms of the volume of this portion of the unscaled Lawlor neck
$N_{a}$.  But to accomplish this, more precise knowledge of the metric
and the volume form on $N_{a}$ is needed.

% Metric on the Lawlor neck

\begin{lemma} \label{lemma:lawlormetric}
    The volume element on $\R \times \mathbf{S}^{n-1}$ corresponding 
    to the metric induced from $\Rtn$ by
    the Lawlor embedding $\Psi_{a}$ satisfies the estimate
     \begin{equation} \label{eqn:lawlorvol}
	\mathrm{Vol}_{N_{a}} \leq \mathrm{Vol}_{0}
    \end{equation}
    in the region $\vert \lambda \vert \geq A^{-1/2}$, where
    $\mathrm{Vol}_{0}$ is the volume form of the metric $g_{0} = 2
    (\dif \lambda)^{2} + 2 \lambda^{2} g_{S}$ on $\R \times S^{n-1}$
    and $g_{S}$ is the standard metric on the unit sphere.
\end{lemma}

\begin{proof}
The metric on $\R \times \mathbf{S}^{n-1}$ induced by the Lawlor
embedding is
\begin{equation} \label{eqn:lawlormetric}
    g_{N_{a}} = \left( \sum_{k=1}^{n} \frac{(\mu^{k})^{2}} {a_{k}^{-1}
    + \lambda^{2} } \right) \left( \lambda^{2} + \frac{1}{P} \right)
    \bigl( \dif \lambda \bigr)^{2} + \sum_{k=1}^{n} \bigl( a_{k}^{-1}
    + \lambda^{2} \bigr) \bigl( \dif \mu^{k} \bigr)^{2} \, ,
\end{equation}   
where the $\mu$ coordinates are restricted to the unit sphere $\sum
(\mu^{k})^{2} = 1$ in $\R^{n}$ and $P = P(a ,\lambda)$.  Thus it is a
straightforward calculation to estimate $\mathrm{Vol}_{N_{a}}$ in
terms of $\mathrm{Vol}_{0}$.
\end{proof}

\noindent The preceding lemma leads to the estimate on the volume of
the neck region of $\Mapprox$.

% Volume bound

\begin{prop} \label{prop:volbound}
    There is a constant $\Cvoltr$ independent of $\alpha$ such that if
    $\alpha \leq \Rimdz$, then the volume of the neck region of
    $\Mapprox$ satisfies $\mathit{Vol}(T_{1} \cup N' \cup T_{2}) \leq
    \Cvoltr \alpha^{n}$.
\end{prop}    
    
\begin{proof}
According to Lemma \ref{lemma:lawlormetric} and the scaling property
\eqref{eqn:volscale} of the volume of the neck region of $\Mapprox$,
the following estimate is valid:
\begin{align*}
    \mathit{Vol}(T_{1} \cup N' \cup T_{2}) &= \mathit{Vol}( T_{1} \cup
    T_{2}) + \mathit{Vol} ( N' ) \\
    &= C \alpha^{n} \left(1 + \alpha \int_{\left[ \scriptstyle{-}
    C \alpha^{-1/n}, C \alpha^{-1/n} \right] \times \mathbf{S}^{n-1}}
    \mathrm{Vol}_{N_{a}} \right) \, ,
\end{align*}
where the result of Proposition \ref{prop:transvol} as well as the
values of $\delta$ and $\eps$ in terms of $\alpha$ have been used. 
According to equation \eqref{eqn:lawlorvol}, the volume form
$\mathrm{Vol}_{N_{a}}$ can be replaced by the volume form
$\mathrm{Vol}_{0}$ in the region where $\vert \lambda \vert \geq
A^{-1/2}$.  Consequently,
\begin{equation*}
    \mathit{Vol}(N') \leq C + C' \int_{ \left[ A^{-1/2}, C
    \alpha^{-1/n} \right] \times \mathbf{S}^{n-1} } \mathrm{Vol}_{0}
    = C  + C' \alpha^{-1}
\end{equation*}
Combining the two results above and modifying the constants yields the
desired estimate.
\end{proof}

\noindent An obvious corollary to the volume bound on the neck region
of $\Mapprox$ is the following result.

\begin{cor} \label{cor:uniformly}
    The volume of $\Mapprox$ is uniformly bounded above and below
    whenever $\alpha \leq \Rvol$.
\end{cor}

\subsection{The First Neumann Eigenvalue of $\Mapprox$}

The approximating submanifolds $\Mapprox$ are converging to a singular
variety as $\alpha \rightarrow 0$, which will have the effect of
causing $\alpha$-dependent quantities on $\Mapprox$ to degenerate as
$\alpha \rightarrow 0$.  The most important of these in the context of
this paper is the first Neumann eigenvalue of the Laplacian on
$\Mapprox$, which tends towards zero as $\alpha \rightarrow 0$.  The
specific functional dependence of this quantity on $\alpha$ is a
result of the geometry of $\Mapprox$, and is established in the
following proposition.

\begin{prop} \label{prop:firsteigen}
    The first Neumann eigenvalue of $\Mapprox$ satisfies $\nu_{1} \leq
    \Cnuone \, \alpha^{n-2}$, where $\Cnuone$ is a constant
    independent of $\alpha$.
\end{prop}

\begin{proof}
Recall that the first Neumann eigenvalue of $\Mapprox$ is equal to
$$\nu_{1} = \inf \left\{ \frac{\int_{\Mapprox} \Vert \nabla u
\Vert^{2}}{\int_{\Mapprox} u^{2}} : u \in H^{1}(\Mapprox) \,
\mbox{ and } \int_{\Mapprox} u = 0 \right\} \, ,$$
where $H^{1}(\Mapprox)$ are the $L^{2}$ functions of $\Mapprox$ whose
first weak derivatives are also in $L^{2}$.  Choose a function $u$ on
$\Mapprox$ which is equal to 1 in $M_{1}'$, equal to $-1$ in $M_{2}'
\cap (\mathcal{U} _{\alpha})^{c}$ and interpolates between these
values in the neck region of $\Mapprox$.  The interpolation can thus
be made in such a way that $\vert u \vert \leq 1$ and $\Vert \nabla u
\Vert \leq C \alpha^{-1}$ for some constant $C$ independent of
$\alpha$, and a constant can be subtracted to ensure that the function
$u$ has integral zero.  Now estimate as follows:
\begin{equation} \label{eqn:eigen1}
    \int_{\Mapprox} \Vert \nabla u \Vert^{2} = \int_{\Mapprox \cap
    \mathcal{U}_{\alpha}} \Vert \nabla u \Vert^{2} \leq C \alpha^{n-2}
    \, ,
\end{equation}
using the volume estimates of Proposition \ref{prop:volbound}.  Next,
since $u$ differs from a constant outside a neighbourhood of size
proportional to $\alpha$, the integral $\int_{\Mapprox} u^{2}$ is
bounded below by a constant independent of $\alpha$.  Taking this fact
together with \eqref{eqn:eigen1} yields the desired estimate.
\end{proof}

\section{Deformations of $\Mapprox$}

\subsection{Parametrizing Lagrangian Embeddings of $\Mapprox$}

The next step in the proof of the Main Theorem is to define the
parametrization of Lagrangian embeddings of $\Mapprox$ near $\happrox$
that will be used to set up the PDE which must be solved using the
Inverse Function Theorem.  Unless the parametrization of Lagrangian
submanifolds near $\Mapprox$ is chosen with care, the strong
dependence of $\nu_{1}$ on $\alpha$ will manifest itself in the dependence
of $C_{L}(\alpha)$ on $\alpha$.  For, suppose that the Banach space
parametrizing nearby Lagrangian submanifolds contained the first
Neumann eigenfunction of the Laplacian on $\Mapprox$.  Denote this
function by $S_{\alpha}$.  Then,
\begin{equation*}
    \Dif F_{\alpha} (0,0) (S_{\alpha}, 0) = \nu_{1} \cos
    (\theta_{\Mapprox}) \, S_{\alpha} - \sin(\theta_{\Mapprox} ) \,
    \big\langle \vec{H}_{\Mapprox} , \nabla S_{\alpha} 
    \big\rangle_{\Mapprox} \, ,    
    \label{eqn:nuest1}
\end{equation*}
according to calculation of the linearization performed in Section
2.3, and indicates that the constant $C_{L}(\alpha)$ would be less
than some quantity proportional to $\nu_{1}$.  Thus the parametrizing
Banach space should exclude $S_{\alpha}$ in order to achieve a better
estimate of $C_{L}(\alpha)$.

With this observation in mind, begin by choosing a Banach space over
which to parametrize embeddings.  Let $Z$ be the inward pointing unit
normal of $\partial \Mapprox$.

\begin{defn} \label{defn:parambanach} 
    Define the Banach space $\mathcal{B}_{\alpha} = \mathcal{B}_{1,
    \alpha} \times \R$ where
    \begin{equation}
	\mathcal{B}_{1, \alpha} = \left\{ H \in C^{2, \beta}(\Mapprox)
	\: : \: Z(H) \big|_{\partial \Mapprox} = 0 \quad \mbox{and}
	\quad \int_{\Mapprox} \!  H \: = \int_{\Mapprox} \!  H \cdot
	S_{\alpha} \, = 0 \right\} \, .
	\label{eqn:banach1}
    \end{equation}
    Here, $\beta \in (0,1)$ is the degree of H\"older continuity and
    will be chosen later.
\end{defn}

\noindent \scshape Remark: \upshape The constant functions (the kernel
of the linearized operator) and the first eigenfunction $S_{\alpha}$
are excluded from $\mathcal{B}_{1,\alpha}$ because the integral
conditions in the definition above ensure that the functions in
$\mathcal{B}_{1,\alpha}$ are $L^{2}$-orthogonal to the constants and
$S_{\alpha}$.  The $\R$ factor in the definition above will be related
to a deformation of $\Mapprox$ taking $\partial \Mapprox$ away from
its scaffold $W$ and is will be used to guarantee the surjectivity of
the linearized operator.

\medskip The parametrization of Lagrangian embeddings of $\Mapprox$
near $\happrox$ will be of the following form.  To each $(H, b) \in
\mathcal{B}_{1,\alpha} \times \R$, associate the embedding
$\phi_{v_{e}}^{b} \!  \circ \phi_{H_{e}}^{1} \!  \circ \happrox$,
where the $\phi^{1} _{H_{e}}$ term is the time-1 Hamiltonian flow of a
suitable extension of the function $H$, and the $\phi^{b}_{v_{e}}$
term is the time-$b$ Hamiltonian flow of the extension of a carefully
chosen function $v : \partial \Mapprox \rightarrow \R$.  In order to
describe the extended functions $H_{e}$ and $v_{e}$ in greater detail,
a lemma concerning the structure of tubular neighbourhoods of
$\partial \Mapprox$ is needed first.  

\begin{lemma} \label{lemma:scaffolddarboux}
    Let $W$ be a symplectic submanifold of codimension 2 in $\Rtn$ and
    suppose that $L$ is a Lagrangian submanifold with boundary
    $\partial L \subset W$.  Then there exists a tubular neighbourhood
    $\mathcal{U}$ of the boundary and a symplectomorphism $\psi :
    \mathcal{U} \longrightarrow \left(T^{\ast}\partial L\right) \times
    \R^{2}$ with the following properties:
    \begin{enumerate}
	\item $\psi \big( W \cap \mathcal{U} \big) \subset T^{\ast}
	(\partial L) \times \{0,0\}$;
	\item $\psi(\partial L) = \partial L \times \{0,0\}$;
	\item $\psi \big( L \cap \mathcal{U} \big) \subset \partial L
	\times \R_{+} \times \{0\}$; and
	\item $\psi_{\ast} Z = \frac{\partial}{\partial s^{1}}$ along
	$\partial L$, where $(s^{1}, s^{2})$ are the coordinates for
	the $\R^{2}$ factor and $Z$ is the inward unit normal of
	$\partial L$.
    \end{enumerate}
\end{lemma}

\begin{proof} The proof of this lemma can be found in \cite{me2}.
\end{proof}

\noindent \scshape Definition of $H_{e}$ \upshape \medskip

\noindent The function $H$ is extended in two stages: $H$ is first
extended in an obvious manner to a neighbourhood of $\Mapprox$; then
this extension is modified near the boundary to ensure that the
Hamiltonian deformation $\phi_{H_{e}}^{t}$ keeps $\partial \Mapprox$
confined to $W$.  The preliminary extension will be carried out in a
tubular neighbourhood of $\Mapprox$ chosen according to the following
considerations.

Recall that the submanifold $\Mapprox$ consists of two large pieces
$M_{1}'$ and $M_{2}'$ which are connected by a thin neck.  Hence there
is a tubular neighbourhood for $\Mapprox$ which is large around each
$M_{i}'$ but small in the vicinity of the neck.  Denote this
neighbourhood by $\mathcal{U}_{1}$.  Without loss of generality, this
tubular neighbourhood is symplectomorphic to a neighbourhood of the
zero section in $T^{\ast} \Mapprox$.  Let $\tau : \Rtn \longrightarrow
\R$ be a smooth function equal to 1 in a tubular neighbourhood
$\mathcal{U}_{1}'$ of $\Mapprox$ contained in $\mathcal{U}_{1}$ and
that vanishes outside $\mathcal{U}_{1}$.  Furthermore, suppose $\tau$
satisfies $\vert \tau \vert + \alpha \Vert \nabla \tau \Vert \leq C$
(the factor $\alpha$ arises because of the narrowness
$\mathcal{U}_{1}$ near the neck region of $\Mapprox$).  Now define
$H_{1} : \mathcal{U}_{1} \longrightarrow \R$ in Lagrangian
neighbourhood coordinates by:
$$H_{1}(q,p) = \tau (q,p) \, H(q)$$
and to extend $H_{1}$ outside $\mathcal{U}_{1}$, simply make it zero.

The previous extension must now be modified near the boundary.  First,
choose a tubular neighbourhood $\mathcal{U}_{2}$ of $\partial
\Mapprox$ in which symplectic coordinates can be chosen as in Lemma
\ref{lemma:scaffolddarboux}.  Suppose that the number $w_{2}$
characterizes the width of $\mathcal{U}_{2}$ in the sense that if $(x,
y ; s^{1}, s^{2})$ denotes a point in $T^{\ast} \partial \Mapprox
\times \R^{2}$, then it belongs to $\mathcal{U}_{2}$ if $\mathit{max}
\left\{ \Vert y \Vert, \vert s^{1} \vert, \vert s^{2} \vert \right\} <
w_{2}$.  Now let $\eta_{0} : \R \longrightarrow \R$ be a smooth,
positive cut-off function that is bounded by 1, vanishes outside the
interval $[0,1]$, and is equal to 1 inside the interval $[0,1/2]$. 
Define the extension $H_{2} : \mathcal{U}_{2} \longrightarrow \R$ in
the $(x,y\, ; s^{1}, s^{2})$ coordinates by
$$ H_{2}(x,y; \, s^{1},s^{2}) = \eta_{0} \!  \left( \frac{\Vert y
\Vert}{w_{2}} \right) \eta_{0} \!  \left( \frac{ |s^{2}|}{w_{2}}
\right) \eta_{0} \!  \left( \frac{|s^{1}|}{w_{2}} \right) H(x,s^{1})
$$
and once again, let it to be zero outside $\mathcal{U}_{2}$. 

The complete extension of the function $H$ that is desired will come
from smoothly interpolating between the extensions $H_{1}$ and
$H_{2}$.  Denote by $\eta_{1}$ the function of $\Rtn$ given by
extending the function $\eta_{0}(\vert s^{1} \vert / w_{2})$ defined
in $\mathcal{U}_{2}$ to all of $\Rtn$ by setting it equal to zero in
$\Rtn \setminus \mathcal{U}_{2}$.

\begin{defn} \label{defn:hamext}
    For any $H \in C^{2,\beta}(\Mapprox)$, the extension of $H$ to
    $\Rtn$ is denoted by $H_{e}$ and is defined by the equation
    \begin{equation} \label{eqn:hamext}
	H_{e}(x) = \big( 1 - \eta_{1} (x) \big) H_{1}(x)
	+ H_{2}(x)
    \end{equation}
    for any point $x$ in $\Rtn$.
\end{defn}

The following proposition shows that, with the proper boundary
conditions on the functions $H$, the deformations $\phi_{H_{e}}^{1}$
deform $\Mapprox$ in the desired manner.

% Extension of H with Neumann boundary conditions

\begin{prop} \label{prop:hamflow}
    Let $H \in C^{2,\beta} (\Mapprox)$ and suppose that $H$ satisfies
    $Z(H) \big|_{\partial \Mapprox} = 0$.  Then the family of
    submanifolds $\phi_{H_{e}}^{\, t}(\Mapprox)$ is a Lagrangian
    deformation of $\Mapprox$ and the family of boundaries
    $\phi_{H_{e}}^{\, t}( \partial \Mapprox )$ remains on the scaffold
    $W$.  Furthermore, the deformation vector field associated to this
    family of submanifolds is the vector field $X_{H}$ which satisfies
    $X_{H} \elbow \omega = \dif H$ on $\Mapprox$.
\end{prop}

\begin{proof}  The proof of this proposition is a straightforward 
algebraic calculation and can be found in \cite{me}.
\end{proof}

\medskip \noindent \scshape Definition of $v_{e}$ \upshape \medskip

\noindent First define the function $v : \partial \Mapprox \rightarrow
\R$.  Recall that $\partial \Mapprox$ consists of the disjoint union
of the two separate boundary components $\partial M_{1}$ and $\partial
M_{2}$.

\begin{defn} \label{defn:bdfunction}
    Define $v: \partial \Mapprox \longrightarrow \R$ by the
    prescription
    \begin{equation}
        v(x) = 
        \begin{cases}
            \frac{1}{\mathit{Vol}(\partial M_{1})} &\qquad x \in 
            \partial M_{1} \\
            \frac{-1}{\mathit{Vol}(\partial M_{2})} &\qquad x \in 
            \partial M_{2} 
        \end{cases}  
        \label{eqn:bdfunction}
    \end{equation}
    Note that this definition implies that $\int_{\partial \Mapprox} v
    = 0$; this fact will be used later.
\end{defn}

Use the symplectic coordinates for the neighbourhood $\mathcal{U}_{2}$
guaranteed by Lemma \ref{lemma:scaffolddarboux} to extend $v$.  Define
the extended function $v_{e} : \mathcal{U}_{2} \longrightarrow \R$ by
\begin{equation}
    v_{e}(x,y; \, s^{1}, s^{2}) = v(x) \, s^{1} \, \eta_{0} \! \left(
    \frac{ \Vert y \Vert}{w_{2}} \right) \eta_{0} \! \left(
    \frac{ |s^{1}|}{w_{2}} \right) \eta_{0} \! \left(
    \frac{ |s^{2}|}{w_{2}} \right) \, .
    \label{eqn:bdfnext}
\end{equation}
Extend $v_{e}$ outside $\mathcal{U}_{2}$ by setting it equal to zero
there.  Note that the function $v_{e}$ satisfies the property $Z
(v_{e}) = \frac{\partial} {\partial s^{1}} v_{e} = v$ at the boundary.

The deformation vector field of the Hamiltonian deformation
$\phi^{b}_{v_{e}}$ is the vector field that satisfies $X_{v_{e}}
\elbow \omega = \dif v_{e}$ on $\Mapprox$.  This vector field can be
obtained from similar calculations to those carried out above in the
definition of $H_{e}$, resulting in the expression
$$X_{v_{e}} \Big|_{\partial \Mapprox} = - v(x) \frac{\partial}
{\partial s^{2}} \, .$$
This is perpendicular to $W$.  Consequently, the deformation
$\phi^{b}_{v_{e}}$ moves $\partial \Mapprox$ away from $W$.

\medskip The differential operator $F_{\alpha}$ describing minimal
Lagrangian submanifolds near $\Mapprox$ comes from combining the
general consideration of Section 2.2 with the specific parametrization
of nearby Lagrangian embeddings constructed in the previous
paragraphs.

\begin{defn} \label{defn:defop}
    Let $\mathcal{B}_{1,\alpha}$ be the Banach space of functions 
    given in equation \eqref{eqn:banach1}.  Define the map 
    $F_{\alpha} : \mathcal{B}_{1,\alpha} \times \R^{2} \rightarrow 
    C^{0,\beta}(\Mapprox)$ by
    \begin{equation}
	F_{\alpha}(H,\theta,b) = \left\langle \left( \phi_{v_{e}}^{b}
	\!  \circ \phi_{H_{e}}^{1} \!  \circ \happrox \right)^{\, \ast} 
	\, \Im (\me^{\mi \theta} \dif z ), \mathrm{Vol}_{\Mapprox} 
	\right\rangle_{\Mapprox} \, .
	\label{eqn:defop}
    \end{equation}
\end{defn}

The linearization of this operator at the point $(0,0,0)$ is now
easily calculated by applying the general result of Proposition
\ref{prop:lin}, and is given by the following proposition.

\begin{prop} \label{prop:linpde}
    Let $(u, a, b) \in \mathcal{B}_{1,\alpha} \times \R^{2}$.  Then 
    the linearization of $F_{\alpha}$ at $(0,0,0)$ acting on the 
    point $(u,a,b)$ is given by the formula
    \begin{align}
	\Dif F_{\alpha}(0,0,0)(u,a,b) &= - \cos(\theta_{\Mapprox}) \,
	\Delta_{\Mapprox} u - \sin(\theta_{\Mapprox}) \, \big\langle
	\vec{H} _{\Mapprox} , \nabla u \big\rangle \notag \\
	&\qquad + a \cos(\thapprox) - b \, \Delta_{\Mapprox} v_{e} \,
	.
	\label{eqn:linpde}
    \end{align}
\end{prop}

\subsection{The Weighted Schauder Norm}
 
Lagrangian submanifolds close to $\Mapprox$ are parametrized over the
Banach space $\mathcal{B}_{1,\alpha} \times \R$ according to
Definition \ref{defn:parambanach}.  It remains to choose the norms for
both $\mathcal{B}_{1,\alpha}$ and $C^{0,\beta} (\Mapprox)$ in which
the measurements of $C_{L}(\alpha)$, $C_{N}(\alpha)$ and $F_{\alpha}
(0,0,0)$ will be made.  The usual $C^{k,\beta}$ norms will turn out to
be unsuitable for exhibiting the explicit dependence on $\alpha$ of
these quantities.  The optimal norms for solving the deformation
problem are \emph{weighted} Schauder norms, whose \emph{weight
function} compensates for the two different scales of $\Mapprox$, one
within the neck region and one outside it.

The two different scales of $\Mapprox$ are encoded in the \emph{radius
of uniformity} of the induced metric $g_{\Mapprox}$ and is defined
as follows.  For each point $p$ of $\Mapprox$, there is a radius
$r(p)$ so that geodesic normal coordinates can be used in the ball
$B_{r(p)}(p)$ and the metric coefficients are uniformly
$C^{1,\beta}$-bounded in these coordinates.  That is,
\begin{equation}
    \left\vert g_{ij} - \delta_{ij} \right\vert_{1,\beta,
    B_{r(p)}(p)}^{\ast} \leq 1 \, ,
    \label{eqn:unifmetric}
\end{equation}
where $\vert \cdot \vert^{\ast}_{1,\beta}$ is the scale-invariant
Schauder norm of $\R^{n}$ given by
\begin{equation}
    \vert u \vert^{\ast}_{k, \beta, B_{R}} = \vert u \vert_{0, B_{R}}
    + R \vert \nabla u \vert_{0, B_{R}} + \cdots + R^{k} \vert
    \nabla^{k} u \vert_{0, B_{R}} + R^{k+\beta} [ \nabla^{k} u
    ]_{\beta, B_{R}} \, .
    \label{eqn:localnorm}
\end{equation}
Here, $| \cdot |_{0, B_{R}}$ denotes the supremum norm over $B_{R}$
and $[ \cdot ]_{\beta, B_{R}}$ denotes the H\"older coefficient in
$B_{R}$.

If $p$ is in the interior region of $\Mapprox$ then $r(p) \geq \eps
r_{1}(\frac{p}{\eps})$ where $r_{1}$ gives the radius of uniformity of
the metric in the unscaled Lawlor neck $N_{1}$.  The behaviour of
$r_{1}$ in $N_{1}$ is as follows.  In $N_{1} \cap B_{a}(0)$ for
some fixed radius $a$, $r_{1}(p)$ is bounded below by some fixed
number $R_{int}$.  In $N_{1} \cap \big(B_{a}(0) \big)^{c}$, the
bound $r_{1}(p) \geq m_{1} \Vert p \Vert$, for some fixed rate $m_{1}$
because if $N_{1} \cap \big( B_{a}(0) \big)^{c}$ is close to a
cone if $a$ is sufficiently large.  Translating this behaviour
back to $\Mapprox$ is a matter of rescaling: $r(p) \geq \eps R_{int}$
in $\eps N_{1} \cap B_{\eps a}(0)$ and $r(p) \geq m \Vert p \Vert$
in $\eps N_{1} \cap \big( B_{\eps a}(0) \big)^{c}$.  Next, the
value of $r(p)$ for $p$ in the exterior region $M_{1}' \cup M_{2}'$ is
independent of $\alpha$ and is bounded below by some value $R_{ext}$. 
Thus $r(p)$ grows from size $\eps R_{int}$ with constant rate until it
reaches the value $R_{ext}$, which it attains in some ball of
radius independent of $\alpha$.

The behaviour of the radius of uniformity of the induced metric of
$\Mapprox$ suggests that the desired weight function should be one
which interpolates between the value $\eps R_{int}$ in a ball of
radius $\eps a$ contained in the interior region $N'$ and the
value $R_{ext}$ in some part of the exterior region $M_{1}' \cup
M_{2}'$.  Let $\Rweight = \min \{ R_{int}, R_{ext} \}$ and choose any
$\beta \in (0,1)$.  It is easy to verify that a smooth and increasing
function of the following form exists on $\Mapprox$.

\begin{defn} \label{defn:weight} 
    Let $\rho : \Mapprox \longrightarrow \R$ be a function of the form:
    \begin{equation*}
        \rho(x) =
        \begin{cases}
	    \Rweight \, \eps & x \in \Mapprox \cap B_{\eps a} (0) \\
	    \mathrm{Interpolation} & x \in \Mapprox \cap \mathit{Ann}
	    _{(\eps a, \, \eps^{\beta} b}(0) )\\
	    \Rweight & x \in \Mapprox \cap \big( B_{\eps^{\beta}
	    b}(0) \big)^{c}
        \end{cases}
    \end{equation*}
where $b \in \R$ and the interpolation can be chosen so that
$\rho$ is smooth and so that the following additional properties hold. 

\medskip \noindent \scshape Property 1: \upshape
The gradient of $\rho$ satisfies the bound $\Vert \nabla \rho \Vert
\leq K \eps^{-\beta}$ for some $K$ independent of $\alpha$ because
$\rho$ grows from size $\eps \Rweight$ to size $\Rweight$ in an
annular region of width on the order of $\eps^{\beta}$.

\medskip \noindent \scshape Property 2: \upshape
Since $\rho = \eps R$ in $B_{\eps a}(0)$ and then grows to size 1
in the annulus $Ann_{(\eps a, \eps^{\beta} b)} (0)$, there is
a constant $C$ independent of $\alpha$ so that $\rho(x) \geq C \Vert x
\Vert$ in this annulus.

\medskip \noindent \scshape Property 3: \upshape 
Suppose $k \geq 1$.  Then $[\rho^{k} ]_{\beta, \Mapprox} \leq C
\eps^{-\beta} $, where $C$ is independent of $\alpha$.

\noindent \itshape Proof: \upshape Choose two points $x$ and 
$x'$ in $\Mapprox$ and let $\gamma$ be the length-minimizing geodesic 
connecting these two points.  Now calculate,
\begin{align}
    \frac{ \left\vert \rho^{k}(x) - \rho^{k}(x') \right\vert }{\left(
    \mathit{dist}(x,x') \right)^{\beta}} &= \frac{\left\vert
    \int_{\gamma} \langle \nabla (\rho^{k}), \dot{\gamma} \rangle \dif
    s \right\vert}{\left( \mathit{dist}(x,x') \right)^{\beta}}\notag
    \\
    &\leq \Vert \nabla \rho^{k} \Vert_{0,\Mapprox} \left(
    \mathit{dist} (x,x') \right)^{1-\beta} \notag \\
    &\leq k \vert \rho^{k-1} \vert_{0, \Mapprox} \, \Vert \nabla \rho
    \Vert_{0, \Mapprox} \left( \mathit{diam} (\Mapprox) \right)
    ^{1-\beta} \, .
    \label{eqn:haveiusedit}
\end{align}
The distance function used here is the distance function on $\Mapprox$
corresponding to the induced metric.  All quantities except the
gradient term in the estimate \eqref{eqn:haveiusedit} are bounded
above by constants independent of $\alpha$, whereas the gradient is
bounded by $K \eps ^{-\beta}$.  \hfill \qedsymbol

\medskip \noindent \scshape Property 4: \upshape 
Suppose $\beta \in (0,1)$.  Then $[\rho^{\beta} ]_{\beta, \Mapprox}
\leq C \eps^{-\beta}$, where $C$ is independent of $\alpha$.

\noindent \itshape Proof: \upshape Choose two points $x$ and $x'$ in 
$\Mapprox$ which satisfy $\mathit{dist}(x,x') \geq \eps$.  Then,
\begin{align*}
    \frac{\vert \rho^{\beta}(x) - \rho^{\beta}(x') \vert}{\big(
    \mathit{dist}(x,x') \big)^{\beta}} &\leq 2 \vert \rho^{\beta} 
    \vert_{0,\Mapprox} \eps^{-\beta} \\
    &\leq 2 R^{\beta} \, \eps^{-\beta} \, .
\end{align*}
Next, suppose $\mathit{dist}(x,x') < \eps$.  Then,
\begin{align*}
    \frac{\vert \rho^{\beta}(x) - \rho^{\beta}(x') \vert}{\big(
    \mathit{dist}(x,x') \big)^{\beta}} &\leq \Vert \nabla 
    \rho^{\beta} \Vert_{0,\Mapprox} \eps^{1-\beta} \\
    &\leq \beta \vert \rho^{\beta - 1} \vert_{0, \Mapprox} \, \Vert
    \nabla \rho \Vert_{0, \Mapprox} \eps^{1-\beta} \\
    &\leq \beta (\Rweight)^{\beta - 1} K \eps^{-\beta} \, ,
\end{align*}
since $\rho$ is bounded below by $\Rweight \eps$.  Taking the supremum
over $\Mapprox$ yields the desired result.  \hfill \qedsymbol
    
\medskip \noindent \scshape Property 5: \upshape
Suppose $p < n$.  Then $\Vert \rho^{-1}\Vert_{L^{p}(\Mapprox)} \leq
C$, where $C$ is a constant independent of $\alpha$.

\noindent \itshape Proof: \upshape Choose $p < n$ and calculate as
follows. 
\begin{align*}
    \int_{\Mapprox} \rho^{-p} &= \int_{\Mapprox \cap B_{\eps
    a}(0)} \rho^{-p} + \int_{\Mapprox \cap \mathit{Ann}_{( \eps
    a, \, \eps^{\beta} b) }(0)} \rho^{-p} + \int_{\Mapprox
    \cap \big( B_{\eps^{\beta} b}(0) \big)^{c}} \rho^{-p} \\
    &\leq \left(\Rweight \, \eps \right)^{-p} \, \mathit{Vol} \big(
    B_{\eps a}(0) \big) + C \int_{\eps a}^{\eps^{\beta} b}
    s^{n-p-1} \dif s + \int_{\big( B_{\eps^{\beta} b}(0)
    \big)^{c}} R^{-p} \qquad \mbox{(by Property 2)} \\
    &\leq C \left( 1 + \eps^{n-p} \right) \, ,
\end{align*}
which is bounded when $p < n$. \hfill \qedsymbol

\medskip \noindent \scshape Property 6: \upshape
There exists a constant $C$ independent of $\alpha$ so that $\rho(x)
\eps^{\beta} \leq C r(x)$ for every $x \in \Mapprox$, where $r(x)$ is
the radius of uniformity of the metric coefficients at the point $x$.

\noindent \itshape Proof: \upshape The result is true by definition of
$\rho$ and $r$ in the ball $B_{\eps a}(0)$ as well as well as
outside some large ball of radius independent of $\alpha$ where $r(p)$
is bounded below.  Between these two regions where $r$ grows linearly,
the gradient bound on $\rho$ implies that $\rho(x) \leq \rho(0) + K
\eps^{-\beta} \Vert x \Vert$ and thus
\begin{align*}
    \rho(x) \eps^{\beta} &\leq \eps^{1+\beta} + K  \Vert x \Vert \\
    &\leq \frac{1}{m} \left( \frac{\eps^{\beta}}{a} + K \right) 
    r(x)
\end{align*}
using the fact that $r(x) \geq m \Vert x \Vert$ in the region in
question.  This leads to the desired estimate.  \hfill \qedsymbol

\end{defn}

The $\rho$-weighted Schauder norms on $\Mapprox$ that will be used to
carry out the estimates of the Main Theorem are defined as follows.

\begin{defn} \label{defn:schauder}
    Let $u$ be any $C^{k,\beta}$ function on $\Mapprox$.  The weighted
    $(k,\beta)$-Schauder norm of $u$ will be denoted $\vert u
    \vert_{C^{k,\beta}_{\rho}(\Mapprox)}$ and is defined as
    $$\big\vert u \big\vert_{C^{k,\beta}_{\rho}(\Mapprox)} = \big\vert
    u \big\vert_{0, \Mapprox} + \big\vert \, \rho \,  \nabla u
     \, \big\vert_{0, \Mapprox} + \cdots + \big\vert \, \rho^{k}
    \,  \nabla^{k} u  \, \big\vert_{0, \Mapprox} + \big[ \,
    \rho^{k + \beta}\,  \nabla^{k} u  \big]_{\beta,
    \Mapprox} \, .$$
    Here, $| \cdot |_{0, \Mapprox}$ denotes the supremum norm and $[
    \, \cdot \, ]_{\beta, \Mapprox}$ denotes the H\"older coefficient
    on $\Mapprox$; the norms and derivatives in this expression are
    those corresponding to the induced metric on $\Mapprox$.
\end{defn}

The differential operator $F_{\alpha}$ maps between the Banach spaces
$\mathcal{B}_{1,\alpha} \times \R^{2}$ and $C^{0,\beta} (\Mapprox)$. 
The analysis that will be performed in order to solve the equation
$F_{\alpha}(u, a, b) = 0$ requires a choice of norm to be made in both
of these spaces.
\begin{itemize}
    \item Use the weighted norm $| u |_{C^{2,\beta}_{\rho}
    (\Mapprox)}$ for functions $u$ in the Banach space
    $\mathcal{B}_{1,\alpha}$.
    \item Use the product norm $\Vert (u,a,b) \Vert_{C^{2,\beta}
    _{\rho}(\Mapprox) \times \R^{2}} = \left(|u|^{2}_{C^{2,
    \beta}_{\rho} (\Mapprox)} + a^{2} + b^{2} \right)
    ^{\scriptscriptstyle 1/2}$ for elements in the Banach space
    $\mathcal{B} _{1,\alpha} \times \R^{2}$.
    \item Use the weighted norm $\left\vert \rho^{2} f \right\vert_{ 
    C^{0,\beta}_{\rho}(\Mapprox)}$ for functions $f$ in the Banach 
    space $C^{0,\beta}(\Mapprox)$.
\end{itemize}

It is easy to verify that $\mathcal{B}$ and $C^{0,\beta} (\Mapprox)$
are indeed Banach spaces with these norms.  Furthermore, $\Dif
F_{\alpha}(0,0,0)$ is a bounded operator in the weighted norms.

% Boundedness of the linearization

\begin{prop} \label{prop:linbounded}
    The linearization $\Dif F_{\alpha}(0,0,0)$ of $F_{\alpha}$ at the
    origin is a bounded operator between the space $\mathcal{B}_{1,
    \alpha} \times \R^{2}$ with norm $ \Vert \cdot \Vert_{C^{2,
    \beta}_{\rho} (\Mapprox) \times R^{2}}$ and the space $
    C^{0,\beta} (\Mapprox)$ with norm $| \rho^{2} \cdot
    |_{C^{0,\beta}_{\rho}(\Mapprox)}$.
\end{prop}

\begin{proof}
    
Let $(u,a,b) \in \mathcal{B}_{1,\alpha} \times \R^{2}$.  Then using 
the form of $\Dif F_{\alpha}(0,0,0)$ derived in Proposition 
\ref{prop:linpde},
\begin{align}
    \big\vert \rho^{2} \, \Dif F_{\alpha}(0,0,0) (u,a,b)
    \big\vert_{C^{0,\beta}_{\rho}(\Mapprox)} & \leq \big\vert \rho^{2}
    \cos(\theta_{\Mapprox}) \Delta_{\Mapprox} u \big\vert_{\cobr}
    \notag \\
    &\qquad + \big\vert \rho^{2} \sin(\theta_{\Mapprox}) \langle
    \vec{H}_{\Mapprox} , \nabla u \rangle_{\Mapprox} \big\vert
    _{C^{0,\beta}_{\rho}(\Mapprox)} \notag \\
    &\qquad + \vert a \vert \, \vert \rho^{2} \cos(\thapprox)
    \vert_{\cobr} - \vert b \vert \, \vert \rho^{2} \Delta_{\Mapprox}
    v_{e} \vert_{\cobr} \,
    \label{eqn:linupperbd}
\end{align}
Each of the four terms in the expression above will be estimated in
turn.  Begin with the first term:
\begin{align*}
    \big\vert \rho^{2} \cos(\theta_{\Mapprox}) \Delta_{\Mapprox} u
    \big\vert_{\cobr} &\leq \vert \rho^{2} \cos (\thapprox)
    \Delta_{\Mapprox} u \vert_{0,\Mapprox} + \vert \cos (\thapprox)
    \vert_{0, \Mapprox} \cdot [ \rho^{2+\beta} \Delta_{\Mapprox} u
    ]_{\beta, \Mapprox} \\
    &\qquad + [\cos(\thapprox) ]_{\beta, \Mapprox} \cdot \vert
    \rho^{\beta} \vert_{0,\Mapprox} \cdot \vert \rho^{2} \Delta
    _{\Mapprox} u \vert_{0, \Mapprox} \\
    &\leq \Big( 1 + [\cos(\thapprox) ]_{\beta, \Mapprox} \cdot \vert
    \rho^{\beta} \vert_{0, \Mapprox} \Big) \vert u \vert_{\ctbr}
\end{align*}
by definition of the $C^{2,\beta}_{\rho}$ norm and the bounds
\eqref{eqn:trigbounds} on the trigonometric functions of $\thapprox$. 
Now,
\begin{equation*}
    \big\vert \rho^{2} \cos(\theta_{\Mapprox}) \Delta_{\Mapprox} u
    \big\vert_{\cobr} \leq C (1 + \alpha^{1- \beta}) \vert u
    \vert_{\ctbr} \, ,
\end{equation*}
by the bounds on the function $\rho$ and on the H\"older norm of the
cosine term.  

The second term in equation \eqref{eqn:linupperbd} can be estimated in
a manner similar to that used for the first term, this time using the
bounds on the mean curvature from Proposition \ref{prop:meancurv}. 
One obtains
\begin{equation}
    \big\vert \rho^{2} \sin(\theta_{\Mapprox}) \langle
    \vec{H}_{\Mapprox} , \nabla u \rangle_{\Mapprox} \big\vert
    _{C^{0,\beta}_{\rho}(\Mapprox)} \leq C \alpha \, .
    \label{eqn:secondterm}
\end{equation}
Finally, deal with the remaining two terms in \eqref{eqn:linupperbd}. 
It is trivial to show that the third term is bounded above by a
constant independent of $\alpha$.  For the last term, recall that
$v_{e}$ is explicitly independent of $\alpha$ and is nonzero only near
the boundary of $\Mapprox$ where $\rho$ is independent of $\alpha$;
thus bounding this term above by a constant independent of $\alpha$ is
trivial as well.

The considerations of the previous two paragraphs lead to the
estimates
$$\big\vert \rho^{2} \, \Dif F_{\alpha}(0,0,0) \big(u,a,b \big)
\big\vert_{C^{0,k}_{\rho}(\Mapprox)} \leq C \Big( \vert u 
\vert_{\ctbr} + \vert a \vert + \vert b \vert \Big) \, .$$
Since the norm on the right hand side above is equivalent to the
$\Vert \cdot \Vert_{\ctbr \times \R^{2}}$ norm, the proof of the
proposition is complete.
\end{proof}

% Global elliptic estimate

The next theorem shows that an elliptic estimate for the operator
$\Delta_{\Mapprox}$ can be found using the weighted Schauder norms and
that the constant appearing there is independent (or nearly so) of
$\alpha$.

\begin{thm}[Elliptic Estimate] \label{thm:schauder}
    There is a constant $\Cell$ independent of $\alpha$ so that the
    elliptic estimate
    \begin{equation}
	\vert u \vert_{C^{2,\beta}_{\rho}(\Mapprox)} \leq \Cell \left(
	\eps^{-2\beta} \vert \rho^{2} \Delta_{\Mapprox} u \vert_{C^{0,
	\beta}_{\rho} (\Mapprox)} + \vert u \vert_{0,\Mapprox} \right)
        \label{eqn:schauder}
    \end{equation}
    holds for any $C^{2,\beta}$ function $u$ on $\Mapprox$ with
    $\alpha \leq \Rvol$ and satisfying the Neumann boundary condition
    $Z(u) = 0$ on $\partial \Mapprox$.
\end{thm}
    
\begin{proof}
    
The strategy for proving the elliptic estimate \eqref{eqn:schauder} is
to piece together local elliptic estimates, valid in coordinate charts
in which the metric coefficients are uniformly bounded.  These local
elliptic estimates can be phrased as follows and their proof can be
found in any standard textbook on elliptic theory, for instance
\cite{gt}.

\medskip
\noindent \scshape Fact: \upshape Suppose $P : C^{2,\beta} (\Omega)
\longrightarrow C^{0,\beta} (\Omega)$ is a second order elliptic
operator defined on a domain $\Omega$ contained in $\R^{n}$.  Then the
following estimates are valid for any $C^{2,\beta}$ function $u$ on
$\Omega$.  Here, $\theta \in (0,1)$ is arbitrary and
$C_{\mathit{loc}}$ is a constant which depends only on $\theta$,
$\beta$, the $C^{0,\beta}$ norm of the coefficients of $P$ and the
dimension $n$.  The norms used here are the local, scale invariant
Schauder norms introduced in \eqref{eqn:localnorm}.
\begin{enumerate}
    
    \item Let $x$ belong to the interior of $\Omega$.  Then,
    \begin{equation}
	\vert u \vert_{2,\beta, B_{\theta R}(x)}^{\ast} \leq
	C_{\mathit{loc}} \left( R^{2} \vert P u \vert_{0,\beta,
	B_{R}(x)} ^{\ast} + \vert u \vert_{0, B_{R}(x)} \right)
	\label{eqn:localest}
    \end{equation}
    where $B_{R}(x)$ is a ball centered at $x$ and contained in the
    interior of $\Omega$. 
    
    \item Let $x$ belong to $\partial \Omega$ and define $B_{R}^{+}(x)
    = B_{R}(x) \cap \Omega$ and $U_{R}(x) = B_{R}(x) \cap \partial
    \Omega$.  Then,
    \begin{equation}
	\vert u \vert_{2,\beta, B_{\theta R}^{+}(x)}^{\ast} \leq
	C_{\mathit{loc}} \left( R^{2} \vert P u \vert_{0,\beta,
	B_{R}^{+}(x)}^{\ast} + \vert u \vert_{0, B_{R}^{+}(x)} + R
	\vert Z (u) \vert_{1,\beta, U_{R}(x)}^{\ast} \right)
	\label{eqn:localbdest}
    \end{equation}
    where $Z$ is the conormal vector field of the boundary of $\Omega$. 
    
\end{enumerate}
\medskip

The proof of the global elliptic estimate for $\Delta_{\Mapprox}$ on
all of $\Mapprox$ will establish in two separate calculations that for
any $x \in \Mapprox$,
\begin{equation}
    \vert u(x) \vert + \rho(x) \Vert \nabla u (x) \Vert + \rho^{2}(x)
    \Vert \nabla^{2} u(x) \Vert \leq Q
    \label{eqn:ellstep1}
\end{equation}
and    
\begin{equation}
    \frac{\big\Vert \, \rho^{2 + \beta}(x) \nabla^{2}u(x) - \rho^{2 +
    \beta}(x') \nabla^{2}u(x') \, \big\Vert}{ \big( \mathit{dist}(x,x')
    \big)^{\beta}} \leq Q \, ,
    \label{eqn:ellstep2}
\end{equation}
where $Q$ refers to the right hand side of the inequality
\eqref{eqn:schauder} above.  If this is established for all $x \in
\Mapprox$ (or $x$ \emph{and} $x'$ in $\Mapprox$ in the second case),
then the theorem follows by taking the supremum.

Begin by choosing a point $x$ in the interior of $\Mapprox$.  For the
first calculation, let $s = \rho(x)$ and recall that $\Vert \nabla
\rho(y) \Vert \leq K \eps^{-\beta} \equiv \Crho$ for all $y \in
\Mapprox$.  Property 6 of the weight function shows that one can
assume, without loss of generality, that the constant $\Crho$ is such
that $\frac{s}{2\Crho}$ is less than $r(x)$, which is the radius of
uniformity of the metric coefficients at $x$.  Consequently, local
coordinates in which the coefficients of the metric are uniformly
bounded can be used within $B_{s/2\Crho}(x)$.  Furthermore, if $y$ is
any point in $B_{s/2\Crho}(x)$, then $\vert \rho(y) - \rho(x) \vert
\leq \Crho \vert x - y \vert$ and this implies
\begin{equation*}
    \vert \rho(y) \vert \geq \vert \rho(x) \vert - \Crho \vert x - y
    \vert \geq \frac{s}{2} \, .
\end{equation*}
By a similar argument, $\vert \rho(y) \vert \leq \frac{3s}{2}$ in
$B_{s/2\Crho}(x)$.  Now choose $\theta \in (0,1)$ and argue as
follows.  First, by simple algebra,
\begin{equation*}
    \vert u(x) \vert + s \Vert \nabla u(x) \Vert + s^{2} \Vert
    \nabla^{2} u(x) \Vert \leq \vert u(x) \vert + C \eps^{-2\beta}
    \left( \frac{\theta s}{2\Crho} \Vert \nabla u(x) \Vert + \left(
    \frac{\theta s}{2\Crho} \right)^{2} \Vert \nabla^{2} u(x) \Vert
    \right)
\end{equation*}
where $C$ depends only on $\theta$ and $K$.  But now, the
definition of the local $\vert \cdot \vert^{\ast}$ Schauder norm and
the local elliptic estimate in the ball $B_{s/2 \Crho}(x)$ can be applied
to give
\begin{align}
    \vert u(x) \vert + s \Vert \nabla u(x) \Vert + s^{2} \Vert
    \nabla^{2} u(x) \Vert \notag \\
    &\hspace{-15ex} \leq C \eps^{-2 \beta} \left( s^{2} \vert
    \Delta_{\Mapprox} u \vert_{0,B_{s/2\Crho}(x)} + s^{2+ \beta} [
    \Delta_{\Mapprox} u ]_{\beta, B_{s/2\Crho}(x)} \right) + \vert u
    \vert_{0,\Mapprox} \, ,
    \label{eqn:elleqn2}
\end{align}
where $C$ now depends on $C_{loc}$.  The first term in equation
\eqref{eqn:elleqn2} is easy to handle:
\begin{equation}
    \vert \rho^{2} \Delta_{\Mapprox} u \vert_{0, \Mapprox} \geq
    \sup_{y \in B_{s/2\Crho}(x)} \vert \rho^{2}(y) \Delta_{\Mapprox}
    u(y) \vert \geq \frac{s^{2}}{4} \vert \Delta_{\Mapprox} u \vert_{0,
    B_{s/2\Crho}(x)}
    \label{eqn:elleqn3}
\end{equation}
by the lower bounds on $\rho$.  For the second term of
\eqref{eqn:elleqn2}, choose $y$ and $y'$ in $B_{s/2\Crho}(x)$ and
estimate
\begin{align}
    \frac{\vert \rho^{2+\beta}(y) \Delta_{\Mapprox}u(y) -
    \rho^{2+\beta}(y') \Delta_{\Mapprox}u(y') \vert} {\vert y -y'
    \vert^{\beta}} & \notag \\
    &\hspace{-30ex} \geq \vert \rho^{2+\beta}(y) \vert \frac{\vert
    \Delta_{\Mapprox}u(y) - \Delta_{\Mapprox}u(y') \vert}{\vert y -y'
    \vert^{\beta}}  \notag \\
    &\hspace{-30ex} \qquad - (2+\beta) \, \vert \rho^{1+\beta}
    \vert_{0, B_{s/2 \Crho}(x)} \, \Vert \nabla \rho \Vert_{0,
    B_{s/2\Crho}(x)} \, \vert y - y' \vert^{1-\beta} \, \vert
    \Delta_{\Mapprox} u(y') \vert \notag \\
    &\hspace{-30ex} \geq C \left( s^{2+\beta} \frac{ \vert
    \Delta_{\Mapprox} u(y) - \Delta_{\Mapprox} u(y') \vert} {\vert y
    -y' \vert^{\beta}} - \Crho^{\beta} \frac{s^{2}}{4} \vert
    \Delta_{\Mapprox} u \vert_{0,B_{s/2\Crho}(x)} \right) \notag \\
    &\hspace{-30ex} \geq C \left( s^{2+\beta} \frac{ \vert
    \Delta_{\Mapprox}u(y) - \Delta_{\Mapprox}u(y') \vert} {\vert y -y'
    \vert^{\beta}} - \eps^{-2\beta} \vert \rho^{2}
    \Delta_{\Mapprox} u \vert_{0,B_{s/2\Crho}(x)} \right)
    \label{eqn:elleqn5}
\end{align}
using the bounds on $\rho$ and $\Vert \nabla \rho \Vert$ as well as
the result in equation \eqref{eqn:elleqn3}.  If the supremum of
\eqref{eqn:elleqn5} over all $y$ and $y'$ in the ball
$B_{s/2\Crho}(x)$ is taken, then the inequality
\begin{equation}
    s^{2+\beta} [\Delta_{\Mapprox} u ]_{\beta, B_{s/2\Crho}(x)} \leq C
    \eps^{-2\beta} \vert \rho^{2} \Delta_{\Mapprox} u \vert_{C^{0,
    \beta}_{\rho}(\Mapprox)}
    \label{eqn:elleqn4}
\end{equation}
follows, again by applying equation \eqref{eqn:elleqn3}.  Substituting
the inequality \eqref{eqn:elleqn4} along with the previous inequality
\eqref{eqn:elleqn3} into equation \eqref{eqn:elleqn2} yields the first
estimate \eqref{eqn:ellstep1}.

For the sake of brevity, the calculation of the second inequality in
\eqref{eqn:ellstep2} will be omitted because it is essentially the
same as the previous calculation.  Nevertheless, all the details can
be found in \cite{me}.  Finally, the case of $x \in \partial \Mapprox$
follows trivially from the local elliptic boundary estimate because
the boundary of $\Mapprox$ is independent of $\alpha$.
\end{proof}    

\section{Analysis of the Linearized Operator}

\subsection{Outline}

% Outline of injectivity estimate

In order to invoke the Inverse Function Theorem to produce solutions
of the equation $F_{\alpha}(H, \theta, b) = 0$, the linearized
operator $\Dif F_{\alpha}(0,0,0)$ must be a bijection satisfying the
estimates required by the Inverse Function Theorem.  The injectivity
of the linearized operator will be established by producing a lower
bound of the form $\vert \rho^{2} \Dif F_{\alpha}(0,0,0) ( u,a,b )
\vert_{C^{0, \beta}_{\rho} (\Mapprox)} \geq C_{L} ( \alpha ) \Vert (u,
a, b) \Vert _{C^{2,\beta}_{\rho} (\Mapprox) \times \R^{2}}$ for any $u
\in \mathcal{B}_{1,\alpha}$ and $(a,b) \in \R^{2}$.  Let 
\begin{gather*}
    P_{\alpha} = \cos(\theta_{\Mapprox}) \, \Delta_{\Mapprox} +
    \sin(\theta_{\Mapprox}) \, \vec{H}_{\Mapprox} \!  \cdot \nabla 
    \, ,\\
    \psizero = \cos(\thapprox) \, , \\
    \psione = - \Delta_{\Mapprox} v_{e} \, .
\end{gather*}
The desired estimate is now equivalent to:
\begin{equation}
	\Big\vert \rho^{2} P_{\alpha} u - a \psizero + b \psione
	\Big\vert_{C^{0, \beta}_{\rho} (\Mapprox)} \geq C_{L}( \alpha)
	\left( \vert u \vert_{C^{2,\beta}_{\rho} (\Mapprox)}^{2} +
	a^{2} + b^{2} \right)^{1/2} \, .
	\label{eqn:injbd}
\end{equation}
Such an estimate will be developed in Sections 5.2 through 5.6 and the
precise dependence of $C_{L}$ on $\alpha$ will of course be determined
as well.  It will be shown that the choice of parametrization made in
the previous Section is enough to ensure that $C_{L}(\alpha)$ is
nearly independent of $\alpha$.  The surjectivity of the linearized
operator will follow in a more or less straightforward manner from the
results of the analysis leading up to the injectivity estimate and
will be presented in Section 5.7.

The lower bound \eqref{eqn:injbd} comes from combining four separate
results, and it is worthwhile to indicate in general terms how this
will be done before proceeding with the details.  The starting point
is an estimate on the second Neumann eigenvalue of the Laplacian on
$\Mapprox$.  This is the key estimate of this paper.  It is purely
global in nature and is of the form
\begin{equation} 
    \Vert \Delta_{\Mapprox} u \Vert_{L^{2}(\Mapprox)} \geq C \Vert u
    \Vert_{L^{2}(\Mapprox)}
    \label{eqn:eigenest}
\end{equation}
for functions $u$ perpendicular to $1$ and to $S_{\alpha}$, the first
eigenfunction of the Laplacian on $\Mapprox$.  Here, $C$ is
independent of $\alpha$.

The second step is to deduce an estimate similar to the one above but
in the weighted $C^{k,\beta}_{\rho}$ Schauder norms; in other words,
one of the form
\begin{equation}
    \vert \rho^{2} \Delta_{\Mapprox} u \vert_{C^{0,\beta}
    _{\rho}(\Mapprox)} \geq C \eps^{2\beta} \vert u
    \vert_{C^{2,\beta}_{\rho}(\Mapprox)} \, .
    \label{eqn:continjest}
\end{equation}
This estimate will be found by combining the elliptic estimate for
$\Delta_{\Mapprox}$ derived in Section 4.2 with the De Giorgi Nash
estimate for $\Delta_{\Mapprox} u$.

The third step is to deduce the desired injectivity estimate
$P_{\alpha}$.  This estimate, namely
\begin{equation}
    \big\vert \rho^{2} \big( \cos(\thapprox) \, \Delta_{\Mapprox} u +
    \sin(\theta_{\Mapprox}) \, \vec{H} _{\Mapprox} \!  \cdot \nabla u
    \big) \big\vert_{C^{0,\beta} _{\rho}(\Mapprox)} \geq C
    \eps^{2\beta} \vert u \vert_{C^{2,\beta}_{\rho}(\Mapprox)} \, ,
    \label{eqn:restrictedest}
\end{equation}
holds because the mean curvature trigonometric terms can be controlled
when $\alpha$ is small.

The final component of the injectivity estimate is to incorporate the
$\psizero$ and $\psione$ factors into the estimate found in the
preceding step.  This is a consequence of the fact that functions in
$\mathcal{B}_{1, \alpha}$ are orthogonal to $1$ and $S_{\alpha}$.

\subsection{Second Eigenvalue Estimate for $\Delta_{\Mapprox}$}

% Second Neumann eigenvalue

The fundamental fact that will guarantee the injectivity of the
linearized operator $\Dif F_{\alpha}(0,0,0)$ is that the \emph{second}
eigenvalue $\nu_{2}$ of $\Delta_{\Mapprox}$ is bounded below
independently of $\alpha$.  This bound will be derived by piecing
together two $\alpha$-independent inequalities that already hold on
$\Mapprox$.  The first of these is essentially the Poincar\'e
inequality of the original special Lagrangian submanifolds $M_{1}$ and
$M_{2}$ comprising $M$.  

\begin{prop} \label{prop:poincare}
    Suppose $u$ belongs to $H^{1}(\Mapprox)$, satisfies
    $\int_{\Mapprox} u = 0$, and has support contained in only one of
    the factors $M_{1}'$ or $M_{2}'$.  Then the function $u$ satisfies
    the inequality
    $$\int_{\Mapprox} \Vert \nabla u \Vert^{2} \geq E_{0} \int_{\Mapprox} 
    u^{2}$$
    where $E_{0}$ is a constant independent of $\alpha$,
\end{prop}

\begin{proof}
    
Let $E_{0} = \min \{ \lambda_{1}, \lambda_{2} \}$, where each
$\lambda_{i}$ is the first Neumann eigenvalue of $M_{i}$.  The result
of the proposition is now an elementary consequence of the Poincar\'e
inequalities $\int_{M_{i}} \Vert \nabla u \Vert^{2} \geq \lambda_{i}
\int_{M_{i}} u^{2}$ and the fact that $M_{i}' \subset M_{i}$.
\end{proof}
    
% Michael Simon inequality

\noindent The second inequality that is valid on $\Mapprox$ holds
independently of $\alpha$ even in the neck region.

% M-S Poincar\'e estimate

\begin{prop} \label{prop:mspoincare}
    There is a number $\Rms$ with $0 < \Rms \leq \Rvol$ and a
    geometric constant $\Cms$ independent of $\alpha$ so that if
    $\alpha \leq \Rms$ then the inequality
    \begin{equation}
	\int_{\Mapprox} \vert u \vert^{2} \leq \Cms \int_{\Mapprox}
	\Vert \nabla u \Vert^{2} 
	\label{eqn:mspoincare}
    \end{equation}
    holds for any $H^{1}$ function $u$ that vanishes in a
    neighbourhood of $\partial \Mapprox$.
\end{prop}

\begin{proof}
    
The Michael-Simon inequality \cite{simon} for a submanifold $M$ with
boundary in $\Rtn$ states
\begin{equation} 
    \left( \int_{M} \vert v \vert^{\tfrac{n}{n-1}}
    \right)^{\tfrac{n-1}{n}} \leq C \int_{M} \Big( \Vert
    \nabla v \Vert + \Vert \vec{H}_{M} \Vert \cdot \vert v \vert \Big)
    \label{eqn:simon}
\end{equation}
for any function $v \in H^{1}(M)$ which vanishes in some neighbourhood
of the boundary $\partial M$.  Here, $\vec{H}_{M}$ is the mean
curvature vector of $M$ and $C$ is a constant depending only on $n$.

According to Proposition \ref{prop:meancurv}, the mean curvature
vector of $\Mapprox$ is pointwise bounded everywhere in the transition
region $T_{1} \cup T_{2}$ and vanishes everywhere else.  Thus integral
norms of the magnitude of $\vec{H}_{\Mapprox}$ can be made as small as
desired by choosing the parameter $\alpha$ small enough.  This makes
it possible to absorb the mean curvature term in \eqref{eqn:simon}
into the left hand side, thereby producing the inequality
\eqref{eqn:mspoincare}.  This calculation will be omitted because it
is a straightforward application of the H\"older inequality upon
substituting $v = u^{2(n-1)/n}$.
\end{proof}

The desired estimate on $\nu_{2}$ now comes from combining the
Poincar\'e inequalities on the separate pieces $M_{i}'$ with the
Michael-Simon result in the neck region $T_{1} \cup N' \cup T_{2}$.

% The second eigenvalue estimate

\begin{thm}[Second Eigenvalue Estimate] \label{thm:2ndeigen}
    If $\alpha \leq \Rms$, then the second Neumann eigenvalue
    $\nu_{2}$ of the Laplacian on $\Mapprox$ is bounded below by a
    geometric constant $\Cnutwo$ independent of $\alpha$. 
    Consequently, the inequality
    $$\Vert \Delta_{\Mapprox} u \Vert_{L^{2}(\Mapprox)} \geq \Cnutwo
    \Vert u \Vert_{L^{2}(\Mapprox)}$$
    holds for all $u \in \Bonealpha$.
\end{thm}

\begin{proof}
    
The lower bound $\Cnutwo$ will be deduced using the max-min 
characterization of the second eigenvalue, namely that
$$\nu_{2} \geq \inf_{\genfrac{}{}{0pt}{}{u \in H^{1}(\Mapprox)}{u \neq
0}} \left\{ \frac{\int_{\Mapprox} \Vert \nabla u \Vert^{2}} {\int_{
\Mapprox} u^{2}} : \int_{\Mapprox} u \cdot \phi_{1} = \int_{\Mapprox}
u \cdot \phi_{2} = 0 \right\} $$
for any specific choice of $\phi_{1}$ and $\phi_{2}$.  Begin by 
making such a choice.

Recall that the submanifold $\Mapprox$ consists of the union of the
exterior components $M_{1}'$ and $M_{2}'$ of $M$, along with the neck
region $T_{1} \cup N' \cup T_{2}$.  Denote this latter region simply
by $N$.  Define the functions $\phi_{i} : \Mapprox \longrightarrow \R$
by the prescription
\begin{equation}
    \begin{gathered}
	\phi_{1} = \chi_{\! \genfrac{}{}{0pt}{2}{}{M_{1}'}} +
	\chi_{\! \genfrac{}{}{0pt}{2}{}{M_{2}'}} \\
        \phi_{2} = \chi_{\! \genfrac{}{}{0pt}{2}{}{M_{1}'}} -
	\chi_{\! \genfrac{}{}{0pt}{2}{}{M_{2}'}} 
    \end{gathered}
    \label{eqn:charfns}
\end{equation}
where $\chi_{\!  \genfrac{}{}{0pt}{3}{}{U}}$ refers to the
characteristic function of the subset $U$ of $\Mapprox$.

Let $u$ be any $C^{1}$ function on $\Mapprox$.  Suppose that it has
$L^{2}$ norm equal to one and satisfies the orthogonality conditions
$\int_{\Mapprox} \phi_{1} \cdot u = 0$ and $\int_{\Mapprox} \phi_{2}
\cdot u = 0$.  These two conditions give
\begin{equation}
    \int_{M_{1}'} u = 0 \qquad \mbox{and} \qquad \int_{M_{2}'} u = 0 
    \, .
    \label{eqn:conditions}
\end{equation}

Choose any $\delta \in (0,1)$ and consider the two separate cases:
either $\int_{M_{1}' \cup M_{2}'} u^{2} \geq 1 - \delta$ and thus
$\int_{N} u^{2} \leq \delta$; or else $\int_{M_{1}' \cup M_{2}'} u^{2}
\leq 1 - \delta$ and thus $\int_{N} u^{2} \geq \delta$.  It will be
possible to show that there is always a choice of $\delta \in (0,1)$
so that in \emph{both} cases above, the inequality $\int_{\Mapprox}
\Vert \nabla u \Vert^{2} \geq \Cnutwo$ holds for some positive
geometric constant $\Cnutwo$ independent of $\alpha$.

\medskip 
\noindent \scshape Case 1: \upshape $\int_{M_{1}' \cup
M_{2}'} u^{2} \geq 1 - \delta$ and $\int_{N} u^{2} \leq \delta$.
\medskip

\noindent Calculate as follows:
\begin{align*}
    \int_{\Mapprox} \Vert \nabla u \Vert^{2} &= \int_{M_{1}'} \Vert
    \nabla u \Vert^{2} + \int_{N} \Vert \nabla u \Vert^{2} +
    \int_{M_{2}'} \Vert \nabla u \Vert^{2} \\
    &\geq \int_{M_{1}} \Vert \nabla \bar{u} \Vert^{2} + \int_{M_{2}}
    \Vert \nabla \bar{u} \Vert^{2}
\end{align*}
where $\bar{u}$ is the $L^{2}$ function on $\Mapprox$ which is equal
to $u$ on $M_{1}' \cup M_{2}'$ and equal to zero in $N$.  The
Poincar\'e inequality of Proposition \eqref{prop:poincare} for
functions supported away from the neck region can be applied (since
$H^{1}(\Mapprox)$ is dense in $L^{2}(\Mapprox)$) to give
\begin{align}
    \int_{\Mapprox} \Vert \nabla u \Vert^{2} &\geq E_{0}
    \int_{\Mapprox} \bar{u}^{2} \notag \\
    &= E_{0} \int_{M_{1}' \cup M_{2}'} u^{2} \notag \\
    &\geq E_{0} (1-\delta) \, .
    \label{eqn:nucase1}
\end{align}

\medskip 
\noindent \scshape Case 2: \upshape $\int_{M_{1}' \cup
M_{2}'} u^{2} \leq 1 - \delta$ and $\int_{N} u^{2} \geq \delta$. 
\medskip

Choose a function $\eta : \Mapprox \longrightarrow \R$ with the
following properties.  Let $\eta$ satisfy $0 \leq \eta \leq 1$, be
equal to one in $N$ and vanish in a neighbourhood of the boundary of
$\Mapprox$.  Moreover, choose $\eta$ such that $\Vert \nabla \eta
\Vert$ is bounded above by a constant $K$ that is independent of
$\alpha$.  Write $u = \eta u + (1-\eta) u$.  The function $\eta u$ is
thus a member of $H^{1}(\Mapprox)$ and vanishes in a neighbourhood of
the boundary.  Consequently, the Poincar\'e inequality of Proposition
\eqref{prop:mspoincare} is valid for $\eta u$, and this yields the
inequalities
\begin{align}
    \int_{\Mapprox} \Vert \nabla u \Vert^{2} &= \int_{\Mapprox}
    \left\Vert \nabla \big( \eta u + (1-\eta) u \big) \right\Vert^{2} 
    \notag \\
    &= \int_{\Mapprox} \Vert \nabla \eta u \Vert^{2} + 2
    \int_{\Mapprox} \left\langle \nabla \eta u , \nabla (1-\eta)u
    \right\rangle + \int_{\Mapprox} \Vert \nabla (1-\eta) u \Vert^{2}
    \notag \\
    &\geq \Cms \int_{N} u^{2} + 2 \int_{\Mapprox} (1 - 2\eta) u \nabla
    u \cdot \nabla \eta - 2 K^{2} \int_{M_{1}' \cup M_{2}'} u^{2}
    \notag \\
    &\geq \Cms \delta + 2 \int_{\Mapprox} (1 - 2\eta) u \nabla
    u \cdot \nabla \eta - 2 K^{2} (1-\delta) \label{eqn:nustep1}
\end{align}
by the Michael-Simon Poincar\'e inequality and the fact that $\nabla
\eta$ is equal to zero in $N$.  It remains only to deal with the cross
term.  By the Cauchy-Schwarz inequality,
\begin{align*}
    2 \int_{\Mapprox} (1 - 2\eta) u \nabla u \cdot \nabla \eta &\geq - 
    2 \int_{\Mapprox} \vert u \vert \, \vert 1 - 2 \eta \vert \, 
    \vert \nabla u \cdot \nabla \eta \vert \\
    &\geq - 2 \int_{\Mapprox} \vert u \vert \, \Vert \nabla u \Vert \,
    \Vert \nabla \eta \Vert \, .
\end{align*}
The fact that $\vert 1 - 2 \eta \vert \leq 1$ has been used here.  Now
apply the H\"older inequality and the Schwarz inequality $2 a b \leq
a^{2} + b^{2}$ to the integral on the right hand side above to obtain:
\begin{align}
    2 \int_{\Mapprox} (1 - 2\eta) u \nabla u \cdot \nabla \eta &\geq - 
    2 \left( \int_{\Mapprox} \Vert \nabla u \Vert^{2} 
    \right)^{\frac{1}{2}} \left( \int_{\Mapprox} \vert u \vert^{2} 
    \, \Vert \nabla \eta \Vert^{2} \right)^{\frac{1}{2}} \notag \\
    &\geq - \int_{\Mapprox} \Vert \nabla u \Vert^{2} - K^{2}
    \int_{M_{1}' \cup M_{2}'} \vert u \vert^{2} \notag \\
    &\geq - \int_{\Mapprox} \Vert \nabla u \Vert^{2} - 
    K^{2}(1-\delta) \, .\label{eqn:nustep2}
\end{align}
Use the inequality \eqref{eqn:nustep2} in the estimate 
\eqref{eqn:nustep1} above and combine the $\int \Vert \nabla u 
\Vert^{2}$ terms on the left hand side.  This results in:
\begin{equation}
    \int_{\Mapprox} \Vert \nabla u \Vert^{2} \geq \frac{\Cms 
    \delta}{2} - \frac{3 K^{2} (1-\delta)}{2} \, .
    \label{eqn:nucase2}
\end{equation}

In order to complete the proof of the theorem, the quantity $\delta$
must be chosen between zero and one to make \emph{both} bounds
\eqref{eqn:nucase1} and \eqref{eqn:nucase2} positive.  This can be
accomplished by choosing $\delta$ near enough to 1 to make the
negative term in \eqref{eqn:nucase2} strictly smaller than the
positive term.
\end{proof}

\subsection{Injectivity Estimate for $\Delta_{\Mapprox}$ in the
Weighted Schauder Norms}

The next step in the proof of the injectivity estimate is to combine
the eigenvalue estimate of the previous section with the elliptic
estimate for $\Delta_{\Mapprox}$ and derive equation
\eqref{eqn:continjest} in the Schauder $C^{k,\beta}_{\rho}$ norms on
$\Mapprox$.  The tool which allows the two estimates to be combined is
the De Giorgi Nash inequality.  This inequality will make many
appearances in the sequel.

% DGN Inequality

\medskip
\noindent \scshape Fact: \upshape Let $u \in H^{1}(\Mapprox)$ be a
weak solution of the equation $\Delta_{\Mapprox} u = f$ where $f \in
L^{q/2}(\Mapprox)$ for some $q > n$.  Furthermore, suppose $u$ 
satisfies Neumann boundary conditions.  Then there is a constant
$C_{DGN}$ depending only on $n$ and $\delta$ so that for every $p \geq
1$,
\begin{equation}
    \vert u \vert_{0, \Omega'} \leq C_{DGN} \left( \big(
    \mathit{Vol}(\Omega) \big)^{-1/p} \Vert u \Vert_{L^{p}(\Omega)} +
    \big( \mathit{Vol}(\Omega) \big)^{2/n - 2/q} \Vert f
    \Vert_{L^{q/2}(\Omega)} \right) \, ,
    \label{eqn:dgn}
\end{equation}
where $\Omega$ is any subset of $\Mapprox$ and $\Omega'$ is any subset
of the interior of $\Omega$.  The reader should consult \cite{schoen},
\cite{simon2} or \cite{taylor} for more details.  \medskip

The constant in the De Giorgi Nash inequality is independent of
$\alpha$ because it depends on the metric of $\Mapprox$ \emph{only}
through the constant appearing in the Michael-Simon inequality
\eqref{eqn:simon}, and this quantity is independent of $\alpha$.  The
estimate thus provides the essential link between the
$\alpha$-independent $L^{2}$ injectivity estimate for $\Delta_
{\Mapprox}$ derived in the previous section and the required
$C^{k,\beta}$ estimates.

% C^{0} injectivity bound

\begin{prop} \label{prop:c0injest}
    If $\alpha \leq \Rnew$, then there exists a constant $\Cczero$
    independent of $\alpha$ so that
    $$\vert u \vert_{0,\Mapprox} \leq \Cczero \eps^{- 2\beta} \vert
    \rho^{2} \Delta_{\Mapprox} u \vert_{0,\Mapprox}$$
    for any function $u \in \mathcal{B}_{1,\alpha}$.
\end{prop}

\begin{proof}
    
Suppose that the proposition is false.  Then there exists a sequence
of parameters $\alpha_{j} \rightarrow 0$ (thus also a sequence of
weights $\rho_{j}$ and scales $\eps_{j}$) and a sequence of functions
$u_{j} \in \mathcal{B}_{1,\alpha_{j}}$, normalized so that $\vert
u_{j} \vert_{0, \Mapproxj} = 1$, that satisfy the following
inequality:
\begin{equation}
    \eps_{j}^{-2\beta} \left\vert \rho_{j}^{2} \Delta_{\Mapproxj}
    u_{j} \right\vert_{0, \Mapproxj} \leq \frac{1}{j} \, ,
    \label{eqn:assumption}
\end{equation}
for each $j$. Let $x_{j} \in \Mapproxj$ be a point where $u_{j}(x_{j}) = 1$
and consider a large ball $B_{R}(x_{j})$ centered at $x_{j}$ --- any
$R$ independent of $j$ will do.  By the De Giorgi Nash Estimate with
$p=2$ and some $q > n$ that will be specified later,
\begin{equation} 
    1 = \big\vert u_{j} \big\vert_{0, B_{\theta R}(x_{j})} 
    \leq C \left( \Vert u_{j} \Vert_{L^{2}(\Mapproxj)} + \Vert
    \Delta_{\Mapproxj} u_{j} \Vert_{L^{q/2}(\Mapproxj)} \right) \, ,
    \label{eqn:fuckyoualready}
\end{equation}
where $C$ depends on $C_{DGN}$ and $R$.  Next, apply the second 
eigenvalue estimate along with  Green's identity to obtain
\begin{equation}
    \int_{\Mapproxj} u_{j}^{2} \leq C \int_{\Mapproxj} \vert u_{j} 
    \cdot \Delta_{\Mapproxj} u_{j} \vert \leq C \Vert 
    \Delta_{\Mapproxj} u_{j} \Vert_{L^{1}(\Mapproxj)} \, ,
    \label{eqn:gotyou}
\end{equation}
using the normalization of $u_{j}$.  Substituting the result of
equation \eqref{eqn:gotyou} into \eqref{eqn:fuckyoualready} and
applying H\"older's inequality yields
\begin{equation} 
    1 \leq C \Vert \Delta_{\Mapproxj} u_{j} \Vert_{L^{q/2}
    (\Mapproxj)} \, .
    \label{eqn:c0est1}
\end{equation}
This will lead to a contradiction if it can be shown that the right
hand side of inequality \eqref{eqn:c0est1} goes to zero as $j
\rightarrow \infty$.  Note that in dimension $n=3$, the fraction $q/2$
can be less than 2.  Thus using the straightforward estimate $\Vert
u_{j} \Vert _{L^{2} (\Mapproxj)} \leq C \Vert \Delta_{\Mapproxj} u_{j}
\Vert _{L^{2}(\Mapproxj)}$ in \eqref{eqn:fuckyoualready} would not
lead to the estimate \eqref{eqn:c0est1}.  But the $L^{q/2}$ norm will
turn out to be a crucial component of the following argument.

To proceed, calculate the $L^{q/2}$ norm of $\Delta_{\Mapproxj}
u_{j}$.
\begin{align*}
    \Vert \Delta_{\Mapproxj} u_{j} \Vert_{L^{q/2}(\Mapproxj)} &=
    \left( \int_{\Mapproxj} \vert \Delta_{\Mapproxj} u_{j} \vert^{q/2}
    \right)^{2/q} \\
    &\leq \Rweight^{-2 \beta} \left( \int_{\Mapproxj} \rho_{j}^{-q(1 -
    \beta) } \vert \eps_{j}^{-2\beta} \rho_{j}^{2} \Delta_{\Mapproxj}
    u_{j} \vert^{q/2} \right)^{2/q} \\
    &\leq \frac{C}{j} \Vert \rho_{j}^{-1} \Vert_{L^{q(1-\beta)}
    (\Mapproxj)}^{2(1-\beta) }\, ,
\end{align*}
using the fact that each $\rho_{j}$ is everywhere bounded below by
$\Rweight \eps_{j}$.  By Property 5 of the weight function, $\big\Vert
\rho_{j}^{-1} \big\Vert _{L^{q \left(1- \beta \right)}(\Mapproxj)}$ is
bounded above independently of $j$ whenever $q < \frac{n}{ (1 -
\beta)}$.  Since $n < \frac{n}{1-\beta}$, there is such a choice of
$q$ compatible with the requirement $q>n$ that can be made in equation
\eqref{eqn:c0est1}.  Therefore, the quantity in the right hand side of
\eqref{eqn:c0est1} can be made to approach zero.
\end{proof}

Given the result of the previous proposition, it is now an easy matter
to transform the elliptic estimate of Theorem \ref{thm:schauder} into
a true injectivity estimate for the operator $\Delta_{\Mapprox}$ on
the space of functions $\mathcal{B}_{1, \alpha}$.

% Injectivity in the weighted norms

\begin{prop} \label{prop:Deltainj}
    Suppose $\alpha \leq \Rnew$.  Then there is a constant $\Cinj$
    independent of $\alpha$ so that the estimate
    \begin{equation}
	\vert u \vert_{C^{2,\beta}_{\rho}(\Mapprox)} \leq \Cinj \,
	\eps^{- 2\beta} \vert \rho^{2} \Delta_{\Mapprox} u
	\vert_{C^{0,\beta}
	_{\rho}(\Mapprox)}
        \label{eqn:Deltainj}
    \end{equation}
    holds for any function $u \in \Bonealpha$.
\end{prop}

\begin{proof}

The global elliptic estimate of Theorem \ref{thm:schauder} implies
that the function $u$ satisfies:
$$\vert u \vert_{C^{2,\beta}_{\rho}(\Mapprox)} \leq \Cell \left(
\eps^{-2 \beta} \vert \rho^{2} \Delta_{\Mapprox} u \vert_{C^{0,
\beta}_{\rho} (\Mapprox)} + \vert u \vert_{0,\Mapprox} \right) \, .$$
The $\vert u \vert_{0,\Mapprox}$ term can be controlled as in
Proposition \ref{prop:c0injest} and yields the desired estimate.
\end{proof}

\subsection{Injectivity Estimate for $P_{\alpha}$ in the Weighted
Schauder Norms}

The key to finding the injectivity estimate \eqref{eqn:restrictedest}
for the operator $P_{\alpha}$ is that $\cos(\thapprox) \approx 1$ and
$\sin (\thapprox) \approx 0$ when $\alpha$ is sufficiently small.

\begin{prop} \label{restrinjest}
    There is a number $\Rrest$ with $0 < \Rrest \leq \Rnew$ and a
    constant $\Crest > 0$ independent of $\alpha$ so that if $\alpha
    \leq \Rrest$, then
    \begin{equation}
	\vert u \vert_{\ctbr} \leq \Crest \eps^{-2 \beta} \vert
	\rho^{2} P_{\alpha} u\vert_{\cobr}
	\label{eqn:restrinjest}
    \end{equation}
    for any $u \in \mathcal{B}_{1,\alpha}$.  Here, $\thapprox$ is the
    angle function of $\Mapprox$ and $\meancurv$ is its mean curvature
    vector.
\end{prop}

\begin{proof}
    
Start with the right hand side of equation \eqref{eqn:restrinjest}.
\begin{equation}
    \vert \rho^{2} P_{\alpha} u \vert \geq \big\vert \rho^{2}
    \cos(\thapprox) \, \Delta _{\Mapprox} u \big\vert_{\cobr} -
    \big\vert \rho^{2} \sin( \thapprox) \, \langle \meancurv , \nabla
    u \rangle \big\vert_{\cobr}
    \label{eqn:restrone}
\end{equation}
The first term of \eqref{eqn:restrone} is easy to deal with by
applying the injectivity estimate from Proposition \ref{prop:Deltainj}
and the properties of the trigonometric functions of
$\theta_{\Mapprox}$ listed in equations \eqref{eqn:trigbounds}:
\begin{align}
    \big\vert \rho^{2} \cos(\thapprox) \, \Delta_{\Mapprox} u
    \big\vert_{\cobr} &= \vert \rho^{2} \cos(\thapprox) \,
    \Delta_{\Mapprox} u \vert_{0,\Mapprox} + [ \rho^{2+\beta}
    \cos(\thapprox) \, \Delta_{\Mapprox} u ]_{\beta, \Mapprox} \notag \\
    &\hspace{-5ex} \geq \vert \rho^{2} \Delta_{\Mapprox} \vert_{0,
    \Mapprox} \cdot \vert \cos(\thapprox) \vert_{0,\Mapprox} + [
    \rho^{2+\beta} \Delta_{\Mapprox} u ]_{\beta,\Mapprox} \cdot \vert
    \cos(\thapprox) \vert_{0,\Mapprox} \notag \\
    &\hspace{-5ex} \qquad - \vert \rho^{2} \Delta_{\Mapprox} u
    \vert_{0,\Mapprox} \cdot \vert \rho^{\beta} \vert_{0,\Mapprox}
    \cdot [ \cos(\thapprox) ]_{\beta,\Mapprox} \notag \\
    &\hspace{-5ex} \geq C \big( 1 - \alpha^{1-\beta} \big) \vert
    \rho^{2} \Delta_{\Mapprox} u \vert_{\cobr} \notag \\
    &\hspace{-5ex} \geq C \eps^{2\beta} \vert u \vert_{\ctbr} \, ,
    \label{eqn:restrtwo}
\end{align}
for some constant $C$ independent of $\alpha$ when $\alpha$ is
sufficiently small.  The second term of \eqref{eqn:restrone} can be
estimated from below using inequality \eqref{eqn:secondterm} from
Proposition \ref{prop:linbounded} and results in
\begin{equation*}
    \big\vert \rho^{2} \cos(\thapprox) \, \Delta_{\Mapprox} u +
    \sin(\thapprox) \, \meancurv \cdot \nabla u \big\vert_{\cobr} \geq
    C \eps^{2 \beta} \big(1 - \alpha \eps^{-2\beta} \big) \vert u
    \vert_{\ctbr} \, .
\end{equation*}
This gives the desired estimate when $\alpha$ is sufficiently
small.
\end{proof}

\subsection{Properties of the First Eigenfunction}

% First eigenfunction

In order to complete the injectivity component of the proof of the
Main Theorem, it remains only to extend the injectivity bounds derived
in the previous section to the full linearized deformation operator
$(u,a,b) \mapsto P_{\alpha}u + a \psizero + b \psione$.  Recall that
the functions $u$ are taken in the Banach space $\Bonealpha$ and are
thus $L^{2}$-orthogonal to the span of $1$ and $S_{\alpha}$, where
$S_{\alpha}$ is the first non-trivial eigenfunction of
$\Delta_{\Mapprox}$.  It will be shown that the injectivity estimate
for the full linearized operator is equivalent to finding a lower
bound, independent of $\alpha$, on the $L^{2}$ inner products
\begin{equation}
    \int_{\Mapprox} \psi_{i,\alpha} \qquad \mbox{and} \qquad
    \int_{\Mapprox} \psi_{i,\alpha} \cdot S_{\alpha}
    \label{eqn:ip}
\end{equation}
for $i= 1$ and $2$.  Unfortunately, too little is known about the
function $S_{\alpha}$ to accomplish this.  The purpose of this section
is to derive the necessary properties of $S_{\alpha}$.  Begin with the
following lemma.

\begin{lemma}
    \label{lemma:eigensupbd}
    There is a number $\Rinj$ with $0 < \Rinj \leq \Rrest$ and a
    constant $C> 0$ independent of $\alpha$ so that if $\alpha \leq
    \Rinj$, then the function $S_{\alpha}$ satisfies $\vert S_{\alpha}
    \vert_{0,\Mapprox} \leq C$.
\end{lemma}

\begin{proof}

Let $x \in \Mapprox$ be a point where $S_{\alpha}$ is maximal and let
$B_{R}(x)$ be any ball about $x$ with radius independent of $\alpha$. 
Apply the De Giorgi Nash inequality with $p=2$ and any $q > n$ to
$S_{\alpha}$ to $S_{\alpha}$:
\begin{align}
    \vert S_{\alpha} \vert_{0,\Mapprox} &= \vert S_{\alpha} \vert_{0,
    B_{R}(x)} \notag \\
    &\leq C \left( \Vert \Delta_{\Mapprox} S_{\alpha} \Vert_{L^{q/2}
    (\Mapprox)} + \Vert S_{\alpha} \Vert_{L^{2}(\Mapprox)} \right)
    \notag \\
    &\leq C \left( \nu_{1} \vert S_{\alpha} \vert_{0, \Mapprox} + 1 
    \right) \, .
    \label{eqn:esbd1}
\end{align}
Recall that $\nu_{1} \leq C \alpha^{n-2}$; thus if $\alpha$ is
sufficiently small, the $\vert S_{\alpha} \vert_{0,\Mapprox}$ terms
can be grouped together in \eqref{eqn:esbd1}, and lead to the desired
estimate.
\end{proof}

The remaining properties of $S_{\alpha}$ will be derived by first
approximating this function and then deriving these properties from
the approximation.  As in Proposition \ref{prop:firsteigen} of Section
3.4, let $\sigma$ be a function which is equal to 1 in $M_{1}'$, equal
to -1 in $M_{2}'$, and which interpolates smoothly between 1 and -1 in
the neck region of $\Mapprox$.  Furthermore, choose $\sigma$ so that
the derivative bounds $ \Vert \nabla \sigma \Vert + \alpha \Vert
\nabla^{2} \sigma \Vert \leq C \alpha^{-1} $ are satisfied within the
neck region.  This function can not be in $\Bonealpha$ because of the
following result.

\begin{lemma}
    \label{lemma:eigenipbd}
    The quantity $\int_{\Mapprox} \sigma \cdot S_{\alpha}$ is bounded
    below by a positive constant independent of $\alpha$.
\end{lemma}

\begin{proof}
    
Assume that the lemma is false.  Then there is a sequence $\alpha_{j}
\rightarrow 0$ and a sequence of interpolating functions $\sigma_{j}$
as above so that $\int_{\Mapproxj} \sigma_{j} \cdot S_{\alpha_{j}}
\rightarrow 0$.  Let $u_{j} = \sigma_{j} - \langle \sigma_{j} , 1
\rangle / \mathit{Vol} (\Mapproxj) - \langle \sigma_{j} ,
S_{\alpha_{j}} \rangle \, S_{\alpha_{j}}$, using the abbreviation
$\langle u, v \rangle = \int_{\Mapproxj} u \cdot v$.  Now, $\langle
u_{j}, 1 \rangle = \langle u_{j}, S_{\alpha_{j}} \rangle = 0$ so
$u_{j} \in \Bonealphaj$.  Consequently, the $\alpha_{j}$-independent
second eigenvalue estimate applies to $u_{j}$ to give
\begin{gather}
    \int_{\Mapproxj} \Vert \nabla u_{j} \Vert^{2} \geq \Cnutwo
    \int_{\Mapproxj} u_{j}^{2} \, ,\notag \\
    \intertext{which implies}
    \Vert \nabla \sigma_{j} \Vert_{L^{2}(\Mapproxj)}^{2} + \langle
    \sigma_{j} , S_{\alpha_{j}} \rangle^{2} \geq C \left( \Vert u_{j}
    \Vert_{L^{2}(\Mapproxj)}^{2} - \langle \sigma_{j} , 1 \rangle^{2}
    \right)
    \label{eqn:eigenip}
\end{gather}
using orthogonality properties of the functions $1$ and
$S_{\alpha_{j}}$ and Green's identity.  By construction of $\sigma
_{j}$, the right hand side of \eqref{eqn:eigenip} is bounded below by
a constant independent of $j$, whereas the $\langle \sigma_{j},
S_{\alpha_{j}} \rangle$ term goes to zero with $j$ by assumption and
$\Vert \nabla \sigma_{j} \Vert_{L^{2}(\Mapproxj)}^{2} \leq C \alpha
_{j}^{n-2} \rightarrow 0$ as calculated in Proposition
\ref{prop:firsteigen}.  This contradiction proves the lemma.
\end{proof}

Since the quantity $\langle \sigma , S_{\alpha} \rangle$ is positive,
one can define the approximate first eigenfunction $\Sapprox$ by
\begin{equation}
    \Sapprox = \frac{\sigma - \langle \sigma , 1 \rangle 
    \mathit{Vol}(\Mapprox)^{-1}}{ \langle \sigma, S_{\alpha} \rangle} 
    \, .
    \label{eqn:eigenapprox}
\end{equation}
Note that $\langle \Sapprox, 1 \rangle = 0$ and $\langle \Sapprox,
S_{\alpha} \rangle = 1$.  The two properties of $S_{\alpha}$ needed to
estimate the inner products \eqref{eqn:ip} will be deduced from the
$C^{0}$ difference between $S_{\alpha}$ and $\Sapprox$.

\begin{prop} Suppose $\alpha \leq \Rinj$.  Then the approximate first 
eigenfunction $\Sapprox$ satisfies the estimate
\begin{equation}
    \vert \Sapprox - S_{\alpha} \vert_{0, \Mapprox} \leq C \alpha 
    ^{(n-2)/2} \, ,
    \label{eqn:eigencoest}
\end{equation}
where $C$ is a constant independent of $\alpha$.
\end{prop}

\begin{proof} Let $x \in \Mapprox$ be a point where $\Sapprox -
S_{\alpha}$ is maximal and let $B_{R}(x)$ be any ball about $x$ with
radius independent of $\alpha$.  The De Giorgi Nash estimate with
$p=2$ and any $q >n$ implies that
\begin{equation}
    \vert \Sapprox - S_{\alpha} \vert_{0, B_{R}(x)} \leq C \left( 
    \Vert \Sapprox - S_{\alpha} \Vert_{L^{2}(\Mapprox)} + \Vert 
    \Delta_{\Mapprox} \Sapprox - \Delta_{\Mapprox} S_{\alpha} 
    \Vert_{L^{q/2}(\Mapprox)} \right) \, .
    \label{eqn:ecoest1}
\end{equation}
Estimate each term on the right hand side of \eqref{eqn:ecoest1} in 
turn.  First, by construction of $\Sapprox$, the difference 
$\Sapprox - S_{\alpha}$ is orthogonal to $1$ and $S_{\alpha}$; thus 
the min-max characterization of the second eigenvalue as well as the 
triangle inequality yields
\begin{align}
    \Vert \Sapprox - S_{\alpha} \Vert_{L^{2}(\Mapprox)} &\leq C \left(
    \Vert \nabla \Sapprox \Vert_{L^{2}(\Mapprox)} + \Vert \nabla
    S_{\alpha} \Vert_{L^{2}(\Mapprox)} \right) \notag \\
    &\leq C \left \Vert \nabla \Sapprox \Vert_{L^{2}(\Mapprox)} +
    \Vert S_{\alpha} \cdot \Delta S_{\alpha} \Vert_{L^{2}(\Mapprox)}
    \right) \tag{by Green's identity} \\
    &\leq C \alpha^{(n-2)/2} \, ,
    \label{eqn:ecoest2}
\end{align}
by the derivative bounds on $\Sapprox$ and the fact that $S_{\alpha}$
is a normalized eigenfunction of $\Delta_{\Mapprox}$ with small
eigenvalue proportional to $\alpha^{n-2}$.  Next,
\begin{align}
    \Vert \Delta_{\Mapprox} \Sapprox - \Delta_{\Mapprox} S_{\alpha}
    \Vert_{L^{q/2}(\Mapprox)} &\leq \Vert \Delta_{\Mapprox} \Sapprox
    \Vert_{L^{q/2}(\Mapprox)} + \Vert \Delta_{\Mapprox} S_{\alpha}
    \Vert _{L^{q/2}(\Mapprox)} \notag \\
    &\leq C \alpha^{-1 + nq/2} \, ,
    \label{eqn:ecoest3}
\end{align}
using the derivative bounds of $\Sapprox$ and the first eigenvalue
estimate once again.  Combining \eqref{eqn:ecoest2} and
\eqref{eqn:ecoest3} in \eqref{eqn:ecoest1} yields the desired
estimate.
\end{proof}

This more specific information about the eigenfunction $S_{\alpha}$
makes it possible now to deduce the lower bounds \eqref{eqn:ip}.  The
most important of these, and the only one which needs to be seen
explicitly, is the following one.

\begin{prop} \label{lemma:lowerbds}
    There is a number $\Rfull$ with $0 < \Rfull \leq \Rinj$ and a
    constant $\Ctrnsvsl > 0$ which is independent of $\alpha$ so that
    if $\alpha \leq \Rfull$ then the function $\psione$ satisfies the
    inequality $\vert \langle \psione, S_{\alpha} \rangle \vert \geq
    \Ctrnsvsl$.
\end{prop}

\begin{proof}

Choose $\alpha \leq \Rinj$ so that all the estimates of the previous
sections are valid.  Using integration by parts, the Neumann condition
and the definition of $v_{e}$ given in equation \eqref{eqn:bdfnext},
one deduces
\begin{align}
    -\left\langle \psione, S_{\alpha} \right\rangle &= \left\langle 
    \Delta_{\Mapprox} v_{e} , S_{\alpha} \right\rangle \notag \\
    &= \int_{\partial \Mapprox} S_{\alpha} \frac{\partial 
    v_{e}}{\partial n} - \nu_{1} \left\langle v_{e} , S_{\alpha} 
    \right\rangle \notag \\
    &= \int_{\partial \Mapprox} \negthickspace S_{\alpha} \, v \: -
    \nu_{1} \left\langle v_{e} , S_{\alpha} \right\rangle \, .
    \label{eqn:fullest4}
\end{align}
The properties of $S_{\alpha}$ that have been deduced from the 
approximation $\Sapprox$ will be used to estimate both these terms.  
First, by the estimate of $\nu_{1}$,
\begin{equation}
    \nu_{1} \left\langle v_{e}, S_{\alpha} \right\rangle \leq \Cnuone
    \alpha^{n-2}
    \label{eqn:fullest4b}
\end{equation}
where $C$ is independent of $\alpha$ because the $L^{2}$ norm of
$v_{e}$ depends only on the geometry of $\Mapprox$ near its boundary. 
Next,
\begin{equation}
    \int_{\partial \Mapprox} \negthickspace S_{\alpha} \, v \geq
    \int_{\partial \Mapprox} \negthickspace \bar{S}_{\alpha} \, v -
    \int_{\partial \Mapprox} \big\vert S_{\alpha} - \Sapprox \big\vert
    \, \vert v \vert \geq \frac{2}{\langle \sigma , S_{\alpha}
    \rangle} - \Capprox \alpha^{(n-2)/2} \, ,
    \label{eqn:fullest5}
\end{equation}
using the supremum norm estimate of $\vert \Sapprox - S_{\alpha}
\vert$ at the boundary and the fact that the boundary values of
$\Sapprox$ are known exactly.  In Lemma \ref{lemma:eigenipbd}, it was
shown that the quantity $\langle \sigma, S_{\alpha} \rangle$ is
bounded below by a constant independent of $\alpha$.  Thus combining
the upper bound \eqref{eqn:fullest4b} with the lower bound
\eqref{eqn:fullest5} in inequality \eqref{eqn:fullest4} yields the
desired result when $\alpha$ is sufficiently small.
\end{proof}

\subsection{The Full Injectivity Estimate for $\Dif F_{\alpha}(0,0,0)$}

The groundwork for the full injectivity estimate for the operator 
$\Dif F_{\alpha}(0,0,0)$ has now been laid and it is possible to prove 
the following theorem.

% Full estimate in the weighted norm

\begin{thm} \label{thm:fullest}
    If $\alpha \leq \Rsurj$, then there is a geometric constant
    $C$ independent of $\alpha$ so that the following estimate is
    valid:
    \begin{equation}
	\left( \vert u \vert_{C^{2,\beta}_{\rho}(\Mapprox)}^{2} +
	a^{2} + b^{2} \right)^{1/2} \leq C \eps^{-2\beta} \big\vert
	\rho^{2} ( P_{\alpha} u + a \psizero + b \psione ) \big\vert
	_{C^{0,\beta}_{\rho}(\Mapprox)}
        \label{eqn:linearbd}
    \end{equation} 
    for any $(u,a,b) \in \Bonealpha \times \R^{2}$.  Consequently, 
    $C_{L}(\alpha) = C \eps^{2\beta}$
\end{thm}

\begin{proof}
    
Suppose that the theorem is false.  Then there is a sequence of
parameters $\alpha_{j} \rightarrow 0$ (thus also a sequence of scales
$\eps_{j}$ and weight functions $\rho_{j}$) and a sequence of elements
$(u_{j}, a_{j}, b_{j}) \in \Bonealphaj$, normalized so that $\Vert
(u_{j}, a_{j}, b_{j}) \Vert_{\ctbrj \times \R^{2}}$ = 1 for every $j$,
which satisfy the estimates
\begin{equation}
    \eps_{j}^{-2\beta} \big\vert \rho_{j}^{2} \left( P_{\alpha_{j}}
    u_{j} + a_{j} \psizeroj + b_{j} \psionej \right) \big\vert_{\cobrj}
    \leq \frac{1}{j} \, .
    \label{eqn:approx2}
\end{equation}
It will be shown that $(u_{j}, a_{j}, b_{j}) \rightarrow (0,0,0)$,
contradicting the normalization.

The first step towards establishing the contradiction is to show that
the sequences $a_{j}$ and $b_{j}$ tend towards zero as $j \rightarrow
\infty$.  Begin by integrating both sides of equation
\eqref{eqn:approx2} over $\Mapproxj$ and estimating terms:
\begin{align}
    \int_{\Mapproxj} \frac{\eps_{j}^{2\beta}}{j} &\geq \vert a_{j}
    \vert \int_{\Mapproxj} \rho_{j}^{2} \psizeroj - \vert b_{j} \vert
    \, \left\vert \int_{\Mapproxj} \rho_{j}^{2} \psionej \right\vert -
    \left\vert \int_{\Mapproxj} \rho_{j}^{2} \cos(\theta_{\Mapproxj})
    \Delta_{\Mapproxj} u_{j} \right\vert \notag \\
    &\qquad - \left\vert \int_{\Mapproxj} \rho_{j}^{2}
    \sin(\theta_{\Mapproxj}) \langle \meancurvj \nabla u_{j} \rangle
    \right\vert \, .
    \label{eqn:approx3}
\end{align}
The first integral in equation \eqref{eqn:approx3} is bounded below
according since $\rho = R$ and $\psizeroj = 1$ in $M_{1}' \cup M_{2}'$. 
The second integral vanishes because $\rho_{j}$ is constant in the
support of $\psionej$ which occurs only near the boundary:
\begin{align}
    \int_{\Mapproxj} \rho_{j}^{2} \psionej &= - R^{2} \int_{\Mapproxj} 
    \Delta_{\Mapproxj} v_{e} \notag \\
    &= R^{2} \int_{\partial \Mapproxj} v \notag \\
    &= 0
    \label{eqn:approx3c}
\end{align}
by the definition of $v$.  

For the third integral in equation \eqref{eqn:approx3}, use the
orthogonality of $\Delta_{\Mapproxj} u_{j}$ to the constant function
$R^{2}$ to estimate
\begin{align}
    \left\vert \int_{\Mapproxj} \rho_{j}^{2} \cos(\theta_{\Mapproxj})
    \Delta_{\Mapproxj} u_{j} \right\vert &= \left\vert
    \int_{\Mapproxj} \left( \rho_{j}^{2} \cos(\theta_{\Mapproxj}) -
    R^{2} \right) \Delta_{\Mapproxj} u_{j} \right\vert \notag \\
    &\hspace{-10ex} \leq C \left\Vert \rho_{j}^{2} - R^{2}
    \right\Vert_{L^{2}(\Mapproxj)} \big\Vert \Delta_{\Mapproxj} u
    \big\Vert_{L^{2} (\Mapproxj)} \tag{using the bounds on 
    $\theta_{\Mapproxj}$} \\
    &\hspace{-10ex} \leq C \eps_j^{n\beta} \, ,
    \label{eqn:approx3b}
\end{align}
because of the fact that $\rho_{j}$ differs from $R$ only in a ball of
radius $\eps_{j}^{\beta}$.  The second eigenvalue estimate as well as
the normalization of $u_{j}$ have been used to control the $\Vert
\Delta_{\Mapproxj} u \Vert_{L^{2} (\Mapproxj)} $ term.  Finally, the
last integral in \eqref{eqn:approx3} is bounded above by $C
\alpha_{j}$ because of the bounds on $\theta_{\Mapproxj}$ from
\eqref{eqn:trigbounds} and on the mean curvature from Proposition
\ref{prop:meancurv}.  Taking this fact along with equations
\eqref{eqn:approx3b} and \eqref{eqn:approx3c} together with equation
\eqref{eqn:approx3} and rearranging terms yields an estimate of the
form
\begin{equation}
    \vert a_{j} \vert \leq C \left( \frac{\eps_{j}^{2\beta}}{j} +
    \eps_{j}^{n \beta} + \alpha_{j} \right) \, ,
    \label{eqn:agoesto0}
\end{equation}
which shows that $a_{j} \rightarrow 0$.  

Next, multiply both sides of equation \eqref{eqn:approx2} by
$S_{\alpha_{j}}$, integrate over $\Mapproxj$, and use the
orthogonality of $\Delta_{\Mapproxj} u_{j}$ to the eigenfunction
$S_{\alpha_{j}}$ as in \eqref{eqn:approx3b} to obtain
\begin{align}
    C \frac{\eps_{j}^{2 \beta} \vert S_{\alpha_{j}}
    \vert_{0,\Mapproxj}}{j} &\geq \vert b_{j} \vert \left\vert
    \int_{\Mapproxj} \rho_{j}^{2} \psionej \cdot S_{\alpha_{j}}
    \right\vert - \left\vert \int_{\Mapproxj} (\rho_{j}^{2}
    \cos(\theta_{\Mapproxj}) - R^{2} ) \big( \Delta_{\Mapproxj} u_{j}
    \big) S_{\alpha_{j}} \right\vert \notag \\
    &\qquad - \left\vert \int_{\Mapproxj} \sin(\theta_{\Mapproxj} )
    \langle \meancurvj, \nabla u_{j} \rangle \right\vert - \vert a_{j}
    \vert \left\vert \int_{\Mapproxj} \rho_{j}^{2} \psizeroj \cdot
    S_{\alpha_{j}} \right\vert \, .
    \label{eqn:bgoesto0}
\end{align}
The $\vert S_{\alpha_{j}} \vert_{0,\Mapproxj}$ term can be bounded
below by using the approximating eigenfunction $\Sapproxj$.  The first
integral on the right hand side of equation \eqref{eqn:bgoesto0} is
bounded below as a result of Lemma \ref{lemma:lowerbds} and the
properties of $\rho$.  The second and third integrals can be bounded
below exactly as in the analysis of equation \eqref{eqn:approx3b}. 
The last integral can be bounded below using the fact that
$\rho_{j}^{2} \psizeroj \approx 1$ in $M_{1}' \cup M_{2}'$ and
$\int_{\Mapproxj} S_{\alpha_{j}} = 0$.  Consequently, equation
\eqref{eqn:bgoesto0} leads to an estimate of the form
\begin{equation}
    \vert b_{j} \vert \leq C \left( \frac{\eps_{j}^{4 \beta}}{j} +
    \vert a_{j} \vert + \eps_{j}^{n\beta} + \alpha_{j} \right) \, ,
    \label{eqn:bgoneto0}
\end{equation}
where $C$ is independent of $j$, which shows that the sequence $b_{j}
\rightarrow 0$ as well.

The fact that the limits of the sequences $a_{j}$ and $b_{j}$ are zero
has two contradictory consequences.  On the one hand, the
normalization of $ (u_{j}, a_{j}, b_{j}) $ implies that if $j$ is
sufficiently large, then $\vert u_{j} \vert_{\ctbrj} \geq
\frac{1}{2}$.  On the other hand, the assumption \eqref{eqn:approx2}
in conjunction with \eqref{eqn:agoesto0} and \eqref{eqn:bgoneto0}
implies that
$$\eps_{j}^{-2\beta} \vert \rho_{j}^{2} \Delta_{\Mapproxj} u_{j}
\big\vert_{\cobrj} \longrightarrow 0 \, .$$
Thus by the injectivity estimate for $\Delta_{\Mapproxj}$ in the
$C^{k,\beta}_{\rho_{j}}$ norms,
$$\frac{1}{2} \leq \vert u_{j} \vert_{\ctbrj} \leq \Cinj \eps_{j}^{-
2\beta} \vert \rho_{j}^{2} \Delta_{\Mapproxj} u_{j} \vert_{\cobrj} 
\longrightarrow 0 \, .$$
This contradiction shows that the assumption \eqref{eqn:approx2} must
be false, thereby completing the proof of the theorem.
\end{proof}

\subsection{Surjectivity of the Linearized Operator}

The preceding result concludes the lengthy proof of the injectivity of
the linearized deformation operator $\Dif F_{\alpha} (0,0,0)$ and
identifies the dependence of the injectivity bound on the parameter
$\alpha$.  The final step in the analysis of the linearized operator
is to show surjectivity.  But given all the work that has been done so
far, this is now a fairly trivial matter.

% Surjectivity

\begin{thm} \label{thm:surj}
    If $\alpha \leq \Rsurj$, then the operator $\Dif F_{\alpha}(0,0,0)
    : \mathcal{B}_{1, \alpha} \times R^{2} \longrightarrow
    C^{0,\beta}_{\rho}(\Mapprox)$ is surjective.
\end{thm}
    
\begin{proof}
    
Recall that when $\alpha$ is sufficiently small, $\Dif
F_{\alpha}(0,0,0)$ is a small perturbation of the operator
\begin{equation}
    (u,a,b) \mapsto \Delta_{\Mapprox} u + a \psizero + b \psione \,
    .
    \label{eqn:mesurj}
\end{equation}
Thus by standard elliptic theory, the surjectivity of $\Dif
F_{\alpha}(0,0,0)$ is equivalent to the surjectivity of the operator 
in \eqref{eqn:mesurj}.

Let $f$ be any function in $C^{0,\beta}_{\rho}(\Mapprox)$.  To show 
surjectivity, one must solve the equation
\begin{equation}
    \Delta_{\Mapprox} u =  a \psizero + b \psione - f 
    \label{eqn:surj1}
\end{equation}
for $(u,a,b) \in \Bonealpha \times \R^{2}$.  By the self-adjointness
of the Laplacian, equation \eqref{eqn:surj1} can be solved for $u \in
\mathcal{B}_{1 ,\alpha}$ if and only if the right hand side of the
equation is orthogonal to both the functions $1$ and $S_{\alpha}$. 
Thus $a$ and $b$ must be chosen to satisfy the system of equations:
\begin{equation}
    \begin{gathered}
        a \int_{\Mapprox} \psizero + b \int_{\Mapprox} \psione =
        \int_{\Mapprox} f \\
        a \int_{\Mapprox} \psizero \, S_{\alpha} + b \int_{\Mapprox} 
        \psione \, S_{\alpha} = \int_{\Mapprox} f S_{\alpha}
    \end{gathered}
    \label{eqn:surj2}
\end{equation}
in order to solve \eqref{eqn:surj1}.  But recall once again that
$\int_{\Mapprox} \psione = \int_{\partial \Mapprox} v = 0$ and that
$\int_{\Mapprox} \psizero$ is strictly positive.  Thus,
$$ a = \left( \int_{\Mapprox} f \right) \left( \int_{\Mapprox}
\psizero \right)^{-1} \, .$$
This value of $a$ can now be substituted into the second of the two
equations in \eqref{eqn:surj2} above and because $\int_{\Mapprox}
\psione \, S_{\alpha}$ is also strictly positive according to
Proposition \eqref{lemma:lowerbds}, $b$ can be found as well.

The integrability conditions for equation \eqref{eqn:surj1} can thus
be satisfied by a suitable choice of $a$ and $b$.  Consequently, a
function $u \in \mathcal{B}_{1 ,\alpha}$ satisfying \eqref{eqn:surj1}
exists.  The regularity of $u$ follows from elliptic regularity theory
for the operator $\Delta_{\Mapprox}$ and the fact that $f$, $\psizero$
and $\psione$ are at least $C^{0,\beta}$.  
\end{proof}

\section{Solving the Nonlinear Problem}

\subsection{The Nonlinear Estimate}

The remaining components of the proof of the Main Theorem are to bound
the nonlinearities of $F_{\alpha}$ and to estimate the size of
$F_{\alpha}(0,0,0)$ according to the outline presented in Section 2.1. 
The difference
\begin{multline}
    \left\vert \rho^{2} \Big( \big(\Dif F_{\alpha}(H, \theta, B) - \Dif
    F_{\alpha}(0,0,0) \big)(u,a,b) \Big) \right\vert_{\cobr} \\ \leq
    C_{N} (\alpha) \, \Vert (H,\theta, B) \Vert_{\ctbr \times \R^{2}}
    \, \Vert (u, a, b) \Vert_{\ctbr \times \R^{2}}
\end{multline}
will be estimated using scaling techniques, and starts with a lemma
concerning the behaviour of the Hamiltonian flow under rescaling

\begin{lemma} \label{lemma:scale}
    Let $H: \Rtn \longrightarrow \R$ be a $C^{2,\beta}$ function and
    suppose $\phi^{t}_{H}$ denotes the Hamiltonian flow associated to
    $H$.  If $\sigma_{\eps} : \Rtn \longrightarrow \Rtn$ is the
    diffeomorphism given by $\sigma_{\eps}(x) = \eps x$ then
    $$\sigma_{\eps}^{-1} \circ \, \phi^{t}_{H} \circ \,
    \sigma_{\eps} = \phi_{H^{\eps}}^{t}$$
    where $H^{\eps}$ is the function $\eps^{-2} H \circ \,
    \sigma_{\eps}$.
\end{lemma}
 
\begin{proof}.  The proof is simply a matter of differentiating both
sides in $t$ and comparing the results.
\end{proof}

\noindent This fact about the Hamiltonian flow makes it possible to
deduce the dependence of $C_{N}(\alpha)$ on $\alpha$.

\begin{prop} \label{prop:nonlinsupest}
    Suppose $\alpha \leq \Rfull$.  Then the linearization of the
    operator $F_{\alpha}$ at the point $(H, \theta, B)$ in $\Bonealpha
    \times \R^{2}$ satisfies the estimate:
    \begin{multline}
	\Big\vert \rho^{2} \big( \Dif F_{\alpha}(H, \theta, B) - \Dif
	F_{\alpha}(0,0,0) \big) \big(u,a,b \big) \Big\vert
	_{C^{0,\beta}_{\rho} ( \Mapprox)} \\ \leq C_{N}
	\eps^{-2-2\beta} \, \Vert (H,\theta, B) \Vert_{\ctbr \times
	\R^{2}} \, \Vert (u, a, b) \Vert_{\ctbr \times \R^{2}}
	\label{eqn:nonlinsup}
    \end{multline}
    for all $(u,a,b) \in \Bonealpha \times \R^{2}$, where $C_{N}$ is a
    constant independent of $\alpha$.
\end{prop}

\begin{proof}
    
Choose $x \in \Mapprox$ and let $s = \rho(x)$.  Denote by $\sigma$ the
rescaling given by $\sigma(x) = s x$.  Denote by $\bar{h}_{1}$ the
rescaled embedding $\sigma^{-1} \circ \, \bar{h}_{\alpha}$ which
embeds the submanifold $\frac{1}{s} \Mapprox$ into $\Rtn$. 
Furthermore, let $\mathrm{Vol} _{1}$ denote the volume form of this
embedding and $\langle \cdot, \cdot \rangle_{1}$ its induced metric. 
Recall that $F_{\alpha}(H, \theta, B)$ is a $C^{0,\beta}$ function of
$\Mapprox$ and so it can be evaluated at $x$.  This leads to the
calculation:
\begin{align}
    F_{\alpha}(H, \theta, B)(x) &= \left\langle \left( 
    \phi_{v_{e}}^{B} \circ \, \phi_{H_{e}} ^{1} \circ \, \bar{h}_{1}
    \right)^{\ast} \, \sigma^{\ast} \, \mathbf{Im} \big( \me^{\mi
    \theta} \dif z \big), \mathrm{Vol}_{\Mapprox}
    \right\rangle_{\Mapprox} (x) \notag \\
    &= \left\langle \left( \left(\sigma^{-1} \circ \, \phi_{v_{e}}^{B}
    \circ \, \sigma \right) \circ \, \left(\sigma^{-1} \circ \,
    \phi_{H_{e}} ^{1} \circ \, \sigma \right) \circ \, \bar{h}_{1}
    \right)^{\ast} \, \sigma^{\ast} \, \mathbf{Im} \big( \me^{\mi
    \theta} \dif z \big), \mathrm{Vol}_{\Mapprox}
    \right\rangle_{\Mapprox} (x) \notag \\
    &= \left\langle \left( \phi_{v_{e}}^{s^{-1} B} \circ \,
    \phi_{H_{e}^{s}} ^{1} \circ \, \bar{h}_{1} \right)^{\ast} \,
    \mathbf{Im} \big( \me^{\mi \theta} \dif z \big), \mathrm{Vol}_{1}
    \right\rangle_{1}(\sigma(x)) \, .
    \label{eqn:hardest}
\end{align}
The right side of equation \eqref{eqn:hardest} can be written as $G
\big( H^{s}, \theta, s^{-1} B\big) (x)$, where $G : \Rtn \rightarrow
\R$ restricted to the $\Mapprox$.  Since, by definition of $\rho$ the
bounds on the induced metric of $\Mapprox$ are uniform in a ball of
radius $\frac{s}{2\Crho}$ about $x$, the map $G$ is independent of
$\alpha$ in this ball.  Consequently,
\begin{align}
    \left\vert \rho^{2} \left( \Dif F_{\alpha}(H, \theta, B) - \Dif
    F_{\alpha}(0,0,0) \right) \right\vert_{0, B_{s/2K_{\eps}}(x)} &
    \notag \\
    &\hspace{-35ex} \leq \left\vert \, \rho^{2} \cdot \left(\left. 
    \frac{\dif G(H^{s} + t u^{s}, \theta + t a, B/s + t b/s)}{\dif t}
    \right\vert_{t=0} - \left.  \frac{\dif G(t u^{s}, t a, tb/s)
    }{\dif t} \right\vert_{t=0} \right) \, \right\vert_{0,
    B_{s/2\Crho}(x)} & \notag \\
    &\hspace{-35ex} \leq C s^{2} \cdot \Vert ( u^{s}, a, s^{-1} b)
    \Vert_{C^{2}(B_{1/2\Crho}(x) ) \times \R^{2}} \cdot \Vert ( H^{s},
    \theta, s^{-1} B ) \Vert_{C^{2} ( B_{1/2\Crho}(x) ) \times \R^{2}}
    \, ,
    \label{eqn:nonlinsup1}
\end{align}
as a result of straightforward continuity bounds on the map $G$ in the
rescaled ball $B_{1/2\Crho}(x)$.  Simple scaling arguments and the
lower bound on $\rho$ then lead to the estimate
\begin{equation}
    \Vert ( u^{s}, a, s^{-1} b) \Vert_{C^{2}(B_{1/2\Crho}(x)) \times
    \R^{2}} \leq C \eps^{-2} \, \Vert (u, a, b) \Vert_{C^{2} (B_{s/2
    \Crho}(x)) \times \R^{2}}
    \label{eqn:nonlinsup2}
\end{equation}
for any $(u,a,b) \in \Bonealpha \times \R^{2}$.  Substituting for the
norms on the right hand side of \eqref{eqn:nonlinsup1} yields the
supremum estimate needed to prove Proposition \ref{prop:nonlinsupest}. 
The H\"older estimate follows from similar, though more involved,
calculations which can be found in \cite{me}.  The $\eps^{-2\beta}$
factor arises during the course of this latter calculation.
\end{proof}

\subsection{The Size of the $F_{\alpha}(0,0,0)$ Term}

The constants $C_{L}(\alpha)$ and $C_{N}(\alpha)$ required to estimate
the size of the neighbourhood of surjectivity of the map $F_{\alpha}$
about $(0,0,0)$ according to the Inverse Function Theorem have now
been found.  It remains to estimate the size of $F_{\alpha}(0,0,0)$ in
the $C^{0,\beta}_{\rho}$ Schauder norm.

% C^{0,\beta}_{\rho} size of \rho^{2} E = F_{\alpha}(0,0,0)

\begin{prop} \label{prop:sizeofE}
    Suppose $\alpha \leq \Rfull$.  Then $F_{\alpha}(0,0,0)$ satisfies
    the bound
    $$\vert \rho^{2} F_{\alpha}(0,0,0) \vert_{C^{0,\beta}_{\rho}
    (\Mapprox)} \leq \CsizeE \alpha^{3 - 2\beta - 2 \beta /n} \, ,$$
    where the constant $\CsizeE$ is independent of $\alpha$.
\end{prop}

\begin{proof}
    
Abbreviate $F_{\alpha}(0,0,0) = \left\langle \bar{h}_{\alpha}^{\ast}
\left( \imdz \right) , \mathrm{Vol}_{\Mapprox} \right \rangle_
{\Mapprox}$ by $E$.  Thus by Proposition \ref{prop:transition},
$$\vert E \vert_{0, \Mapprox} \leq \Cimdz \, \alpha \, .$$
But it is also true that $E$ vanishes everywhere except inside the
transition region $T_{1} \cup T_{2}$.  Property 1 of the weight
function $\rho$ implies that $\vert \rho(x) \vert \leq C \eps^{-\beta}
\Vert x \Vert$ in the transition region.  Since $\Vert x \Vert \leq C
\alpha$ there,
\begin{equation}
    \vert \rho^{2} E \vert_{0,\Mapprox} \leq C \alpha^{3} 
    \eps^{-2 \beta} \, ,
    \label{eqn:E0}
\end{equation}
where $C$ is a constant independent of $\alpha$.  Using the
relationship $\eps \leq C \alpha^{1 + 1/n}$ yields supremum estimate
required to prove the proposition.  The H\"older estimate can be
estimated in a similar, though more involved manner.  The relevant
calculations can be found in \cite{me}.
\end{proof}

\subsection{Invoking the Inverse Function Theorem}

All four components of the proof of the Main Theorem are now in place:
the injectivity of the linearized operator and the injectivity bound;
the surjectivity of the linearized operator; the estimate on the
nonlinearities of $F_{\alpha}$; and the size of $F_{\alpha}$ at zero. 
All that remains is to assemble these results and to conclude the
proof of the Main Theorem.

\medskip \noindent \scshape Existence and Regularity of $M_{\alpha}$
\upshape \medskip

Choose $\alpha \leq \Rsurj$.  The calculations of Section 5 show that
the linearized operator $\Dif F_{\alpha} (0,0,0)$ is a bijection. 
Furthermore, the functional dependence of $C_{L}(\alpha)$ and
$C_{N}(\alpha)$ on $\alpha$ have been estimated as well, in Sections 5
and 6, respectively.  Substitute these estimates, along with the
estimate the size of $F_{\alpha}(0,0,0)$ found in Section 6.1, into
the inequality \eqref{eqn:fundest} of the Inverse Function Theorem. 
One concludes that a solution of the equation $F_{\alpha}(u,a,b) = 0$
in the space $\Bonealpha \times \R^{2}$ can be found if
\begin{equation}
    C \alpha^{3 - 2\beta - 2 \beta /n} \leq \frac{\big( C_{L}
    \eps^{2\beta} \big)^{2}}{4 C_{N} \eps^{-2 - 2\beta}} \, .
    \label{eqn:magicnumber}
\end{equation}
Rewrite \eqref{eqn:magicnumber} in terms of $\alpha$ alone: thus a
solution exists whenever
$$\alpha^{1 - 2/n - 8 \beta (1+1/n)} \leq C \, .$$
By examining the power of $\alpha$ in this inequality, it is easy to
see that the inequality can be satisfied by \emph{all} sufficiently
small $\alpha$ so long as $n \geq 3$ and $\beta$ is itself chosen less
than some $\alpha$-independent upper bound.  Therefore, there always
exists $(H_{\alpha}, \theta_{\alpha}, b_{\alpha})$ solving the
equation $F_{\alpha}(H,\theta,b) =0$ provided $\alpha$ is smaller than
some upper bound.  Note that the \eqref{eqn:magicnumber} can
\emph{not} be satisfied for all sufficiently small $\alpha$ in the
case $n=2$.

The solution $(H_{\alpha}, \theta_{\alpha}, b_{\alpha})$ that has been
found belongs to $\ctbr \times \R^{2}$ and thus $H_{\alpha}$ is smooth
by standard elliptic regularity theory \cite{taylor}.  By the
definition of the map $F_{\alpha}$, this means that the submanifold
$$M_{\alpha} \equiv \phi^{1}_{b_{\alpha}} \circ \, \phi^{1}
_{(H_{\alpha}) _{e}} \big( \Mapprox \big) $$
is a smooth, minimal Lagrangian submanifold calibrated by the form
$\mathbf{Re} \left( \me^{\mi \theta_{\alpha}} \dif z \right)$ and has
boundary lying on the scaffold $W_{\alpha} = \phi^{1}_{b}(W)$. 
\medskip

\noindent \scshape Embeddedness of $M_{\alpha}$ \upshape \medskip

According to the Inverse Function Theorem, the solution of the
equation $F_{\alpha}(H_{\alpha}, \theta_{\alpha}, b_{\alpha}) = 0$
satisfies the estimate:
\begin{equation}
    \left\Vert H_{\alpha}, \theta_{\alpha}, b_{\alpha}
    \right\Vert_{\ctbr \times \R^{2}} \leq \frac{C_{L} \eps^{2\beta}}{2
    C_{N}\eps^{-2-2\beta}} = C \eps^{2 + 4\beta}
    \label{eqn:nearly}
\end{equation}
for some geometric constant $C$ independent of $\alpha$.  This, in
turn, implies that
\begin{equation}    
    \label{eqn:smallsize}
    \vert \rho \nabla H_{\alpha} \vert_{0, \Mapprox} \leq C \eps^{2 +
    4\beta} \qquad \mbox{and} \qquad \vert b_{\alpha} \vert \leq C
    \eps^{2+ 4 \beta} \, .
\end{equation}
Both these quantities tend towards zero as $\alpha \rightarrow 0$.

In order to determine whether $M_{\alpha} = \phi_{b_{\alpha}}^{1}
\circ \, \phi_{(H_{\alpha})_{e}}^{1} \big( \Mapprox \big)$ is an
embedded submanifold of $\Rtn$, it is enough to show that the
Hamiltonian diffeomorphism $\phi_{b_{\alpha}}^{1} \circ \,
\phi_{(H_{\alpha})_{e}}^{1}$ deforms $\Mapprox$ in a sufficiently
small manner in the $C^{1}$ sense.  But this is now a consequence of
the fact that the Hamiltonian flow is a smooth function of the
Hamiltonian and its first derivatives.  The proof of the Main Theorem
is now complete.  \hfill \qedsymbol

\bigskip

\noindent \scshape Acknowledgements: \upshape I would like to thank my
Ph.D. advisor at Stanford University, Rick Schoen, for his patience,
insight and confidence in me while I was carrying out the research for
this paper.  I would also like to thank Rafe Mazzeo and Leon Simon for
their generous assistance, as well as Justin Corvino and Vin de Silva
for their inspirations, ideas, and careful proofreading.  Finally, I
would like to thank Yng-Ing Lee for discovering and helping me to
overcome certain mathematical difficulties.  

\newpage

\renewcommand{\baselinestretch}{1.0}
\normalsize

\bibliography{paper}

\providecommand{\bysame}{\leavevmode\hbox to3em{\hrulefill}\thinspace}
\begin{thebibliography}{10}

\bibitem{amr}
R.~Abraham, J.~E. Marsden, and T.~Ratiu, \emph{Manifolds, {T}ensor {A}nalysis,
  and {A}pplications}, second ed., Springer-Verlag, New York, 1988.

\bibitem{me2}
Adrian Butscher, \emph{Deformations of minimal {L}agrangian submanifolds with
  boundary}, Submitted to Proceedings of the AMS, August 2001. Preprint
  math.DG/0110052.

\bibitem{me}
\bysame, \emph{Deformation {T}heory of {M}inimal {L}agrangian {S}ubmanifolds},
  Ph.D. thesis, Stanford University, 2000.

\bibitem{gt}
David Gilbarg and Neil~S. Trudinger, \emph{Elliptic {P}artial {D}ifferential
  {E}quations of {S}econd {O}rder}, second ed., Springer-Verlag, Berlin, 1983.

\bibitem{harvey}
Reese Harvey, \emph{Spinors and {C}alibrations}, Academic Press, 1990.

\bibitem{hl1}
Reese Harvey and H.~Blaine Lawson, Jr., \emph{Calibrated geometries}, Acta
  Math. \textbf{148} (1982), 47--157.

\bibitem{haskins}
Mark Haskins, \emph{Constructing {S}pecial {L}agrangian {C}ones}, Ph.D. thesis,
  University of Texas at Austin, 2000.

\bibitem{hitchin}
Nigel~J. Hitchin, \emph{The moduli space of special {L}agrangian submanifolds},
  Ann. Scuola Norm. Sup. Pisa Cl. Sci. (4) \textbf{25} (1997), no.~3-4,
  503--515 (1998), Dedicated to Ennio De Giorgi.

\bibitem{joyce1}
Dominic Joyce, \emph{Lectures on {C}alabi-{Y}au and special {L}agrangian
  geometry}, math.DG/0108088.

\bibitem{joyce}
\bysame, \emph{Compact {R}iemannian $7$-manifolds with holonomy ${G}\sb 2$.
  {I}, {I}{I}}, J. Differential Geom. \textbf{43} (1996), no.~2, 291--328,
  329--375.

\bibitem{nikos2}
Nikolaos Kapouleas, \emph{Constant mean curvature surfaces in {E}uclidean
  spaces}, Proceedings of the International Congress of Mathematicians, Vol.\
  1, 2 (Z\"urich, 1994) (Basel), Birkh\"auser, 1995, pp.~481--490.

\bibitem{nikos1}
\bysame, \emph{On desingularizing the intersections of minimal surfaces},
  Proceedings of the 4th International Congress of Geometry (Thessaloniki,
  1996), Giachoudis-Giapoulis, Thessaloniki, 1997, pp.~34--41.

\bibitem{lawlor2}
Gary Lawlor, \emph{The angle criterion}, Invent. Math. \textbf{95} (1989),
  437--446.

\bibitem{ynging}
Yng-Ing Lee, \emph{Embedded special {L}agrangians in calabi-yau manifolds},
  Preprint, 2001.

\bibitem{rafe}
Rafe Mazzeo, Frank Pacard, and Dan Pollack, \emph{{C}{M}{C} {S}urfaces}, To
  appear.

\bibitem{mclean}
Robert~C. McLean, \emph{Deformations of calibrated submanifolds}, Comm. Anal.
  Geom. \textbf{6} (1998), no.~4, 705--747.

\bibitem{morrison}
David~R. Morrison, \emph{Mathematical aspects of mirror symmetry}, Complex
  Algebraic Geometry (Park City, UT, 1993), Amer. Math. Soc., Providence, RI,
  1997, pp.~265--327.

\bibitem{sandw}
R.~Schoen and J.~Wolfson, \emph{Minimizing volume among {L}agrangian
  submanifolds}, Differential {E}quations: {L}a {P}ietra 1996 (Shatah Giaquinta
  and Varadhan, eds.), Proc. of Symp. in Pure Math., vol.~65, 1999,
  pp.~181--199.

\bibitem{schoen}
Richard Schoen, \emph{Lecture {N}otes in {G}eometric {P}{D}{E}s on
  {M}anifolds}, Course given in the Spring of 1998 at Stanford University.

\bibitem{simon2}
Leon Simon, \emph{Lecture {N}otes in {P}{D}{E} {T}heory}, Course given in 1997
  at Stanford University.

\bibitem{simon}
\bysame, \emph{Lectures on {G}eometric {M}easure {T}heory}, Australian National
  University Centre for Mathematical Analysis, Canberra, 1983.

\bibitem{taylor}
Michael~E. Taylor, \emph{Partial {D}ifferential {E}quations. {I}{I}{I}},
  Springer-Verlag, New York, 1997, Nonlinear equations.

\end{thebibliography}
\bibliographystyle{amsplain}

\end{document}